\documentclass[11pt]{amsart}
\usepackage{amsmath,amssymb,amsthm,amscd}
\usepackage[all]{xy}
\usepackage{color}
\CompileMatrices
\numberwithin{equation}{section}
\newtheorem{prop}{Proposition}[section]
\newtheorem{theo}[prop]{Theorem}
\newtheorem{lemm}[prop]{Lemma}
\newtheorem{coro}[prop]{Corollary}
\newtheorem{rema}[prop]{Remark}
\newtheorem{exam}[prop]{Example}
\newtheorem{defi}[prop]{Definition}

\newdimen\tableauside\tableauside=1.0ex
\newdimen\tableaurule\tableaurule=0.4pt
\newdimen\tableaustep
\def\phantomhrule#1{\hbox{\vbox to0pt{\hrule height\tableaurule width#1\vss}}}
\def\phantomvrule#1{\vbox{\hbox to0pt{\vrule width\tableaurule height#1\hss}}}
\def\sqr{\vbox{%
  \phantomhrule\tableaustep
  \hbox{\phantomvrule\tableaustep\kern\tableaustep\phantomvrule\tableaustep}%
  \hbox{\vbox{\phantomhrule\tableauside}\kern-\tableaurule}}}
\def\squares#1{\hbox{\count0=#1\noindent\loop\sqr
  \advance\count0 by-1 \ifnum\count0>0\repeat}}
\def\tableau#1{\vcenter{\offinterlineskip
  \tableaustep=\tableauside\advance\tableaustep by-\tableaurule
  \kern\normallineskip\hbox
    {\kern\normallineskip\vbox
      {\gettableau#1 0 }%
     \kern\normallineskip\kern\tableaurule}%
  \kern\normallineskip\kern\tableaurule}}
\def\gettableau#1 {\ifnum#1=0\let\next=\null\else
  \squares{#1}\let\next=\gettableau\fi\next}
\numberwithin{equation}{section}

\newcommand{\be}{\begin{equation}}
\newcommand{\ee}{\end{equation}}

\newcommand{\IP}{\mathbb{P}}

\newcommand\IZ{\mathbb {Z}}

\newcommand{\IC}{\mathbb{C}}
\newcommand{\IE}{\mathbb{E}}
\newcommand{\IR}{\mathbb{R}}
\newcommand{\CU}{{\mathcal U}}
\newcommand{\ba}{\begin{array}}
\newcommand{\ea}{\end{array}}

\newcommand{\CX}{{\mathcal X}}
\newcommand{\IF}{{\mathbb F}}

\newcommand{\CB}{{\mathcal B}}
\newcommand{\CK}{{\mathcal K}}
\newcommand{\CQ}{{\mathcal Q}}

\newcommand{\IH}{{\mathbb H}} 
\newcommand{\bal}{\begin{aligned}}
\newcommand{\eal}{\end{aligned}}

\newcommand{\CZ}{{\mathcal Z}}

\newcommand{\wPhi}{{\widetilde \Phi}}
\newcommand{\wphi}{{\widetilde \phi}}
\newcommand{\wpsi}{{\widetilde \psi}}

\newcommand{\mfm}{{\mathfrak{M}}}

\newcommand{\ho}{{\mathrm{Hom}}}
\newcommand{\longto}{\longrightarrow}
\newcommand{\lochom}{{\mathcal Hom}}
\newcommand{\ch}{{\mathrm{ch}}}

\newcommand{\CO}{{\mathcal O}}
\newcommand{\cE}{{\mathcal E}}
\newcommand{\CE}{{\mathcal E}}

\newcommand{\CH}{{\mathcal H}}

\newcommand{\CJ}{{\mathcal J}}

\newcommand{\CC}{{\mathcal C}}
\newcommand{\CI}{{\mathcal I}}
\newcommand{\CL}{{\mathcal L}}
\newcommand{\wCC}{{\widetilde {\mathcal C}}}

\newcommand{\mfe}{{\mathfrak E}}
\title{Moduli of ADHM Sheaves and Local Donaldson-Thomas Theory}
\author{D.-E. Diaconescu}
\address{NHETC, Rutgers University, 126 Frelinghuysen Rd, Piscataway NJ 08854}
\dedicatory{Dedicated to The Memory of John Ionesei -- John Prin\c{t}ul}

\begin{document}

\begin{abstract} 
The ADHM construction establishes a one-to-one correspondence 
between framed torsion free sheaves on the 
projective plane and stable framed  
representations of a quiver with relations 
in the category of complex vector spaces. 
This paper studies the geometry of moduli spaces of representations 
of the same quiver with relations in the abelian category of 
coherent sheaves on a smooth complex projective curve $X$. In particular it is proven  
that this moduli space is virtually smooth and related by
relative Beilinson spectral sequence to the curve  
counting construction via stable pairs of Pandharipande and Thomas. 
This yields a new conjectural construction 
for the local Donaldson-Thomas theory of 
curves as well as a natural higher rank generalization. 
\end{abstract} 
\maketitle

\tableofcontents

\section{Introduction}\label{intro}

An ADHM quiver $Q_{ADHM}$ is a quiver of the form   
\be\label{eq:adhmquiverA}
\xymatrix{
e\ \bullet \ar@(ur,ul)|{a_1}\ar@(dl,dr)|{a_2} \ar@/^/[rr]|{b}
&& \bullet\ e_\infty \ar@/^/[ll]|{c}}
\ee
with one relation specified by the linear combination of paths 
\be\label{eq:adhmquiverB}
a_1a_2-a_2a_1 + cb.
\ee
It is well known that moduli spaces of cyclic representations of the ADHM quiver 
in the category of complex vector spaces are isomorphic to moduli spaces 
of framed torsion free sheaves on the projective plane \cite[Ch. 2]{hilblect}. 
Similarly, moduli spaces of cyclic modules over a deformation 
of the ADHM quiver algebra \cite{Nekrasov:1998ss} are isomorphic to moduli spaces of 
 locally free sheaves on the noncommutative  projective 
plane  \cite{Kapustin:2000ek}. 

ADHM sheaves are framed twisted representations 
of an ADHM quiver in the abelian category of coherent 
$\CO_X$-modules of a smooth projective variety $X$ over $\IC$.
The goal of the present paper is to study the geometry 
of the moduli space of ADHM sheaves on smooth projective 
complex varieties. Note that quiver sheaves have been previously considered 
in the literature in 
\cite{HKquivers,dimred,homquiv,szendroi-2005,benzvi-2006}.
In particular ADHM sheaves have been previously studied in 
\cite{szendroi-2005} in relation to the relative Beilinson monad 
for noncommutative threefolds over curves. A similar relation
between quiver sheaves and relative Beilinson monads has been 
employed in the context of integrable systems in \cite{benzvi-2006}. 
Moduli problems for such objects have been also previously 
considered in \cite{decorated,tensors,modquivers,holchains},
as well as \cite{benzvi-2006}.
We will later elaborate on the relation between some of the 
above papers and the present work.
Before presenting the main results note that this problem is 
motivated by the following questions  in string theory and quantum field theory. 

Moduli spaces of $p-(p+4)$ D-brane configurations have
been intensively studied in the string theory literature starting with the 
work of Douglas and Moore \cite{DM-quivers}. In particular the moduli 
space of supersymmetric flat directions of such configurations 
have been shown in \cite{DM-quivers}
to be isomorphic to the quiver varieties constructed in \cite{ALEquiver-I,ALEquiver-II}. 
These are moduli spaces of supersymmetric constant field configurations on 
the D$p$-brane world-volume. In this paper we consider moduli spaces of 
topologically nontrivial supersymmetric field configurations for $p-(p+4)$ D-brane
systems when the D$p$-branes are wrapped on a complex projective variety $X$. 
Given the physical construction of Donaldson-Thomas theory \cite{INOV}, 
such moduli spaces are expected to provide a mathematical framework for counting 
problems of D6-D2-D0 supersymmetric
bound states on noncompact Calabi-Yau threefolds. Concrete applications 
to local  Gopakumar-Vafa invariants and wallcrossing 
are presented in \cite{chamberI}. 

The same moduli spaces occur in a different physical context, namely 
renormalization group flow in (topologically twisted) 
two dimensional nonabelian gauged linear sigma models. 
The gauge theory in question is the effective field theory on a system of D$1$-D$5$ 
euclidean branes, the lower dimensional branes being wrapped on a smooth 
projective curve $X$. For generic values of the Fayet-Iliopoulos parameters, 
this theory exhibits a Higgs branch isomorphic to the moduli space of 
framed torsion free sheaves on the projective plane. 
The moduli spaces constructed in this paper can be identified with moduli spaces of 
UV gauge theory instantons, or, equivalently, topologically nontrivial $Q$-fixed 
field configurations modulo gauge transformations. In this context, they provide a 
partial compactification of the moduli space of holomorphic maps to the 
moduli space of framed torsion free sheaves on the projective plane. 
This is analogous to the similar statement concerning $Quot$ schemes 
parameterizing coherent quotients on $\IP^1$ and holomorphic maps 
to the Grassmannian \cite{BDW,BCK-I}. 
Note that the later have been shown to occur in a similar 
context in \cite{Donagi:2007hi,HT}. A renormalization 
group flow analogous to \cite{MP,HV} yields a conjectural relation between 
the residual ADHM invariants defined in this paper 
and the $J$-function of the moduli space of 
framed torsion-free sheaves generalizing the results of \cite{fixed-quot,hyperquot,BCK-I,BCK-II}
for $J$-functions of Grassmannians and flag varieties. This is currently under 
investigation \cite{nonabhilb}. The 
dynamics of nonabelian gauged linear sigma models 
has been studied recently in \cite{Nekrasov:2009ui, Nekrasov:2009uh}, 
emphasizing the relation between quantum cohomology rings and 
integrable systems. Our approach is complementary since it emphasizes the 
relation between  UV instanton series and 
$J$-functions.

Let us summarize the main results and plan of the paper.
Let $X$ be a smooth projective variety over $\IC$ 
equipped with a very ample line bundle $\CO_X(1)$. 
Let $M_1,M_2$ be fixed  invertible $\CO_X$-modules, called twisting data, 
and $E_\infty$ a fixed locally free 
$\CO_X$-module called framing data. 
In order to keep the notation short, we will denote by 
$\CX$ the data $(X,M_1,M_2,E_\infty)$. ADHM sheaves are twisted representations 
of the quiver \eqref{eq:adhmquiverA} in the abelian category of 
coherent sheaves on $X$. 
(See (\ref{adhmsheaf}) for a precise definition.)
We will consider ADHM sheaves 
subject to a stability condition (\ref{adhmstab}).
A routine argument shows that flat families of 
stable ADHM sheaves with fixed Hilbert polynomial 
$P$ on $X$ form a fibered groupoid 
$\mfm_{ADHM}(\CX,P)$ 
over the category ${\mathfrak S}$ of 
schemes of finite type over $\IC$. Moreover, according to Lemma (\ref{automlemma}), 
the closed points of $\mfm_{ADHM}(\CX,P)$ have trivial stabilizers. 
Then following theorem follows directly from the results of Schmitt  
\cite[Thm. 2.9.2.44, 2.9.2.45]{GIT-decorated}.
\begin{theo}\label{DMstack} 
The groupoid $\mfm_{ADHM}(\CX,P)$ 
is a quasi-projective scheme over $\IC$.  
\end{theo}

\begin{rema}\label{GITstuff} 
$(i)$ Note that the stability condition $(ii)$ in (\ref{adhmstab})
is very similar to the nondegeneracy condition 
formulated in \cite[Thm 3.1]{szendroi-2005}. However the 
later disallows all proper subsheaves $E'\subset E$ 
satisfying the conditions listed there, not only the saturated 
ones. It will become clear in 
section (\ref{adhmadmsect}) that this difference 
has important consequences for 
the connection to local Donaldson-Thomas theory.

$(ii)$ 
In the proof of theorem \cite[Thm. 2.9.2.45]{GIT-decorated} 
the moduli scheme $\mfm_{ADHM}(\CX,P)$ is constructed 
using geometric invariant theory \cite{GIT}. 
In particular it is equipped with a natural polarization $[{\mathfrak L}]$ 
according to \cite[Thm 1.10]{GIT}, \cite[Rem. 1.4.3.9]{GIT-decorated}. 

$(iii)$ Moreover
\cite[Cor. 2.9.2.42]{GIT-decorated} proves that the stability condition (\ref{adhmstab}) is in fact the asymptotic 
form of a more general GIT stability condition for ADHM sheaves. 
Analogous stability conditions for quiver sheaves were previously formulated and studied  
in \cite{decorated,tensors,modquivers}, and in the context of 
Hitchin-Kobayashi correspondence in \cite{HKquivers,dimred}. This is a natural 
generalization of previous work on decorated sheaves, including \cite{pairsI,verlinde, surfhiggs, triples,  logflips,
stpairs, framed, frHitchin}. 

$(iii)$ Moduli spaces of twisted ADHM sheaves have also been recently studied 
by Marcos Jardim in \cite{instsheaves}, subject to a more general stability 
condition. In particular \cite[Thm 24]{instsheaves} establishes that the 
moduli space of stable twisted ADHM sheaves on a complex projective space is 
a quasi-projective variety. 

\end{rema}

Next, suppose $X$ is a smooth projective curve of genus 
$g$ over $\IC$. In this case the Hilbert polynomial 
of an ADHM sheaf $\CE$ is determined by the rank $r\in \IZ,\ r \geq 1$ and 
the degree $e\in \IZ$. We will denote by $\mfm_{ADHM}(\CX,r,e)$ 
the moduli space of stable ADHM sheaves on $X$ with fixed $(r,e)$.
Our first result establishes equivariant virtual smoothness 
of the moduli space of ADHM sheaves with respect to 
algebraic torus actions preserving the polarization $[{\mathfrak L}]$ obtained 
from the GIT construction. 
\begin{theo}\label{virsmooththm} 
Suppose $X$ is a smooth projective curve over $\IC$ 
and suppose there exists a torus action 
${\bf T}\times \mfm_{ADHM}(\CX,r,e) \to \mfm_{ADHM}(\CX,r,e)$
such that the polarization $[{\mathfrak L}]$ is {\bf T}-equivariant. 
Then  the moduli space $\mfm_{ADHM}(\CX,r,e)$ has a 
${\bf T}$-equivariant perfect tangent-obstruction theory. 
\end{theo} 

\noindent
Here we employ Definition \cite[Def. 2.1]{degGW} 
for a perfect tangent-obstruction theory of 
a Deligne-Mumford stack. 
Theorem (\ref{virsmooththm}) is proven in section (\ref{virsmooth}).

Note that there are natural algebraic torus actions 
on $\mfm_{ADHM}(\CX,P)$ satisfying the conditions 
of Theorem (\ref{virsmooththm}) presented in examples 
(\ref{torusactA}), (\ref{torusactB}). 
Proposition (\ref{fixedloci}) 
shows that if $X$ is a smooth projective curve, 
the fixed locus $\mfm(\CX,r,e)^{\bf T}$ of the torus 
action defined in Example (\ref{torusactB}) is a projective 
scheme over $\IC$. 
Then Theorem (\ref{virsmooththm}),
Proposition (\ref{fixedloci}) and \cite{GP}
yield the following corollary
\begin{coro}\label{vircycle}
Let $X$ be a smooth projective curve over $\IC$ and 
$E_\infty = \oplus_{a=1}^{r_\infty}L_a$ where $r_\infty\in \IZ$, 
$r_\infty\geq 1$ and $L_a$, $a=1,\ldots, r_\infty$, are 
line bundles on $X$. Let ${\bf T}\times \mfm_{ADHM}(\CX,r,e) \to 
\mfm_{ADHM}(\CX,r,e)$ be the torus action defined in Example 
(\ref{torusactB}).  Then the fixed locus  
$\mfm_{ADHM}(\CX,r,e)^{\bf T}$ is a 
projective scheme over $\IC$ equipped with a 
virtual cycle $[\mfm_{ADHM}(\CX,r,e)^{\bf T}]$ 
and a {\bf T}-equivariant virtual normal bundle
$N^{vir}_{\mfm_{ADHM}(\CX,r,e)^{\bf T}/\mfm_{ADHM}(\CX,r,e)}$
defined as an element of the K-theory of  
coherent locally free sheaves on $ \mfm_{ADHM}(\CX,r,e)^{\bf T}$.
\end{coro} 

Corollary (\ref{vircycle}) allows us to define an equivariant 
ADHM theory of curves by analogy with \cite{BP,OP}. 
\begin{defi}\label{localtheory}
Suppose $X$ is a smooth projective curve over $\IC$ and 
$E_\infty = \oplus_{a=1}^{r_\infty}L_a$ where $r_\infty\in \IZ$, 
$r_\infty\geq 1$ and $L_a$, $a=1,\ldots, r_\infty$, are 
line bundles on $X$. Let ${\bf T}\times \mfm_{ADHM}(\CX,r,e) \to 
\mfm_{ADHM}(\CX,r,e)$ be the torus action defined in Example 
(\ref{torusactB}).
Then the rank $r$ equivariant ADHM theory of the 
data $\CX=(X,M_1,M_2,E_\infty)$ is defined by 
\be\label{eq:localADHMtheory} 
\begin{aligned}
Z_{ADHM}& (\CX)_r(q) = \sum_{e\in \IZ} 
q^e \int_{[\mfm_{ADHM}(\CX,r,e)]^{\bf T}} 1
\end{aligned}
\ee
where $[\mfm_{ADHM}(\CX,r,e)]^{\bf T}$ is the equivariant 
virtual cycle of the moduli space. 
\end{defi}

\begin{rema}
$(i)$ Note that the partition function \eqref{eq:localADHMtheory} 
is a natural generalization of the instanton counting 
function defined in \cite{instcountA} to the relative setting.
One can also generalize 
the K-theoretic partition function of \cite{instcountB} 
to the relative setting relying on Corollary (\ref{vircycle}) 
and the results \cite{vf-dg,FG-virtual}.  

$(ii)$ Applying the virtual localization theorem \cite{GP}, the right hand side of equation 
\eqref{eq:localADHMtheory} can be written as 
\[
\bal 
& \sum_{e\in \IZ} 
q^e \int_{[\mfm_{ADHM}(\CX,r,e)^{\bf T}]} {1\over e_{\bf T} 
(N^{vir}_{\mfm_{ADHM}(\CX,r,e)^{\bf T}/\mfm_{ADHM}(\CX,r,e)})},\\
\end{aligned}
\]
where $[\mfm_{ADHM}(\CX,r,e)^{\bf T}]$ is the induced virtual cycle of the fixed locus, and 
 $e_{\bf T} 
(N^{vir}_{\mfm_{ADHM}(\CX,r,e)^{\bf T}/\mfm_{ADHM}(\CX,r,e)})$ 
is the {\bf T}-equivariant Euler class of the virtual normal bundle 
$N^{vir}_{\mfm_{ADHM}(\CX,r,e)^{\bf T}/\mfm_{ADHM}(\CX,r,e)}$.\end{rema}

As explained in the beginning of the 
introduction, string theoretic arguments predict that the ADHM theory of curves 
should be related to the local Donaldson-Thomas theory of curves
defined in \cite{MNOP-I,MNOP-II,BP,OP}.
In sections (\ref{admpairs}), (\ref{adhmadmsect}) 
we show that the relative Beilinson 
spectral sequence \cite{orlov} identifies stable 
ADHM sheaves with $E_\infty=\CO_X$ on a smooth projective curve $X$ with
stable pairs \cite{LePoitierA,LePoitierB,stabpairs-I} on the 
projective plane bundle $Y=\mathrm{Proj}(\CO_X\oplus M_1\oplus M_2)$ 
over $X$.  
The later have been recently employed by Pandharipande 
and Thomas \cite{stabpairs-I} as a new curve counting device
on smooth projective threefolds (see \cite{stabpairs-II,stabpairs-III} 
for further results.) and are conjecturally related to Donaldson-Thomas 
theory. 

The identification  between ADHM sheaves with trivial framing and 
stable pairs is most conveniently formulated in terms of a new class 
of objects -- admissible pairs -- introduced in section (\ref{admpairs}). 
Admissible pairs are characterized by a pair $(d,n)\in \IZ_{\geq 1}\times \IZ$ of 
numerical invariants defined in (\ref{admissible}) and satisfy the following 
properties  derived from 
\cite{LePoitierA,LePoitierB,stabpairs-I}.
\begin{theo}\label{admstack} 
$(i)$
There exists a quasi-projective moduli scheme $\mfm_{Adm}(Y,d,n)$
over $\IC$
of admissible pairs on $Y$ with fixed numerical invariants $(d,n)$. 

$(ii)$ There is a ${\bf T}=\IC^\times \times \IC^\times $
action on $\mfm_{Adm}(Y,d,n)$ induced by a scaling action on $Y$,
such that $\mfm_{Adm}(Y,d,n)$ is 
equipped with a natural 
${\bf T}$-equivariant perfect tangent-obstruction theory.

$(iii)$  The fixed locus $\mfm_{Adm}(Y,d,n)^{\bf T}$ is 
a projective scheme over $\IC$ equipped with a 
virtual cycle $[\mfm_{Adm}(Y,d,n)^{\bf T}]$ and an equivariant 
virtual normal bundle $N^{vir}_{\mfm_{Adm}(Y,d,n)^{\bf T}/\mfm_{Adm}(Y,d,n)}$.
\end{theo}

\noindent
A proof of this theorem based on \cite{LePoitierA,LePoitierB,stabpairs-I}
is outlined in section (\ref{admpairs}). 

Theorem (\ref{admstack}) allows us to define the local equivariant 
Pandharipande-Thomas \cite{stabpairs-I} theory of curves. 
\begin{defi}\label{localPT} 
Let $X$ be a smooth projective curve of genus $g$, $M_1, M_2$ 
be invertible sheaves on $X$ and $Y = \mathrm{Proj}(\CO_X\oplus 
M_1\oplus M_2)$. Then we define the degree $d$ 
local Pandharipande-Thomas 
theory of the triple  $(X,M_1,M_2)$ to be 
\[
Z_{PT}(X,M_1,M_2)_d(q) =\sum_{n\in \IZ} q^n \int_{[\mfm_{Adm}(Y,d,n)]^{\bf T}} 1.
\]
\end{defi}

Our third result is the following.
\begin{theo}\label{correspB}
Let $X$ be a smooth projective curve of genus $g$, $M_1, M_2$ 
be invertible sheaves on $X$ and $E_\infty=\CO_X$.  
Let $Y = \mathrm{Proj}(\CO_X\oplus M_1\oplus M_2)$
and let 
${\bf T}=\IC^\times\times\IC^\times$ 
act on $\mfm_{ADHM}(\CX,r,e)$ as in Example (\ref{torusactA}). 
Then there is a 
${\bf T}$-equivariant isomorphism of  moduli 
spaces
\[
{\mathfrak f}:\mfm_{Adm}(Y,d,n)\simeq 
\mfm_{ADHM}(\CX,d,n+d(g-1))
\]
such that the perfect tangent-obstruction theories of 
the two moduli spaces are compatible with respect to ${\mathfrak f}$, 
the relative perfect obstruction theory being trivial.   
\end{theo}

\noindent Here compatible perfect obstruction theories is meant in 
the sense of \cite[Def. 4.1]{degGW}. Theorem (\ref{correspB}) 
is proven in sections (\ref{isomsection}), (\ref{compattangobs}). 

Theorem (\ref{correspB})  implies
\begin{coro}\label{correspE} 
Let $X$ be a smooth projective curve genus $g$, $M_1, M_2$ 
be invertible sheaves on $X$ and $Y = \mathrm{Proj}(\CO_X\oplus 
M_1\oplus M_2)$. Then we have 
\be\label{eq:equivlocalth}
Z_{PT}(X, M_1, M_2)_{d}(q) = q^{-d(g-1)}Z_{ADHM}(\CX)_d(q)
\ee
for any $d\geq 1$, where $\CX=(X,M_1,M_2,\CO_X)$.
\end{coro}

According to \cite[Conj. 3.3]{stabpairs-I} we have the following conjectural 
relation between the Pandharipande-Thomas and the local Donaldson-Thomas 
theory of curves  \cite{MNOP-I,MNOP-II,BP,OP}
\[
Z_{PT}(X,M_1,M_2)_d(q) = Z'_{DT}(X,M_1,M_2)_d(q).
\]
Therefore Corollary (\ref{correspE}) immediately yields a similar 
conjectural relation between the ADHM theory of curves with trivial 
framing and the local 
Donaldson-Thomas theory of curves. 

For more general framing,
$E_\infty=\oplus_{a=1}^{r_\infty}L_a$ in Definition 
(\ref{localtheory}), the ADHM theory of 
$Z_{ADHM}(\CX)_r(q)$ represents a rigorous construction of instanton 
moduli spaces in the nonabelian gauge theory of D6-D2-D0 brane systems 
on local Calabi-Yau threefolds as above. Conjecturally it should be 
related to a  higher rank version of local Donaldson-Thomas theory
of curves provided the later were defined. 
This construction provides a mathematical framework for some of the results 
obtained in the string theory literature \cite{freefermions, Cirafici:2008sn}.

{\bf Notation and Conventions}. 
Throughout this paper, we will denote by ${\mathfrak S}$ 
the category of schemes of finite type over 
$\IC$. For any such schemes $X,S$  
we  set $X_S = S\times X$ and $X_s=k(s)\times_S X$ 
for any point 
$s\in S$, where $k(s)$ is the residual field of $s$. Let also 
$p_S: X_S \to S$, $p_X:X_S \to X$ denote the canonical projections. 
We will also set $(F)_S = p_X^*F$ for 
any $\CO_X$-module $F$. Given a morphism $f:S'\to S$, 
we will denote by $f_X = f\times 1_X:X_{S'}\to X_S$. 
Any 
morphism $f:S'\to S$, yields a  commutative diagram of the form 
\[
\xymatrix{
X_{S'} \ar[r]^-{p_X'} \ar[d]_-{f_X} & X \ar[d]^-{1_X} \\
X_S \ar[r]^-{p_X} & X. \\
}
\]
Then for any $\CO_X$-module $F$ there is a  
canonical isomorphism 
$F_{S'}\simeq f_X^* F_S$
which will be implicit in the following.
Finally, we will employ the conventions of \cite[Ch. 1.2, Ch. 1.3]{duality}
for computations in derived categories. 

{\bf  Acknowledgments.}
The physical idea of this construction originated in discussions 
with Robbert Dijkgraaf in spring 2006 concerning the 
physics of Donaldson-Thomas theory. It was reinforced in independent 
discussions with Nick Halmagy in summer 2006 
concerning different aspects of D-brane BPS states. 
I am very grateful to them for sharing their insights with me. 

During the completion of this 
project I greatly benefited from discussions and correspondence 
with Richard Thomas
who kindly made a preliminary version of \cite{stabpairs-I} available 
to me prior to publication. The results of sections 
(\ref{admpairs}), (\ref{adhmadmsect}) were 
prompted by his comments on an earlier version of the 
present  paper. The proof of Lemma (\ref{perfect})
is also due to him,  eventual errors being my own. 

I am very grateful to Tony Pantev for patient explanations 
about moduli stacks and virtual cycles, as well as 
invaluable help in proving Theorem 
(\ref{DMstack}), and to Brent Doran for illuminating discussions 
on quotients by algebraic group actions. 

I am also very grateful to Ionut Ciocan-Fontanine, Bumsig Kim and Davesh Maulik 
for very helpful discussions and collaboration on a related project. 
In particular, I owe special thanks to Ionut Ciocan-Fontanine for pointing out 
references \cite{Keel-Mori, vf-dg, FG-virtual} to me.

I would also like to thank Aaron Bertram, 
Giulio Bonelli, Lev Borisov, 
Ugo Bruzzo, Ron Donagi, Bogdan Florea, Alberto Garcia-Raboso,
Sheldon Katz, Alexander Kuznetzov, 
Y.P. Lee, Sergio Lukic, Greg Moore, Alessandro Tanzini, 
and especially Cumrun Vafa
for very helpful conversations and encouragement.  

Part of this work was done during a visit at SISSA Trieste in 
spring 2007,  
as well as during the ``Derived Categories in Mathematics and Physics - 
Snowbird, June 17-23, 2007'' workshop, the 
``Workshop on Algebraic Geometry and Physics - Seoul, June 25-29, 2007'',
the ``String Theory and Quantum Geometry - Aspen 2007'' workshop 
and the ``Simons Workshop on Mathematics and Physics, Stony Brook 2007''
which provided a very stimulating mathematical atmosphere. 

I would also like to acknowledge the partial support 
of NSF grant PHY-0555374-2006.

\section{ADHM Sheaves}\label{adhmsection}

\subsection{Stable ADHM Sheaves}

Quiver sheaves are representations of a quiver with relations 
in the abelian category of modules over a scheme or more generally 
a ringed space \cite{HKquivers,dimred,homquiv}. 
Let $X$ be an arbitrary scheme, $M_{1},M_{2}$ be fixed invertible sheaves 
on $X$ and $E_{\infty}$ be a coherent $\CO_X$-module.  
In this paper we will consider framed twisted representations 
of an ADHM quiver 
in the abelian category of coherent $\CO_X$-modules.
As stated in the introduction, we will denote the data $(X,M_1,M_2,E_\infty)$ 
by $\CX$. 
\begin{defi}\label{adhmsheaf}
$(i)$ An ADHM sheaf on $X$ is defined by the data \\
$\CE=(E, \Phi_{1},\Phi_2, 
\phi, \psi)$ where $E$ is a coherent $\CO_X$-module
and 
\[
\Phi_{i}:E\otimes_X M_{i} \to E, \quad
\phi: E\otimes M_{1}\otimes_X M_{2} 
\to E_{\infty}, \quad \psi :E_{\infty} \to E ,
\] 
with $i=1,2$, are morphisms of $\CO_X$-modules satisfying the ADHM relation
\be\label{eq:ADHMrelX}
\Phi_{1}\circ 
(\Phi_{2}\otimes 1_{M_{1}}) - \Phi_{2} \circ (\Phi_{1}\otimes
1_{M_{2}}) + \psi \circ \phi =0. 
\ee
We will refer to $M_{1},M_{2}$ 
as twisting data and $E_{\infty}$ as framing data.

$(ii)$ A morphism $\xi:\CE\to \CE'$ of ADHM sheaves on $X$ 
is a pair $(\xi,z)$ consisting of a morphism of $\CO_X$-modules $\xi:E \to E'$ 
and a complex number $z\in \IC$ 
such that the following diagrams are commutative
\be\label{eq:isomdiag}
\begin{aligned}
& \xymatrix{ 
E\otimes_X M_{i} \ar[r]^-{\Phi_{i}} 
\ar[d]_-{\xi\otimes 1_{M_{i}}} 
& E \ar[d]^-{\xi} \\
E' \otimes M_{i} \ar[r]^-{\Phi'_{i}} & E' \\}
\qquad \qquad 
\xymatrix{ 
E_{\infty} \ar[r]^-{\psi} \ar[d]_-{z 1_{E_{\infty}}} & E \ar[d]^-{\xi}\\
E_{\infty} \ar[r]^-{\psi'} & E' \\}\\
& \qquad \qquad \xymatrix{
E \otimes_Y M_{1}\otimes_X M_{2} \ar[r]^-{\phi}
\ar[d]_-{\xi\otimes 1_{M_{1}\otimes_X M_{2}}} & 
E_{\infty}\ar[d]^-{z 1_{E_{\infty}}}\\
E' \otimes_X M_{1}\otimes_X M_{2} \ar[r]^-{\phi'} & E_{\infty}. \\}\\
\end{aligned}
\ee
A morphism $(\xi,z)$ of ADHM sheaves will be called proper if $z=1$, in
which case $z$ will be omitted. 

$(iii)$ A morphism $(\xi,z): \CE\to \CE'$ will be called invertible 
if $\xi:E\to E'$ 
is an isomorphism of $\CO_X$-modules and $z\neq 0$. A morphism $(\xi,z)$ 
will be called an isomorphism of ADHM sheaves if it is proper and invertible.
\end{defi}

Suppose $X$ is a smooth projective variety over a field $k$ over $\IC$. 
In addition to the twisting data 
$M_1,M_2$ and framing data $E_\infty$, we will also fix a very ample 
line bundle $\CO_X(1)$ on $X$. 
Let $M_{12}$ denote $M_1\otimes_X M_2$.

\begin{defi}\label{adhmstab}
An ADHM sheaf $\CE=(E,\Phi_{1},\Phi_2, \phi,\psi)$ 
on $X$ with twisting data $M_1,M_2$ and framing data 
$E_\infty$ is stable if the following conditions are satisfied
\begin{itemize}
\item[$(i)$] $E$ is a torsion-free $\CO_X$-module 
\item[$(ii)$] $\psi$ is not identically zero 
\item[$(iii)$]  There exists no nontrivial proper saturated subsheaf 
$0\subset E'\subset E$ such that $\Phi_i(E'\otimes_X M_i)\subseteq E'$ 
for $i=1,2$ and $\mathrm{Im}(\psi)\subseteq E'$. 
\end{itemize}
The stable ADHM sheaf $\CE$ will be said to have Hilbert polynomial 
$P$ if $E$ has Hilbert polynomial $P$.  
\end{defi}

\subsection{Basic properties}
Next we establish some properties of ADHM sheaves which will be 
needed later in the paper. We will work under the same assumptions 
as above i.e. $X$ is a smooth projective variety over a field $k$ over $\IC$ 
equipped with a very ample line bundle $\CO_X(1)$. $M_1, M_2$ will be 
fixed twisting data and the locally free $\CO_X$-module 
$E_\infty$ will be fixed framing data. 

\begin{lemm}\label{automlemma}
$(i)$ Let $\xi: \CE\to \CE'$ be a morphism of stable ADHM sheaves on $X$, 
where $\CE, \CE'$ have identical Hilbert polynomials. 
Then $\xi$ is either trivial or invertible. 

$(ii)$
If the base field $k$ is algebraically closed, 
the automorphism group of a stable ADHM sheaf on $X$ is trivial. 
\end{lemm}

{\it Proof.} Let $(\xi,z):E\to E'$ be a morphism of stable ADHM sheaves  
so  the diagrams 
\eqref{eq:isomdiag} for $\CE, \CE', \xi$ are commutative. 

We first claim that $z=0$ $\Rightarrow$ $\xi=0$. If $z=0$, diagrams \eqref{eq:isomdiag}
imply that 
\be\label{eq:automA}
\Phi_i(\mathrm{Ker}(\xi)\otimes_X M_i) \subseteq \mathrm{Ker}(\xi),\qquad
\mathrm{Im}(\psi)\subseteq \mathrm{Ker}(\xi).
\ee
Since $E$ and $E'$ are torsion free, the kernel is automatically saturated, 
$E/\mathrm{Ker}(\xi)\simeq \mathrm{Im}(\xi)\subseteq E'$. 
Since $\CE$ is assumed to be a stable ADHM sheaf, $\mathrm{Ker}(\xi)$ must be 
either $0$ or $E$ in order to avoid a contradiction. The first case, 
$\mathrm{Ker}(\xi)=0$ is ruled out by the stability condition, which requires 
$\psi$ to be nontrivial. Therefore $\mathrm{Ker}(\xi)=E$, and $\xi=0$. 

Suppose $z\neq 0$. 
Then diagrams \eqref{eq:isomdiag} imply that 
\be\label{eq:automA}
\Phi'_i(\mathrm{Im}(\xi)\otimes_X M_i) \subseteq \mathrm{Im}(\xi),\qquad
\mathrm{Im}(\psi')\subseteq \mathrm{Im}(\xi). 
\ee
Hence $\xi$ cannot be identically 0 since $\psi'$ is not, by assumption. 
Moreover, the first equation in \eqref{eq:automA} 
implies that the saturation ${\overline {\mathrm{Im}(\xi)}}$ is preserved 
by $\Phi'_{1},\Phi'_2$ as well. If this were not the case, at least one of the morphisms 
$\Phi'_i$, $i=1,2$, would induce 
a nontrivial morphism 
\[
({\overline {\mathrm{Im}(\xi)}}/ \mathrm{Im}(\xi))\otimes_X M_i \to E'/{\overline {\mathrm{Im}(\xi)}}
\]
which  leads to a contradiction. 
Therefore, since $\CE'$ is stable by assumption, 
 the saturation 
${\overline{\mathrm{Im}(\xi)}}$
must be equal to $E'$. Since $E,E'$ have identical Hilbert polynomials, 
this implies that $\xi$ must be an isomorphism. 
This proves Lemma (\ref{automlemma}.$i$). 
Lemma  (\ref{automlemma}.$ii$) then follows by a standard argument 
(see for example 
\cite[Cor. 1.2.8]{huylehn}).

\hfill $\Box$

\begin{lemm}\label{unstablemma}
Let $\CE=(E,\Phi_1,\Phi_2,\phi,\psi)$ be an ADHM sheaf on $X$ (not necessarily stable)
with $E$ a nontrivial torsion free $\CO_X$-module. Then 
there exists a nontrivial saturated subsheaf 
$0\subset E_0\subseteq E$ determined by the data $(\Phi_1,\Phi_2,\psi)$ such that 
\begin{itemize} 
\item[$(i)$] $\Phi_i(E_0\otimes_X M_i) \subseteq E_0$ for $i=1,2$, and 
$\mathrm{Im}(\psi) \subset E_0$. 
\item[$(ii)$] Any nontrivial saturated subsheaf $0\subset E'\subseteq E$ satisfying analogous conditions contains $E_0$.  
\end{itemize}
\end{lemm}

{\it Proof.} 
If $\Phi_i=0$, 
$i=1,2$, we take $E_0$ to be the saturation 
of $\mathrm{Im}(\psi)$ in $E$. 

Now suppose $\Phi_1,\Phi_2$ are not simultaneously trivial. 
Let $I=(i_1,i_2,\ldots, i_k)$ be an ordered collection of 
indices $i_l\in \{1,2\}$ for $l=1,\ldots, k$. 
For each such ordered collection, let 
\[
M_I = M_{i_1}\otimes_X \cdots \otimes_X M_{i_k} 
\]
and 
\[
\xymatrix{
\Phi_I : E\otimes M_I \ar[rrrrrr]^-{\Phi_{i_k}\circ(\Phi_{i_{k-1}}\otimes 
1_{M_{k}}) \circ \cdots \circ (\Phi_{i_1} \otimes 1_{M_{i_2}\otimes 
\cdots \otimes M_{i_{k}}})} & & & & & & E.\\}
\]
To the empty collection $I=\emptyset$ we formally assign 
$M_I=\CO_X$ and $\Phi_I= 1_E$.
Then set 
\be\label{eq:invsubmod}
E'_0 = \sum_{I} \Phi_I(\mathrm{Im}(\psi)\otimes_X M_I).
\ee
where the sum is over all finite collections $I$ of arbitrary length.
According to \cite[5.3.4]{EGAI}, $E'_0$ is a coherent 
submodule of $E$. We also have 
\[
\Phi_i(E'_0\otimes_X M_i) \subseteq E'_0, \qquad 
\mathrm{Im}(\psi)\subseteq E'_0
\]
by construction.

Let $E_0$ be the saturation of $E_0'$ in $E$. 
Obviously, $E_0$ contains $\mathrm{Im}(\psi)$ as a subsheaf.
If $E_0=E$, then obviously $\Phi_i(E_0\otimes_X M_i) \subseteq 
E_0$ for $i=1,2$.
If $E_0\neq E$, 
consider the following composition of morphisms of $\CO_X$-modules 
\[
\Phi_{i,0}:E_0\otimes_X M_i {\buildrel \Phi_i\over \longto} 
E \twoheadrightarrow  E/E_0.
\]
for $i=1,2$. Since $E_0$ is the saturation of $E_0'$ it follows that 
$\Phi_{i,0}$ vanishes generically on $X$, hence $\mathrm{Im}(\Phi_{i,0})$ 
must be a torsion $\CO_X$-module. However $E/E_0$ is torsion 
free since $E_0$ is saturated, hence $\mathrm{Im}(\Phi_{i,0})$ 
must be trivial. Therefore $\Phi_{i,0}=0$ and we have 
$\Phi_i(E_0\otimes_X M_i) \subseteq 
E_0$ for $i=1,2$. 

Property $(ii)$ in Lemma (\ref{unstablemma}) is obvious by construction.

\hfill $\Box$

\begin{lemm}\label{zerophi}
Suppose $E_\infty$ is a simple locally free sheaf on $X$.   
Then any stable ADHM sheaf $\CE$ on $X$ must have $\phi=0$.
\end{lemm}

{\it Proof.}
Since $E_\infty$ is simple, 
the endomorphism algebra $\mathrm{End}(E_\infty)$ is isomorphic 
to $\IC$. In particular, for any endomorphism 
$\xi_\infty:E_\infty \to E_\infty$, $\mathrm{Tr}(\xi_\infty)=0$ 
implies $\xi_\infty=0$. 
Then the inductive argument used in the proof of 
lemma \cite[Lemma 2.8]{hilblect} 
applies without modification to any stable ADHM sheaf $\CE$
with $E_\infty$ simple.  
This proves that 
the restriction $\phi|_{E_0'\otimes_X M_{12}}$ is trivial, 
where $E_0'$ is the subsheaf of $E$ constructed in \eqref{eq:invsubmod}.

Then we claim that the restriction $\phi|_{E_0\otimes_X M_{12}}$
also vanishes, where $E_0$ is the saturation of $E_0'$ as in Lemma 
(\ref{unstablemma}). By construction, it is clear that the restriction 
$\phi|_{E_0\otimes_X M_{12}}$ vanishes generically on $X$. 
Therefore its image $\mathrm{Im}\left(\phi|_{E_0\otimes_X M_{12}}
\right)$ must be a torsion sheaf on $X$. This is impossible 
since $E_\infty$ is locally free. Therefore 
$\phi|_{E_0\otimes_X M_{12}}$ must be trivial. 
However, the stability condition implies that $E_0=E$, hence 
$\phi=0$.

\hfill $\Box$

\subsection{Moduli spaces} 
Let $X$ to be a smooth complex projective variety 
equipped with a very ample line bundle $\CO_X(1)$. Let  
$M_1,M_2$, $E_\infty$ be fixed twisting, respectively 
framing data on $X$. Let $\CX$ denote the data $(X,M_1,M_2,E_\infty)$. 
Flat families of ADHM sheaves are defined in the usual way. 

\begin{defi}\label{flatfam} 
$(i)$ A flat family of ADHM sheaves on $X$ parameterized by $S$ 
is an ADHM sheaf $\CE_S$ on $X_S$ with twisting data 
$(M_1)_S, (M_2)_S$ and framing data $(E_\infty)_S$
such that $E_S$ is flat over $S$. $\CE_S$ is a flat family of 
stable ADHM sheaves on $X$ if  and the 
the restriction $\CE_S|_{X_s}$ is a 
stable ADHM sheaf on $X_s$ 
for any point $s\in S$.

$(ii)$ Two flat families of ADHM sheaves $\cE_S=(E_S,\Phi_{S,1,2},
\phi_S, \psi_S)$ $\cE'_S=(E'_S,\Phi'_{S,1,2},\phi'_S, \psi'_S)$
are isomorphic if they are isomorphic as ADHM sheaves on $X_S$. 
\end{defi}

As stated in Theorem (\ref{DMstack}), the groupoid defined by flat families of ADHM sheaves 
on $X$ with fixed Hilbert polynomial $P$ is isomorphic to a quasi-projective scheme $\mfm_{ADHM}(\CX,P)$. This follows from \cite[Thm. 2.9.2.44, Thm. 2.9.2.45]{GIT-decorated} where it is proven that
there is quasi-projective coarse moduli space parameterizing isomorphism classes of stable 
ADHM sheaves on $X$ with fixed Hilbert polynomial. Lemma (\ref{automlemma}) proves that stable 
ADHM sheaves have trivial stabilizers, which implies 
that this is a fine moduli space according to \cite[Thm. 4.2.15]{huylehn}.

\section{Torus actions on the moduli space}\label{stacksection}

Since the moduli spaces of stable ADHM sheaves are not proper, an algebraic 
torus action {\bf T} with proper fixed loci will be required in order to obtain a well defined residual theory 
in the sense of \cite{BP}. Moreover, the construction of a {\bf T}-equivariant 
perfect tangent-obstruction theory for the moduli space   
as defined in \cite[Def. 2.1]{degGW} also requires an \'etale cover 
$\{S_\alpha\}$ consisting of finitely many {\bf T}-equivariant affine charts.
The following lemma shows that such actions can be naturally obtained using the 
GIT construction of the moduli space given in \cite[Sect. 2.9]{GIT-decorated}. 

\begin{lemm}\label{Tactioncoro} 
Let $Z$ be a complex quasi-projective scheme equipped with an ample line bundle 
${\mathcal L}$. Suppose ${\bf T}\times Z\to Z$ is an algebraic torus action on $Z$ and 
${\mathcal L}$ is equipped with a {\bf T}-linearization. Then there exists a Zariski open cover of 
$Z$ consisting of finitely many {\bf T}-invariant affine open charts. 
\end{lemm} 

{\it Proof.} This result has been proven in \cite[Sect. 3, Cor. 2]{equiv-compl}. 
There it is assumed that $Z$ is a normal variety in order to infer the existence 
of an equivariant embedding into a projective space. This assumption is not needed 
in the present case since the existence of a {\bf T}-linearized ample line bundle 
is already part of the assumptions. Then the proof is identical to the proof of 
\cite[Sect. 3, Cor. 2]{equiv-compl}.

\hfill $\Box$ 

Now recall that the moduli space of stable ADHM sheaves has been 
constructed in \cite[Thm. 2.9.2.44, Thm. 2.9.2.45]{GIT-decorated} using GIT 
methods. In fact, according to \cite[Lemm. 2.9.2.43]{GIT-decorated}, 
the moduli problem for ADHM sheaves reduces to special case of 
a more general moduli problem for decorated bundles called twisted affine 
bumps studied in detail 
in \cite[Sect.2.8]{GIT-decorated}. The moduli problem is formulated in 
\cite[Sect. 2.8.1]{GIT-decorated} and the construction of the parameter space 
for such decorated bundles is presented in detail in \cite[Sect. 2.8.3]{GIT-decorated}. 
In particular there is a quasi-projective parameter 
scheme $R$ for stable ADHM sheaves equipped with an affine reductive 
group action $G\times R\to R$ and a $G$-linearized line bundle $L$. 
Moreover the data $(R,G,L)$ determines an ample line bundle ${\mathfrak L}$ 
on the GIT quotient ${\mathfrak M}_{ADHM}(\CX,P)$ according to 
\cite[Thm 1.10]{GIT}, \cite[Rem. 1.4.3.9]{GIT-decorated}.
Since this construction is very technical and relies on results in the previous 
sections of \cite{GIT-decorated}, the details will not be reviewed here.
It suffices to note that there are natural torus actions on $R$ 
acting linearly on the decoration data $(\Phi_{1},\Phi_2,\phi,\psi)$ which lift naturally to 
{\bf T}-linearizations of  $L$. 
Any such action then descends to a torus action on the 
moduli space satisfying the conditions of lemma (\ref{Tactioncoro}) with respect to the 
ample line bundle  ${\mathfrak L}$. Two examples are presented below.

\begin{exam}\label{torusactA} 
Let ${\bf T}=\IC^\times \times \IC^\times$ and define an 
action ${\bf T}\times \mfm_{ADHM}(\CX,P)\to \mfm_{ADHM}(\CX,P)$ by
\be\label{eq:toractA}
(E_S, \Phi_{S,1}, \Phi_{S,2}, \phi_S, \psi_S) 
\to (E_S, t_1\Phi_{S,1}, t_2\Phi_{S,2}, t_1t_2\phi_S, \psi_S) 
\ee
for any $(t_1,t_2)\in {\bf T}(S)$ and any flat family 
$\CE_S= (E_S, \Phi_{S,1},\Phi_{S,2}, \phi_S, \psi_S)$ 
of stable ADHM sheaves. 
\end{exam}   

\begin{exam}\label{torusactB}
Suppose $E_\infty = \oplus_{a=1}^{r_\infty} L_a$, where 
$r_\infty\in \IZ$, $r_\infty\geq 1$, and 
$L_a$, $a=1,\ldots, r_\infty$, are invertible $\CO_X$-modules. 
Let ${\bf T} = \IC^\times \times \IC^\times \times 
(\IC^\times)^{r_\infty-1}$. There is a natural torus action 
\[
  (\IC^\times)^{r_\infty-1} \times E_\infty \to E_\infty 
\]
given by the identification of $(\IC^\times)^{r_\infty-1}$
with the hypersurface 
\[
\prod_{a=1}^{r_\infty} z_a =1 
\]
in $(\IC^\times)^{r_\infty}$. 
This provides a $(\IC^{\times})^{r_\infty-1}$ 
action $\xi_\infty: (\IC^{\times})^{r_\infty-1} \to \mathrm{Aut}(E_\infty)$ 
given by 
\[
\xi_\infty({\underline z})|_{L_a} = z_a 1_{L_a} 
\]
for any $a=1,\ldots,  r_\infty$, where ${\underline z}
=(z_a)_{a=1,\ldots,r_\infty}\in (\IC^{\times})^{r_\infty-1}$.
Using the above presentation of the torus $(\IC^\times)^{r_\infty-1}$, 
its characters will be denoted in the following by 
${\underline m}=(m_1,\ldots, m_{r_\infty})\in \IZ^{r_\infty}$.

If $S$ is a scheme of finite type over $\IC$, we will denote by 
$(\xi_\infty)_S: {(\IC^{\times})}^{r_\infty-1}(S) \to \mathrm{Aut}(E_S)$ 
the relative version of this action. 
Then we have a natural torus action on $\mfm_{ADHM}(\CX,P)$ given by 
\be\label{eq:toractB}
(E_S, \Phi_{S,1}, \Phi_{S,2}, \phi_S, \psi_S) 
\to (E_S, t_1\Phi_{S,1}, t_2\Phi_{S,2}, t_1t_2(\xi_\infty)_S({\underline z})
\circ 
\phi_S , \psi_S\circ \xi_\infty({\underline z})^{-1}) 
\ee
for any $(t_1,t_2, {\underline z})\in {\bf T}(S)$ and any flat family 
$\CE_S= (E_S, \Phi_{S,1},\Phi_{S,2}, \phi_S, \psi_S)$ 
of stable ADHM sheaves.
\end{exam}

\begin{rema} 
Note that the actions defined in Examples (\ref{torusactA}), (\ref{torusactB}) 
coincide if $E_\infty$ is a line bundle on $X$. 
\end{rema} 

The following proposition shows that the fixed loci 
of the action defined in (\ref{torusactB}) are proper if $X$ is a curve. 
In this case the Hilbert polynomial $P$ is determined by the 
rank $r\geq 1$ and the degree $e\in \IZ$, and the  
moduli spaces will be denoted by $\mfm_{ADHM}(\CX,r,e)$.

\begin{prop}\label{fixedloci} 
Let $X$ be a smooth projective curve over $\IC$, 
and suppose $E_\infty = \oplus_{a=1}^{r_\infty}L_a$ as in 
Example (\ref{torusactB}). Then the fixed locus 
$\mfm(\CX,r,e)^{\bf T}$ of the torus action given in Example (\ref{torusactB})
is a projective scheme ver $\IC$. 
\end{prop} 

{\it Proof.} 
The stack theoretic fixed locus for the {\bf T}-action on $\mfm_{ADHM}(\CX,r,e)$
can be determined using \cite[Prop. 2,6]{gract-stacks}, which will be briefly recalled below. 
Let $S$ be a scheme of finite type over $\IC$. 
A fixed point in $\mfm_{ADHM}(\CX,r,e)(S)$ 
is a flat family of stable ADHM sheaves $\CE_S$ on $X_S$ 
equipped with an isomorphism of ADHM sheaves 
\[
\begin{aligned}
\xi_S(t_1,t_2,{\underline z}) : (E_S,\Phi_{S,1},\Phi_{S,2},\phi_S,\psi_S) 
\to (E_S, t_1\Phi_{S,1}, t_2\Phi_{S,2}, t_1t_2\phi_S, \psi_S)\\
\end{aligned}
\]
for each $(t_1,t_2,{\underline z})\in {\bf T}(S)$ 
satisfying a natural cocycle identity. 

In particular, we have a group homomorphism $\xi_S: {\bf T}(S) 
\to \mathrm{Aut}(E_S)$ i.e. a ${\bf T}(S)$-action on $E_S$. 
Therefore, according to \cite{Demazure}, 
$E_S$ decomposes in a direct sum of $\CO_{X_S}$-modules 
\be\label{eq:chardecomp}
E_S = \bigoplus_{(n_1,n_2)\in \IZ^2} 
\bigoplus_{{\underline m}\in \IZ^{r_\infty}}E_S{(n_1,n_2,{\underline m})},
\ee
where only finitely many terms are nontrivial. 
An element $(t_1,t_2,{\underline z})\in {\bf T}(S)$ 
acts by multiplication by $t_1^{n_1}t_2^{n_2} z_1^{m_1} 
\ldots z_{r_\infty}^{m_{r_\infty}}$ 
on $E_S{(n_1,n_2,{\underline m})}$.
Moreover each direct summand $E_S{(n_1,n_2,{\underline m})}$ 
must be flat over $S$ since $E_S$ is flat over $S$.

By definition the restriction of $\CE_S$ to any fiber $X_s$, 
$s\in S$, is a stable ADHM sheaf $\CE_s$ on the smooth projective curve $X_s$ 
over the residual field $k(s)$. Since $E_S$ is flat over $S$, the restriction 
$\xi_S(t_1,t_2,{\underline z})|_{X_s}$ is an isomorphism 
of stable ADHM sheaves 
\[
\xi_s(t_1,t_2,{\underline z}) : (E_s,\Phi_{s,1},\Phi_{s,2},\phi_s,\psi_s) 
\to (E_s, t_1\Phi_{s,1}, t_2\Phi_{s,2}, t_1t_2\phi_s, \psi_s).
\]
Then the  commutative diagrams \eqref{eq:isomdiag} imply that the only nontrivial 
components of $\Phi_{s,1},\Phi_{s,2}, \phi_s, \psi_s$ with respect to the character 
decomposition \eqref{eq:chardecomp} 
are 
\be\label{eq:fixedlociB}
\begin{aligned}
\Phi_{s,1}(n_1,n_2,{\underline m}) : E_s(n_1,n_2,{\underline m}) 
\otimes_{X_s} (M_1)_s & \to E_s(n_1+1,n_2,{\underline m}) \\
\Phi_{s,2}(n_1,n_2,{\underline m}) : 
E_s(n_1,n_2,{\underline m}) \otimes_{X_s} (M_2)_s & 
\to E_s(n_1,n_2+1,{\underline m}) \\
\phi_s(-1,-1, {\underline \delta}^a) : E_s(-1,-1,{\underline \delta}^a) 
\otimes_{X_s} (M_{12})_s & 
\to (L_a)_s \\
\psi_s(a):(L_a)_s & \to  E_s(0,0,{\underline \delta}^a),\\
\end{aligned}
\ee 
where 
\[
{\underline \delta}^a = (\underbrace{0,0,\ldots, 0}_{a-1},  -1,
\underbrace{0,\ldots,0}_{r_\infty -a})
\]
for $a=1,\ldots, r_\infty$.
In particular note that 
\[
\mathrm{Im}(\psi_s) \subseteq \bigoplus_{a=1}^{r_\infty} E_s(0,0,{\underline 
\delta}^a).
\]
Moreover, the following relations must also be satisfied 
\be\label{eq:fixedadhmrel}
\begin{aligned}
& \Phi_{s,2}(n_1+1,n_2,{\underline m})\circ 
(\Phi_{s,1}(n_1,n_2,{\underline m})\otimes 1_{(M_2)_s}) = \\
& \Phi_{s,1}(n_1,n_2+1,{\underline m}) 
\circ (\Phi_{s,2}(n_1,n_2,{\underline m})\otimes 1_{(M_1)_s})
\end{aligned}
\ee
for any $(n_1,n_2)\in \IZ_{\geq 0}^2$ and any ${\underline m}\in 
\IZ^{r_\infty}$.

Then it follows from the first two equations in \eqref{eq:fixedlociB} 
that the canonical destabilizing subsheaf $(E_{s})_0$ constructed in 
Lemma (\ref{unstablemma}) is a subsheaf of 
\[
\bigoplus_{n_1,n_2\geq 0} \bigoplus_{a=1}^{r_\infty} E_s(n_1,n_2,{\underline 
\delta}^a).
\]
Given the direct sum decomposition \eqref{eq:chardecomp}, 
where the terms are torsion free sheaves,  the stability condition implies 
that 
\[
E_s(n_1,n_2,{\underline m})=0
\]
if $n_1<0$ or $n_2<0$ or ${\underline m} \neq {\underline \delta}^a$. 
Otherwise, the saturation 
${\overline {(E_s)}_0}$ would be a proper subsheaf of $E_s$. 
Therefore we must have $\phi_s=0$ for any fixed ADHM sheaf. 

Taking into account relations \eqref{eq:fixedadhmrel},
 same argument implies that the saturations  
\[
\begin{aligned}
{\overline{ 
\mathrm{Im}(\Phi_{s,1}(n_1,n_2,{\underline \delta}^a))}}
& \subseteq E_s(n_1+1,n_2,
{\underline \delta}^a)\\
{\overline{ 
\mathrm{Im}(\Phi_{s,2}(n_1,n_2,{\underline \delta}^a))}}
& \subseteq
E_s(n_1,n_2+1,{\underline \delta}^a)\\
\end{aligned}
\] 
must be equal to $E_s(n_1+1,n_2,
{\underline \delta}^a)$, respectively $E_s(n_1,n_2+1,{\underline \delta}^a)$ 
for any $(n_1,n_2)\in \IZ_{\geq 0}^2$, and any $a=1,\ldots, r_\infty$. 
Moreover the saturation 
\[
{\overline {\mathrm{Im}(\psi_s(a))}}\subseteq
E_s(0,0,{\underline \delta}^a)
\]
must also be equal to 
$E_s(0,0,{\underline \delta}^a)$ for any $a=1,\ldots, r_\infty$. 
Since $(L_a)_s$, $a=1,\ldots, r_\infty$ are line bundles on $X_s$, 
and $X_s$ is a 
smooth projective curve over $s$, this implies that 
the data $(E_s(n_1,n_2,{\underline \delta}^a), \Phi_{s,i}(n_1,n_2,a), 
\psi_s(a))$ satisfy the following conditions 
\begin{itemize}
\item[$(a)$] All nontrivial direct summands 
$(E_s(n_1,n_2,{\underline \delta}^a)$ in 
\eqref{eq:chardecomp} are invertible $\CO_{X_s}$-modules.
\item[$(b)$] 
The morphism of $\CO_{X_s}$-modules $\psi_s(a)$ must be injective for all 
$a=1,\ldots,r _\infty$. 
\item[$(c)$] The morphisms of $\CO_{X_s}$-modules 
$\Phi_s(n_1,n_2,{\underline \delta}^a)$ 
must be injective 
for any $a=1,\ldots, r_\infty$
and any $(n_1,n_2)\in \IZ_{\geq 0}^2$ such that the target is nontrivial. 
\item[$(d)$] Relations \eqref{eq:fixedadhmrel} must 
be satisfied. 
\end{itemize}

Note that the degree of each line bundle $E_s(n_1,n_2,{\underline \delta}^a)$
is constant along each connected component of $S$ since 
$E_S(n_1,n_2,{\underline \delta}^a)$ are flat over $S$. For simplicity
we will assume in the following that $S$ is connected. 
The arguments below generalize obviously to the case when 
$S$ has several connected components. 
Let $e(n_1,n_2,a)= \mathrm{deg}\, E_S(n_1,n_2,{\underline \delta}^a)$. 
For each $a=1,\ldots, r_\infty$ 
let  $\Delta_a$ the set of pairs $(n_1,n_2)\in \IZ^{\times}$ 
such that the direct summand 
$E_S(n_1,n_2,{\underline \delta}^a)$  
in \eqref{eq:chardecomp} is nontrivial. Note that properties $(b)$ and $(c)$ above 
imply that for fixed $a$ the set $\Delta_a\subset \IZ^2$ of integral points 
$(n_1,n_2)$ labeling the summands $E_S(n_1,n_2,a)$ is naturally 
identified with a Young diagram embedded in the first quadrant. 

Since Theorem (\ref{DMstack}) proves that 
the stack  $\mfm_{ADHM}(\CX,r,e)$ 
is isomorphic to a quasi-projective scheme, 
the stack-theoretic fixed locus is isomorphic to the scheme theoretic fixed locus 
\cite[Rem. (3.4.ii)]{gract-stacks}. According to \cite[Prop. A.8.19(1)]{pseudored}, 
the later is a closed subscheme 
of $\mfm_{ADHM}(\CX,r,e)$ since $\mfm_{ADHM}(\CX,r,e)$ is 
quasi-projective, in particular separated, over $\IC$.
Given the above structure results, it follows that 
the functor of points of each connected component of $\mfm_{ADHM}(\CX,r,e)^{\bf T}$ 
is isomorphic to a direct product of moduli functors $\CZ_{\Delta, {\sf e}}$ as defined below. 

For any Young diagram $\Delta \subset \IZ^2$, any fixed line bundle $L$ on $X$ 
of degree $d$, and any set of integers ${\sf e}= 
\{e(n_1,n_2) \in \IZ\, |\, (n_1,n_2)\in \Delta\}$ 
as above let  
$\CZ_{\Delta,{\sf e}}$ be the moduli functor  assigning to a scheme $S$ over 
$\IC$ the set of isomorphism classes of data  
\[
\CL_S=\left(L_S(n_1,n_2), \Phi_{S,1}(n_1,n_2), \Phi_{S,2}(n_1,n_2),\psi_S\right),\]
with $(n_1,n_2)\in \Delta$, satisfying the following conditions .
\begin{itemize} 
\item[(c.1)] $L_S(n_1,n_2)$ is a line bundle on $X_S$, 
such that the degree of $L_S(n_1,n_2)|_{X_s}$ is $e(n_1,n_2)$ for any $s\in S$.
\item[(c.2)] $\Phi_{S,1}(n_1,n_2):L_S(n_1,n_2)\otimes_{X_S} (M_1)_S \to 
L_S(n_1+1,n_2)$, 
$\Phi_{S,2}(n_1,n_2):L_S(n_1,n_2)\otimes_{X_S} (M_2)_S \to L_S(n_1,n_2+1)$, 
$\psi_S : (L)_S \to L_S(0,0)$ are 
morphisms of $\CO_{X_S}$-module 
\item[(c.3)] The restrictions of the morphisms $\Phi_{S,1}(n_1,n_2), \Phi_{S,2}(n_1,n_2),\psi_S$, $(n_1,n_2)\in \Delta$, to $X_s$ are injective morphisms of $\CO_{X_s}$ modules whenever their target is nontrivial, for any point $s\in S$. 
\item[(c.4)] $\Phi_{S,2}(n_1+1,n_2)\circ (\Phi_{S,1}(n_1,n_2)\otimes 1_{(M_2)_S}) = 
\Phi_{S,1}(n_1,n_2+1) \circ (\Phi_{S,2}(n_1,n_2)\otimes 1_{(M_1)_S})$ for
any $(n_1,n_2)\in \Delta$. Here a morphism $\Phi_{S,i}(n_1,n_2)$, $i=1,2$ is by convention 
identically zero if $(n_1,n_2)\notin \Delta$. 
\end{itemize} 
Two such data $\CL_S,\CL'_S$ are isomorphic if there exist isomorphisms 
\[
\xi_S(n_1,n_2): L_S(n_1,n_2,a) 
{\buildrel \sim \over \longto} L'_S(n_1,n_2,a)
\] of line bundles on $X_S$ satisfying the 
obvious compatibility conditions with the decoration data. Note that condition 
$(c.1)$ does not explicitely require $L_S$ to be flat over $S$ since this is automatically satisfied. 
$L_S$ is flat over $X_S$ since it is locally free, and $X_S$ is flat over $S$ since 
$X_S= X\times S$. 

Therefore in order to complete the proof of Proposition \ref{fixedloci} it suffices to 
show that each moduli functor $\CZ_{\Delta, {\sf e}}$ is isomorphic to a 
projective scheme. In fact it will be identified with a closed subscheme of a 
product of symmetric powers of $X$.

Let $d_i=\mathrm{deg}(M_i)$, $i=1,2$. 
Note that the data $\CL_S$ satisfying $(c.1)$-$(c.4)$ is empty unless $d\leq e(0,0)$, $e(n_1,n_2)\leq e(n_1+1,n_2)-d_1$ for all $(n_1,n_2)\in \Delta$ 
such that $(n_1+1,n_2)\in \Delta$, and 
$e(n_1,n_2)\leq e(n_1,n_2+1)-d_2$ for all $(n_1,n_2)\in \Delta$ 
such that $(n_1,n_2+1)\in \Delta$. These conditions will be assumed in the following. 
Then let  $S_{\Delta,{\sf e}}(X)$ denote the following direct product 
\be\label{eq:symmdirprod}
\bal 
S_{\Delta,{\sf e}}(X)  = S^{e(0,0)-d}(X)\times &
\mathop{\sf X}_{\substack{(n_1,n_2)\in \Delta\\ (n_1+1,n_2)\in \Delta}} S^{e(n_1+1,n_2)-e(n_1,n_2)-d_1}(X) \times \\
& \mathop{\sf X}_{\substack{(n_1,n_2)\in \Delta\\ (n_1,n_2+1)\in \Delta}} S^{e(n_1,n_2+1)-e(n_1,n_2)-d_2}(X)
\eal
\ee
where $S^n(X)$, $n\geq 0$,  denotes the $n$-th symmetric product of $X$. 
$S_{\Delta,{\sf e}}(X)$ is a moduli space of collections of effective divisors 
$(D,D_1(n_1,n_2), D_2(n_1,n_2))$ on $X$, where the notation is self-explanatory. 
More precisely, $S_{\Delta,{\sf e}}(X)$ is a fine moduli space for the moduli functor 
assigning to any scheme $S$ over $\IC$ a collection $(D_S, D_{S,1}(n_1,n_2), D_{S,2}(n_1,n_2))$ of effective Cartier divisors on  $X_S$ such that the structure sheaves 
$\CO_{D_S}$, $\CO_{D_{S,i}}(n_1,n_2)$, 
$i=1,2$, are flat over $S$, and their restrictions to any fiber $X_s$, $s\in S$, 
yield a collection of effective Cartier divisors of $X_S$ of fixed degrees determined by   
\eqref{eq:symmdirprod}. 

Next note that there is a morphism of moduli functors $\varrho:\CZ_{\Delta, {\sf e}} 
\to S_{\Delta,{\sf e}}(X)$ where $S_{\Delta,{\sf e}}(X)$ is identified as usual with 
its functor of points. 
Given any set of data $\CL_S$ satisfying conditions $(c.1)$-$(c.4)$ 
consider the Cartier divisors 
$(D_S, D_{S,1}(n_1,n_2), D_{S,2}(n_1,n_2))$ on  $X_S$ such that 
$D_S$ is determined by 
$$\psi_S \in H^0(X_S, (L)_S^{-1}\otimes L_S(0,0))\setminus \{0\},$$ 
$D_{S,1}(n_1,n_2)$ is determined by 
$$
\Phi_{S,1}(n_1,n_2) \in H^0(X_S, (L_S(n_1,n_2)\otimes_{X_S}(M_1)_S)^{-1}\otimes L_S(n_1+1,n_2))\setminus 
\{0\}
$$ 
for any $(n_1,n_2)\in \Delta$ such that $(n_1+1,n_2)\in \Delta$, 
and $D_{S,2}(n_1,n_2)$ is determined analogously. It is clear that the resulting divisors 
depend only on the isomorphism class of $\CL_S$. 
Moreover, for any point $s\in S$, the morphisms $\Phi_{S,1}(n_1,n_2)|_{X_s}, 
\Phi_{S,2}(n_1,n_2)|_{X_s},\psi_S|_{X_s}$, $(n_1,n_2)\in \Delta$, are injective. 
This implies that the structure sheaves $\CO_{D_S}$, $\CO_{D_{S,i}}(n_1,n_2)$, 
$i=1,2$ are flat over $S$. (See for example \cite[Lemm. 2.1.4]{huylehn}.)
Therefore this assignment defines a morphism $\varrho$ as claimed above. 

Furthermore, for any $\CL_S$, the divisors 
$(D_S, D_{S,1}(n_1,n_2), D_{S,2}(n_1,n_2))$ satisfy the relations
\be\label{eq:divrel}
D_{S,2}(n_1+1,n_2) + D_{S,1}(n_1,n_2) = D_{S,1}(n_1,n_2+1) + D_{S,2}(n_1,n_2).
\ee
Since the relations \eqref{eq:divrel} are closed conditions, they determine a closed 
subscheme $S'_{\Delta,{\sf e}}(X) \subset S_{\Delta,{\sf e}}(X)$. It will be shown next that 
the morphism $\varrho$ determines an isomorphism 
$\CZ_{\Delta, {\sf e}} {\buildrel \sim \over \longto}  S'_{\Delta,{\sf e}}(X)$.

First note that if two families $\CL_S, \CL'_S$ determine the same collection of divisors 
on $X_S$, it is straightforward to check that they must be isomorphic. 
Next let $(D_S, D_{S,1}(n_1,n_2), D_{S,2}(n_1,n_2))$ be a collection of Cartier 
divisors on $X_S$ satisfying 
\eqref{eq:divrel}, such that their structure sheaves are flat over $S$, and their restrictions 
to each fiber $X_s$ are Cartier divisors of degrees specified by \eqref{eq:symmdirprod}. 
Then 
one can reconstruct the family $\CL_S$ up to isomorphism as follows.
Let $(n_1,n_2)\in \Delta\setminus \{(0,0)\}$. Choose an ordered sequence $(n_{1,j}, n_{2,j})\in \Delta$, 
$j=0,\ldots,n_1+n_2$ such that $(n_{1,0}, n_{2,0})=(0,0)$, 
$(n_{1,n_1+n_2}, n_{2,n_1+n_2})=(n_1,n_2)$, 
and 
$$
n_{1,j+1}\geq n_{1,j}, \qquad n_{2,j+1}\geq n_{2,j},\qquad
|n_{1,j+1}-n_{1,j}|+|n_{2,j+1}-n_{2,j}|=1$$
 for any $j=0,\ldots,n_1+n_2-1$. 
For each $j=0,\ldots, n_1+n_2-1$, let 
\[
\epsilon_j = \left\{\begin{array}{ll} 
1& \mathrm{if}\ |n_{1,j+1}-n_{1,j}|=1, \ |n_{2,j+1}-n_{2,j}|=0\\
2 & \mathrm{if}\ |n_{1,j+1}-n_{1,j}|=0, \ |n_{2,j+1}-n_{2,j}|=1.\\
\end{array}\right. 
\]
Intuitively, such a collection determines an ascending zig-zag path from $(0,0)$ to 
$(n_1,n_2)$ 
in $\IR^2$ consisting of a sequence of horizontal or vertical unit segments joining 
successive points in $\Delta$. 
Then set  
\[
{\widetilde D}_{S}(n_1,n_2) = D_S + \sum_{j=0}^{n_1+n_2-1} D_{S,\epsilon_j}(n_{1,j}, n_{2,j}).
\]
If $(n_1,n_2)=(0,0)$ set ${\widetilde D}_S(0,0)=D_S$.
Note that relations \eqref{eq:divrel} imply that ${\widetilde D}_{S}(n_1,n_2)$ is independent 
of the collection $(n_{1,j}, n_{2,j})\in \Delta$
used in the construction.  Note also that there are $n_1$ values of 
$j\in \{0,\ldots, n_1+n_2-1\}$ such that $\epsilon_j=1$ and 
$n_2$ values of $j$ such that $\epsilon_j=2$. 

For any $(n_1,n_2)\in \Delta\setminus\{(0,0)\}$ 
let $$
L_{S}(n_1,n_2) = \CO_{X_S}({\widetilde D}_{S}(n_1,n_2))\otimes_{X_S} (M_1)_S^{n_1} 
\otimes_{X_S}  (M_2)_S^{n_2}\otimes_{X_S} (L)_S,
$$ 
and let $L_S(0,0)=\CO_{X_S}({\widetilde D}_{S}(0,0))\otimes_{X_S} (L)_S$. 
By construction, for each $(n_1,n_2)\in \Delta$ such that 
$(n_1+1,n_2)\in \Delta$ there exists a nontrivial section 
$$ 
\sigma_{S,1}(n_1,n_2)\in 
H^0(X_S,{L}_{S}(n_1+1,n_2)\otimes_{X_S}  { L}_{S}(n_1,n_2)^{-1})
$$  
whose zero locus is $D_{S,1}(n_1,n_2)$, and 
for each $(n_1,n_2)\in \Delta$ such that 
$(n_1,n_2+1)\in \Delta$ there exists a nontrivial section 
$$ 
\sigma_{S,2}(n_1,n_2)\in 
H^0(X_S,{L}_{S}(n_1,n_2+1)\otimes_{X_S}  { L}_{S}(n_1,n_2)^{-1})
$$  
whose zero locus is $D_{S,2}(n_1,n_2)$. Moreover, given relations 
\eqref{eq:divrel}, these sections can be chosen such that 
\[
\sigma_{S,2}(n_1+1,n_2)\otimes \sigma_{S,1}(n_1,n_2)) =
 \sigma_{S,1}(n_1,n_2+1) \otimes \sigma_{S,2}(n_1,n_2).
\]
Choose also a section $\sigma_S(0,0)$ of $\CO_{X_S}({\widetilde D}_S(0,0))$ defining 
${\widetilde D}_S(0,0)$. 
The sections $\sigma_S(0,))$, $\sigma_{S,i}(n_1,n_2)$, $i=1,2$ determine 
morphisms $\psi_S$, $\Phi_{S,i}(n_1,n_2)$, $i=1,2$ 
satisfying properties $(c.2)-(c.4)$
Therefore this 
construction determines a set of data $\CL_S$. 
The choice of sections is not canonical, but it is straightforward to check that different 
choices result in isomorphic data $\CL_S$. Therefore one obtains indeed an 
isomorphism between $\CZ_{\Delta {\sf e}}$ and $S'_{\Delta, {\sf e}}(X)$. 

In conclusion, 
each connected components of the fixed locus 
 $\mfm_{ADHM}(\CX,r,e)^{\bf T}$ is isomorphic to a 
product 
\[
\mathop{\sf X}_{a=1}^{r_\infty} S'_{\Delta_a, {\sf e}_a}(X)
\]
 for some Young diagrams $\Delta_a$, $a=1, \ldots, r_\infty$ and some numerical 
 invariants ${\sf e}_a=\{e(n_1,n_2,a) \, |\, (n_1,n_2)\in \Delta_a\}$
such that 
\[
\sum_{a=1}^{r_\infty} |\Delta_a|= r, \qquad 
\sum_{a=1}^{r_\infty} \sum_{(n_1,n_2)\in \Delta_a} e(n_1,n_2,a) = e.
\]

\hfill $\Box$

\section{Deformation Theory of ADHM Sheaves on Curves}

In this section we provide some basic results in the deformation theory 
of ADHM sheaves on curves. Therefore in the following $X$ will be a 
fixed smooth projective curve over $\IC$ of genus $g$. 
Our main 
references for deformation theory are \cite{artin-moduli}, 
\cite{artin-versal}, \cite{LT}. 

The following observation will be very useful throughout this section.  

\begin{lemm}\label{locfree}
Suppose $X$ is a smooth complex projective curve. 
Let $\CE_S$ be a flat family of stable ADHM sheaves on $X$ parameterized 
by a scheme $S$ of finite type over $\IC$. Then $E_S$ is a locally 
free $\CO_{X_S}$-module. 
\end{lemm}

{\it Proof.} Since $\CE_S$ is a flat family of stable ADHM 
sheaves, it follows that $E_S|_{X_s}$ is a torsion 
free $\CO_X$-module, for any point $s\in S$. Since 
$X$ is a smooth projective curve, 
this implies that $E_S|_{X_s}$
is locally free. Then Lemma (\ref{locfree}) follows 
from \cite[Lemma 2.1.7]{huylehn} since $E_S$ is flat over $S$. 

\hfill $\Box$

\subsection{Deformation Complex}
Let $S$ be a scheme of finite type over $\IC$. 
In the following we will call a nilpotent thickening of $S$ 
a scheme $S'$ of finite type over $\IC$ such that 
\begin{itemize}
\item $S'$ is a scheme over $S$  
\item there is  closed embedding $S\subset S'$ of schemes over $S$
such that the defining 
ideal sheaf $I$ of $S$ in $S'$ is nilpotent. In particular there is a commutative 
diagram of morphisms of schemes 
\be\label{eq:schemediag}
\xymatrix{
S \ar@{^{(}->}[r]\ar[dr]_-{1_S} & S' \ar[d]^-{\pi_{S'/S}} \\
& S \\}
\ee
where the horizontal top arrow is a closed embedding and $1_S:S\to S$ is the identity.
\end{itemize}
Note that if these conditions are satisfied, the projection $\pi_{S'/S}:
S'\to S$ is a finite morphism. Therefore it follows from the base change
theorem that higher direct images vanish and $\pi_{S'/S*}$ is an 
exact functor from $\CO_{S'}$-modules to $\CO_S$-modules. 
Moreover, the induced projection $(\pi_{S'/S})_X : X_{S'} \to X_S$ 
satisfies the same properties. For simplicity we will not explicitly 
write $((\pi_{S'/S})_{X})_*$ in the following. Given a $\CO_{X_{S'}}$-module 
$F$ we will denote $((\pi_{S'/S})_{X})_* F$ also by $F$ regarded as an 
$\CO_{X_{S}}$-module. The same conventions will apply to morphisms 
of $\CO_{X_{S'}}$-modules. 
In particular we have a canonical exact sequence of $\CO_S$-modules 
\[
0\to I \to \CO_{S'} \to \CO_S\to 0.
\]
Diagram \eqref{eq:schemediag} implies that 
this exact sequence is split and also determines a splitting 
$\CO_{S}$-modules $\CO_{S'}\simeq \CO_S\oplus I$.
In the following we will identify $\CO_{S'}=\CO_S\oplus I$ as 
$\CO_S$-modules.  

Given a coherent $\CO_S$-module $I$, the trivial nilpotent thickening 
of $S$ determined by $I$ is the nilpotent thickening $S'$ so 
that $\CO_{S'}= \CO_S \oplus I$ as $\CO_{S}$-modules, the ring structure 
of $\CO_{S'}$ being determined by 
\be\label{eq:ringstruct}
(s,x) \cdot (s',x') = (ss', sx+s'x')
\ee
for any local sections $s,s'$ of $\CO_S$ and $x,x'$ of $I$. 
In particular $I^2=0$ with respect to the 
multiplicative structure of $\CO_{S'}$.

\begin{defi}\label{adhmfamB} 
$(i)$ Let $\CE_S$ be a flat family of stable ADHM sheaves on $X$ 
parameterized by a scheme $S$ of finite type 
over $\IC$ and let $S'$ be a nilpotent thickening of $S$. 
An extension of $\CE_S$ to $S'$ is a flat family $\CE_{S'}$ 
of stable ADHM sheaves on $X$ such that $\CE_{S'}\otimes_{S'} \CO_{S}
= \CE_S$.

$(ii)$ 
Two extensions $\CE_{S'}^{(1)}$, $\CE_{S'}^{(2)}$ 
of $\CE_S$ to $S'$ are equivalent if there exists an isomorphism 
$\xi_{S'}: \CE_{S'}^{(1)} \to \CE_{S'}^{(2)}$ of ADHM sheaves on 
$X_{S'}$ such that $\xi_{S'}|_{X_S} = 1_{\CE_S}$. 
\end{defi}

\begin{defi}\label{defcomplex}
Let $\CE_S^{(1)}, \CE_S^{(2)}$ be two  
flat families of ADHM sheaves on $X$ 
parameterized by a scheme $S$ of finite type 
over $\IC$ and let $I$ be a coherent $\CO_S$-module. 
We define the  complex $\CC(\CE_S^{(1)},\CE_S^{(2)},I)$ 
to be the following complex of $\CO_{X_S}$-modules  
\be\label{eq:hypercohB} 
\begin{aligned} 
0 & \to \begin{array}{c} \lochom_{X_S}(E_S^{(1)},E_{S}^{(2)}\otimes_{X_{S}} 
I)\\ 
\end{array}
{\buildrel d_{1}\over \longto} 
\begin{array}{c}  \lochom_{X_S}
(E_S^{(1)}\otimes_{X_S}(M_1)_S,E_{S}^{(2)}\otimes_{X_{S}} I)
\\ \oplus \\
\lochom_{{X_S}}(E_{S}^{(1)}\otimes_{{X_S}} (M_2)_S,
E_{S}^{(2)}\otimes_{X_{S}} I)\\ \oplus \\ 
\lochom_{X_S}(E_S^{(1)}\otimes_{X_S} (M_{12})_S,
(E_{\infty})_{S} \otimes_{X_{S}} I)\\ \oplus \\ 
\lochom_{X_S}((E_{\infty})_S,E_{S}^{(2)}\otimes_{X_S} I) \\ \end{array}\\ 
& {\buildrel d_{2}\over \longto} 
\lochom_{X_S}(E_S^{(1)}\otimes_{X_S}(M_{12})_S,E_{S}^{(2)}\otimes_{X_{S}} I) \to 0, \\
\end{aligned}
\ee
where
\[
\begin{aligned}
d_{1}(\alpha) = 
{}^t(& -\alpha \circ \Phi^{(1)}_{S,1} + (\Phi^{(2)}_{S,1}\otimes 1_I)
\circ (\alpha\otimes 1_{(M_{1})_S}), \\
&-\alpha \circ \Phi^{(1)}_{S,2}+(\Phi^{(2)}_{S,2}\otimes 1_I)
\circ (\alpha\otimes 1_{(M_{2})_S}), \\
& (\phi_{S}^{(2)}\otimes 1_I) \circ 
(\alpha \otimes 1_{(M_{12})_S}), 
-\alpha\circ \psi_S^{(1)})  \\
\end{aligned}
\]
for any local section $\alpha$ of 
$\lochom_{X_S}(E_S^{(1)},E_{S}^{(2)}\otimes_{X_S}I)$, and 
\[
\begin{aligned}
d_{2}(\beta_1,\beta_2,\gamma, \delta) 
=& \beta_1 \circ (\Phi^{(1)}_{S,2}\otimes 1_{(M_1)_S}) -
(\Phi^{(2)}_{S,2}\otimes 1_I)\circ (\beta_1\otimes 1_{(M_{2})_S})\\ 
& - \beta_2\circ (\Phi^{(1)}_{S,1}\otimes 1_{(M_{2})_S})
+ (\Phi^{(2)}_{S,1}\otimes 1_I)\circ (\beta_2\otimes 1_{(M_{1})_S}) \\
& + 
(\psi_{S}^{(2)}\otimes 1_I)\circ \gamma + \delta \circ \phi_S^{(1)}\\
\end{aligned}
\]
for any local sections $(\beta_1,\beta_2,\gamma, \delta)$ 
of 
\[
\begin{aligned}
& \lochom_{X_S}(E^{(1)}_S\otimes_{X_S}(M_1)_S,E_{S}^{(2)}\otimes_{X_S}I) 
\oplus 
\lochom_{X_S}(E^{(1)}_S\otimes_{X_S}(M_2)_S, E_{S}^{(2)}\otimes_{X_S}I)\\
&   \oplus  
\lochom_{X_S}(E_S^{(1)}\otimes_{X_S}(M_{12})_S 
,(E_{\infty})_{S}\otimes_{X_S}I) \oplus 
\lochom_{X_S}((E_{\infty})_S,E_{S}^{(2)}\otimes_{X_S}I ).\\
\end{aligned}
\]
The degrees of the terms of \eqref{eq:hypercohB} are $0,1,2$
respectively. 

In the following we will use the following notation conventions.
\begin{itemize}
\item 
If $I=\CO_S$, we will denote $\CC(\CE_S^{(1)},\CE_S^{(2)},I)$ by 
$\CC(\CE_S^{(1)},\CE_S^{(2)})$. 
\item 
If $\CE_S^{(1)}=\CE_S^{(2)}=\CE_S$, we will denote 
$\CC(\CE_S^{(1)},\CE_S^{(2)},I)$ by $\CC(\CE_S,I)$. 

\item
If $\CE_S^{(1)}=\CE_S^{(2)}=\CE_S$ and $I=\CO_S$, we will denote 
$\CC(\CE_S^{(1)},\CE_S^{(2)},I)$ by $\CC(\CE_S)$.
\end{itemize}
\end{defi}

\begin{rema}\label{defcomprem}
Under the conditions of Definition (\ref{defcomplex}), suppose 
$S'$ is a nilpotent thickening of $S$ and $I$ is the ideal 
sheaf of $S$ in $S'$. Let $\CE_S$ be a flat family of stable ADHM 
sheaves on $X$ parameterized by $S$, and let $\CE_{S'}$ be 
an extension of $\CE_S$ to $S$. Then the ADHM sheaf $\CE_{S'}$ 
has a natural structure of ADHM sheaf on $X_S$ as explained above Definition 
(\ref{adhmfamB}). Therefore we can take $\CE^{(1)}_S=\CE_S$, 
$\CE^{(2)}_S= \CE_{S'}$ in Definition (\ref{defcomplex}) 
obtaining a complex $\CC(\CE_S, \CE_{S'}, I)$ of 
$\CO_{X_S}$-modules. 
\end{rema}

\begin{prop}\label{infdef}
$(i)$ Let $\CE_S$ be a flat family of stable ADHM sheaves on $X$ parameterized 
by a separated scheme $S$ of finite type over $\IC$, let $S'$ 
be a nilpotent thickening of $S$ and let $\CE_{S'}$ be an extension 
of $\CE_S$ to $S'$. Then $\CE_{S'}$ determines a 
hypercohomology class 
${\mathfrak e}(\CE_S,\CE_{S'},I)\in \IH^1(X_S,\CC(\CE_S, \CE_{S'},I))$.

$(ii)$ Let $\CE_{S'}^{(1)}, \CE_{S'}^{(2)}$ be two extensions
of $\CE_S$ to $S'$. Let $\xi_{S'}: \CE_{S'}^{(1)} \to 
\CE_{S'}^{(2)}$ be a morphism of ADHM sheaves on 
$X_{S'}$ such that $\xi_{S'}|_{X_S} =1_{\CE_S}$. Then there are natural induced 
morphisms of hypercohomology groups 
\be\label{eq:hypercohmorph}
\xi^k_{S'*}: \IH^k(X_S, \CC(\CE_S,\CE_{S'}^{(1)}, I))\to 
 \IH^k(X_S, \CC(\CE_S,\CE_{S'}^{(2)}, I)),
\ee
with $k\in \IZ_{\geq 0}$, such that 
\[
\xi^1_{S'*}({\mathfrak e}(\CE_S,\CE^{(1)}_{S'},I)) =
{\mathfrak e}(\CE_S,\CE^{(2)}_{S'},I) 
\]
\end{prop}

{\it Proof.} 
Suppose $\CE_{S'}=(E_{S'}, \Phi_{S',1,2}, \phi_{S'}, \psi_{S'})$ is such an 
extension. Since $E_{S'}$ is a locally free $\CO_{X_{S'}}$-module
according to Lemma (\ref{locfree}), 
we have an exact sequence 
of $\CO_{X_{S'}}$-modules 
\[
0\to E_{S'} \otimes_{X_{S'}} I \to E_{S'}\to 
E_{S} \to 0 
\]
which yields an exact sequence of $\CO_{X_S}$-modules 
\be\label{eq:extensionA} 
0\to E_{S'} \otimes_{X_{S}} I {\buildrel \iota_S\over \longto} E_{S'}\to 
E_{S} \to 0.
\ee
Since $\CO_{S'}=\CO_S\oplus I$, as $\CO_{S}$-modules, we also 
have a canonical identification
\[
(E_\infty)_{S'}= (E_\infty)_S\oplus (E_\infty)_{S'}\otimes_{X_S} I 
\]
as $\CO_{X_S}$-modules. 
By construction the following  
morphisms of $\CO_{X_S}$-modules are obviously trivial 
\be\label{eq:trivcomp}
\begin{aligned}
& \xymatrix{
E_{S'}\otimes_{X_S} (M_i)_S \otimes_{X_S}I \ar@{^{(}->}[r]
& E_{S'} \otimes_{X_S}
(M_i)_S  \ar[r]^-{ \Phi_{S',i}}   & E_{S'}  \ar@{>>}[r] & E_S}\\
&\xymatrix{
E_{S'}\otimes_{X_S} (M_{12})_{S} \otimes_{X_S} I 
\ar@{^{(}->}[r]& E_{S'}\otimes_{X_S}(M_{12})_{S} 
\ar[r]^-{\phi_{S'}}   & (E_{\infty})_{S'}  \ar@{>>}[r] 
& (E_\infty)_S\\}\\
& \xymatrix{
(E_\infty)_{S'} \otimes_{X_S}I \ar@{^{(}->}[r] 
& (E_\infty)_{S'} \ar[r]^-{\psi_{S'}} & 
E_{S'}\ar@{>>}[r] & E_S,\\}\\
\end{aligned}
\ee
where $i=1,2$. 

The extension of $\CO_{X_S}$-modules \eqref{eq:extensionA} 
determines an extension class in 
$\mathrm{Ext}^1_{X_S}(E_{S}, E_{S'}\otimes_{X_S}I)$.
Since $E_S$ is locally free according to Lemma (\ref{locfree}), we have  
\[
\mathrm{Ext}^1_{X_S}(E_{S}, E_{S'}\otimes_{X_S}I) \simeq H^1(X_S, 
\lochom_{X_S}(E_S, E_{S'}\otimes_{X_S}I)).
\]

Since $S$ is separated and of finite type over $\IC$, 
there exists an affine $\check{\mathrm C}$ech cover 
$\CU_S=\{{U_{S,\alpha}}\}$ of $X_S$. 
Given an open embedding, $U_{S,\alpha_1\ldots\alpha_q} 
\subset U_{S,\alpha_1\ldots\alpha_p}$ with $p<q$, we will 
denote the restriction of any $\CO_{U_{S,\alpha_1\ldots\alpha_p}}$-module 
or morphism of $\CO_{U_{S,\alpha_1\ldots\alpha_p}}$-modules to 
$U_{S,\alpha_1\ldots \alpha_q}$ by a subscript $(\alpha_1\ldots\alpha_q)$. 

By restriction to each affine open subset ${U_{S,\alpha}}\subset X_S$, 
the morphisms 
of $\CO_{X_S}$-modules $\Phi_{S,i}, \phi_S, \psi_S$
yield $\check{\mathrm C}$ech 0-cocycles
$\{(\Phi_{S,i})_\alpha\}, \{(\phi_S)_\alpha\}, \{(\psi_S)_\alpha\}$.

Next note that the restriction of the extension \eqref{eq:extensionA} to any 
affine open subset ${U_{S,\alpha}}$ is trivial since $E_S$ is a locally free 
$\CO_{X_S}$-module.  
Therefore we can choose local splittings 
\be\label{eq:locsplitA}
\xymatrix{
0\ar[r]& (E_{S'}\otimes_{X_S}I)_\alpha \ar[r] 
& (E_{S'})_\alpha \ar@/^/[rr] & & (E_S)_\alpha \ar@/^/[ll]^{\eta_\alpha}
\ar[r] & 0,\\}
\ee
where $\eta=\{\eta_\alpha\} \in C^0(\CU_S, \lochom_{X_S}(E_S, E_{S'}))$. 
Then 
\[
e=\{e_{\alpha\beta}\}=\{(\delta\eta)_{\alpha\beta}\} \in 
C^1(\CU_S, \lochom_{X_S}(E_S, E_{S'}\otimes_{X_S}I))
\]
is a 1-cocycle representing the extension class determined by 
\eqref{eq:extensionA}.

The local splittings $\{\eta_\alpha\}$ determine local isomorphisms 
\[
\eta_\alpha \oplus (\iota_S)_\alpha : E_S \oplus (E_{S'}\otimes_{X_S}I)_\alpha 
{\buildrel \sim \over \longto} (E_{S'})_\alpha 
\]
where $\iota_S: E_{S'}\otimes_{X_S} I \to E_{S'}$ is the canonical 
injection of $\CO_{X_S}$ modules in \eqref{eq:extensionA}. 
Moreover, we  
have commutative diagrams of $\CO_{{U_{S,\alpha}}}$ modules 
of the form 
\be\label{eq:locextA}
\begin{aligned} 
& \xymatrix{ 
(E_S\oplus E_{S'}\otimes_{X_S}I)_\alpha 
\otimes_{{U_{S,\alpha}}} ((M_i)_S)_\alpha 
\ar[rrr]^-{(\eta_\alpha \oplus (\iota_S)_\alpha)\otimes 1_{((M_i)_S)_\alpha}} 
\ar[d]_-{(\Lambda_\alpha)_i} 
& & &(E_{S'})_\alpha\otimes_{{U_{S,\alpha}}} ((M_i)_{S})_\alpha
\ar[d]^-{(\Phi_{S',i})_\alpha} \\
(E_S\oplus E_{S'}\otimes_{X_S}I)_\alpha 
\ar[rrr]^-{\eta_\alpha \oplus (\iota_S)_\alpha} 
& & & (E_{S'})_\alpha\\}\\
& \xymatrix{ 
(E_S\oplus E_{S'}\otimes_{X_S}I)_\alpha \otimes_{{U_{S,\alpha}}} 
((M_{12})_S)_\alpha 
\ar[rrr]^-{(\eta_\alpha \oplus (\iota_S)_\alpha)\otimes 
1_{((M_{12})_S)_\alpha}} 
\ar[d]_-{\lambda_\alpha} 
& & &(E_{S'})_\alpha\otimes_{{U_{S,\alpha}}} ((M_{12})_{S})_\alpha
\ar[d]^-{(\phi_{S'})_\alpha} \\
((E_\infty)_S)_\alpha \oplus  ((E_\infty)_{S'}\otimes_{X_S} I)_\alpha
\ar[rrr]^-{1_{((E_\infty)_S)_\alpha} \oplus 
1_{((E_\infty)_{S'}\otimes_{X_S} I)_\alpha}}
& & & ((E_\infty)_S)_\alpha \oplus  ((E_\infty)_{S'}\otimes_{X_S} I)_\alpha
\\}\\
& \xymatrix{ 
((E_\infty)_S)_\alpha \oplus  ((E_\infty)_{S'}\otimes_{X_S} I)_\alpha
\ar[rrrr]^-{1_{((E_\infty)_S)_\alpha} \oplus 
1_{((E_\infty)_{S'}\otimes_{X_S} I)_\alpha}} \ar[d]_-{\rho_\alpha}
& & & &((E_\infty)_S)_\alpha \oplus  ((E_\infty)_{S'}\otimes_{X_S} I)_\alpha
\ar[d]^-{(\psi_{S'})_\alpha}\\
(E_S\oplus E_{S'}\otimes_{X_S}I)_\alpha
\ar[rrrr]^-{\eta_\alpha \oplus (\iota_S)_\alpha} 
& & & & (E_{S'})_\alpha\\}
\end{aligned}
\ee
for $i=1,2$. Since the morphisms \eqref{eq:trivcomp} are 
trivial, $\{\Lambda_\alpha\}, \{\lambda_\alpha\}, \{\rho_\alpha\}$ 
are 0-cochains of the form 

\be\label{eq:locADHM}
\begin{aligned} 
(\Lambda_\alpha)_i =\left(\begin{array}{cc} (\Phi_{S,i})_\alpha & 0\\
({\wPhi}_{\alpha})_i & (\Phi_{S',i})_\alpha\otimes 1_{I_\alpha} \\
\end{array}\right) & \in C^0(\CU_S, {\mathcal Hom}_{X_S}
((E_S\oplus E_{S'}\otimes_{X_S}I)\otimes_{X_S}(M_i)_S,\\
& \qquad \qquad \qquad \qquad  \ \ E_S\oplus E_{S'}\otimes_{X_S}I))\\
\lambda_\alpha = \left(\begin{array}{cc} (\phi_{S})_\alpha & 0\\
{\wphi}_\alpha & (\phi_{S'})_\alpha \otimes 1_{I_\alpha}\\
\end{array}\right) & \in C^0(\CU_S, \lochom_{X_S}
((E_S\oplus E_{S'}\otimes_{X_S}I)\otimes_{X_S}(M_{12}))_S,\\
& \qquad \qquad \qquad \qquad  \ \ 
(E_\infty)_S\oplus (E_\infty)_{S'}\otimes_{X_S}I) )\\
\rho_\alpha = \left(\begin{array}{cc} (\psi_S)_\alpha & 0 \\
{\wpsi}_\alpha & (\psi_{S'})_\alpha \otimes 1_{I_\alpha}\\
\end{array}\right) & \in C^0(\CU_S, {\mathcal Hom}_{X_S}(
(E_\infty)_S\oplus (E_\infty)_{S'}\otimes_{X_S}I, \\
& \qquad \qquad \qquad \qquad  \ \ E_S\oplus E_{S'}\otimes_{X_S}I)),
\end{aligned}
\ee
with $i=1,2$, 
satisfying the conditions listed below. 
First, we have gluing conditions 
\be\label{eq:cechA} 
\begin{aligned}
(g_{S})_{\alpha\beta} ((\Lambda_\beta)_i)_{\alpha\beta}
& = ((\Lambda_\alpha)_i)_{\alpha\beta} 
(g_{S,i})_{\alpha\beta}\\
((g_\infty)_{S})_{\alpha\beta} 
(\lambda_\beta)_{\alpha\beta} & 
= (\lambda_\alpha)_{\alpha\beta} (g_{S,12})_{\alpha\beta} \\
(g_S)_{\alpha\beta} (\rho_\beta)_{\alpha\beta} 
& = (\rho_\alpha)_{\alpha\beta} ((g_\infty)_S)_{\alpha\beta},\\
\end{aligned}
\ee
where 
\[
\begin{aligned}
(g_{S})_{\alpha\beta} & = \left(\begin{array}{cc} 
1_{(E_S)_{\alpha\beta}} & 0 \\
(e_{S})_{\alpha\beta} & 1_{(E_{S'}\otimes_{X_S}I)_{\alpha\beta}} \\
\end{array} \right) \\
(g_{S,i})_{\alpha\beta} & = \left(\begin{array}{cc} 
1_{(E_S)_{\alpha\beta}}\otimes 1_{((M_i)_S)_{\alpha\beta}} & 0 \\
(e_{S})_{\alpha\beta}\otimes 1_{((M_i)_S)_{\alpha\beta}} 
& 
1_{(E_{S'}\otimes_{X_S}I)_{\alpha\beta}}\otimes 1_{((M_i)_S)_{\alpha\beta}} \\
\end{array} \right) \\
(g_{S,12})_{\alpha\beta} & = \left(\begin{array}{cc} 
1_{(E_S)_{\alpha\beta}}\otimes 1_{((M_{12})_S)_{\alpha\beta}} & 0 \\
(e_{S})_{\alpha\beta}\otimes 1_{((M_{12})_S)_{\alpha\beta}} 
& 1_{(E_{S'}\otimes_{X_S}I)_{\alpha\beta}}\otimes 
1_{((M_{12})_S)_{\alpha\beta}} \\
\end{array} \right) \\
((g_{\infty})_{S})_{\alpha\beta} & = \left(\begin{array}{cc} 
1_{((E_{\infty})_{S})_{\alpha\beta}} & 0 \\
0 & 1_{((E_{\infty})_{S'}\otimes_{X_S}I)_{\alpha\beta}} \\
\end{array} \right) \\
\end{aligned}
\]
for $i=1,2$. 
In addition, the ADHM relation 
\[
\Phi_{S',1} \circ (\Phi_{S',2}\otimes 1_{(M_1)_S}) - 
\Phi_{S',2} \circ (\Phi_{S',1}\otimes 1_{(M_2)_S}) + 
\psi_{S'} \circ \phi_{S'} =0
\] 
is equivalent to the following condition
 \be\label{eq:cechC}
\begin{aligned} 
& \wPhi_{1} (\Phi_{S,2}\otimes 1_{(M_1)_S}) -
(\Phi_{S',2}\otimes 1_I)(\wPhi_{1} \otimes 1_{(M_2)_S}) \\
& -\wPhi_{2} (\Phi_{S,1}
\otimes 1_{(M_2)_S}) 
+ (\Phi_{S',1}\otimes 1_I)(\wPhi_{2} \otimes 1_{(M_1)_S})\\
& + (\psi_{S'}\otimes 1_I)\wphi + \wpsi \phi_S=0.\\
\end{aligned}
\ee
A simple computation shows that the conditions \eqref{eq:cechA} 
are equivalent to 
\be\label{eq:cechB}
\begin{aligned}
\delta(\wPhi_{S,i}) & = -e \Phi_{S,i} + (\Phi_{S',i}\otimes 1_I) (e\otimes 
1_{(M_i)_S}) \\
\delta(\wphi_S) & = (\phi_{S'}\otimes 1_I) 
(e\otimes 1_{(M_{12})_S}) \\
\delta(\wpsi_S) & = -e \psi_S,\\
\end{aligned}
\ee
where $i=1,2$. 

Note that equations \eqref{eq:cechC}, \eqref{eq:cechC} imply that 
the collection of cochains 
$\left(e, \wPhi_{1}, \wPhi_{2}, \wphi, \wpsi\right)$ defines a 
1-cocycle in the total hypercohomology complex associated to the 
complex of $\CO_{S}$-modules \eqref{eq:hypercohB}, i.e. we have 
\[
\delta(e)=0\qquad d_1(e)-\delta(\wPhi_{1,2},\wphi,\wpsi)=0\qquad 
d_2(\wPhi_{1,2},\wphi,\wpsi)=0.
\]

Next suppose we make a different choice of local 
splittings $\{\eta'_\alpha\}$ in \eqref{eq:locsplitA}.
Then note that the difference 
$\{\gamma_\alpha =\eta'_\alpha-\eta_\alpha\}$ 
is a 0-cochain 
\[
\{\gamma_\alpha\} \in C^0(\CU_S, \lochom_{X_S} 
(E_S, E_{S'}\otimes_{X_S}I)).
\]
An elementary computation shows that the cocycle 
$e$ and the cochains $(\wPhi_i, \wphi, \wpsi)$ change as follows 
\[
\begin{aligned}
e'-e & = \delta(\gamma) \\
\wPhi_i' -\wPhi_i , 
& = - \gamma \Phi_{S,i} + (\Phi_{S',i}\otimes 1_I)(\gamma\otimes 
1_{(M_i)_S}) \\
\wphi'-\wphi & = (\phi_{S'}\otimes 1_I)(\gamma\otimes 1_{(M_{12})_S}) \\
\wpsi'-\wpsi & = -\gamma \psi_S,\\
\end{aligned}
\]
for $i=1,2$. Therefore it follows that the hypercohomology 1-cocycle 
$\left(e, \wPhi_{1}, \wPhi_{2}, \wphi, \wpsi\right)$ changes by a 
coboundary in the total hypercohomology complex associated to 
\eqref{eq:hypercohB}, that is 
\[
\begin{aligned}
e'-e & = \delta(\gamma) \\
\left(\wPhi_{1,2},\wphi, \wpsi\right)' - 
\left(\wPhi_{1,2}, \wphi, \wpsi\right) & = 
d_1(\gamma).\\
\end{aligned} 
\]
This proves Proposition (\ref{infdef}.$i$). The required hypercohomology 
class is 
\[
{\mathfrak e}(\CE_S,\CE_{S'},I)=
\left[\left(e, \wPhi_{1}, \wPhi_{2}, \wphi, \wpsi\right)\right].
\]

In order to prove (\ref{infdef}.$ii$)
let $\CE^{(1)}_{S'}$, $\CE^{(2)}_{S'}$ be  
extensions of $\CE_S$ to $S'$, and let 
$\xi_{S'}: \CE^{(1)} \to \CE^{(2)}$ be a morphism 
of ADHM sheaves on $X$ such that $\xi_{S'}|_{S}=1_{\CE_S}$.
Since $\xi_{S'}$ is a morphism of ADHM sheaves, it induces a natural 
morphism of complexes of $\CO_{X_S}$-modules 
\[
\xi_{S'*} : \CC(\CE_S, \CE^{(1)}_{S'}, I) \to \CC(\CE_S, 
\CE^{(2)}_{S'},I) 
\]
which induces in turn the morphisms \eqref{eq:hypercohmorph}. 

It remains to check that $\xi_{S'*}^{1}$ maps the 
class $e(\CE_S,\CE^{(1)}_{S'},I)$ to $e(\CE_S, \CE^{(2)}_{S'},I)$. 
Let $\{\eta^{(1)}_\alpha\}\in C^0(\CU_S, \lochom_{X_S}(E_S, E^{(1)}_{S'}))$  
be local splittings as in \eqref{eq:locsplitA} for 
the extension \eqref{eq:extensionA} with $E_{S'}=E^{(1)}_{S'}$. 
 Since $\xi_{S'}|_S = 
1_{\CE_S}$, the 0-cochain 
$\{\eta^{(2)}_\alpha\}\in C^0(\CU_S, \lochom_{X_S}(E_S, E^{(2)}_{S'}))$, 
\[
\eta^{(2)}_\alpha = (\xi_{S'})_\alpha \circ \eta^{(1)}_\alpha 
\]
provides analogous local splittings for the extension \eqref{eq:extensionA} 
with $E_{S'}= E^{(2)}_{S'}$. 
In particular this implies that 
\[
e^{(2)} = \xi_{S'}\circ e^{(1)}
\]
Moreover we have commutative diagrams of $\CO_{{U_{S,\alpha}}}$ modules 
\be\label{eq:locextB}
\xymatrix{ 
(E_S \oplus E^{(1)}_{S'}\otimes_{X_S}I)_\alpha 
\ar[rrr]^-{\eta_\alpha^{(1)}\oplus \iota^{(1)}_\alpha} 
\ar[d]_-{1_{(E_S)_\alpha} \oplus ((\xi_{S'})_\alpha \otimes 1_{I_\alpha})} 
& & & E_{S'}^{(1)} \ar[d]^-{(\xi_{S'})_\alpha} \\
(E_S \oplus E^{(2)}_{S'}\otimes_{X_S}I)_\alpha 
\ar[rrr]^-{\eta_\alpha^{(2)}\oplus \iota^{(1)}_\alpha}  
& & & E_{S'}^{(2)} \\}
\ee
in which the horizontal arrows are isomorphisms. 
Then, using the commutative diagrams \eqref{eq:locextA}, \eqref{eq:locextB}
and the fact that $\xi_{S'}$ is a morphism of ADHM sheaves, 
a routine computation yields the following relations 
\[
\begin{aligned} 
((\xi_{S'})_\alpha \otimes 1_{I_\alpha})\circ (\wPhi^{(1)}_\alpha)_i  
& = (\wPhi^{(2)}_\alpha)_i\\
\wphi^{(1)}_\alpha 
& = \wphi^{(2)}_\alpha\\
((\xi_{S'})_\alpha \otimes 1_{I_\alpha}) \circ 
\wpsi^{(1)}_\alpha
& = \wpsi^{(2)}_\alpha.\\
\end{aligned}
\]
Obviously, these relations hold for any choice of local splittings 
$\{\eta^{(1)}_\alpha\}$. Since $\xi_{S'}$ is a global morphism 
of ADHM sheaves, this implies that indeed 
$\xi^1_{S'*}({\mathfrak e}(\CE_S, \CE_{S'}^{(1)}, I)) = 
{\mathfrak e}(\CE_S, \CE_S^{(2)}, I)$. 

\hfill $\Box$

\begin{rema}\label{trivthick} 
Under the hypotheses of Proposition (\ref{infdef}.$i$), 
suppose $S'$ is the trivial nilpotent thickening of $S$ determined by 
a coherent 
$\CO_S$-module $I$. Then we have a canonical isomorphism 
$\CC(\CE_S,\CE_{S'},I)\simeq \CC(\CE_S,I)$ 
since $I^2=0$ in the ring structure of $\CO_{S'}$, 
therefore $\CE_{S'}\otimes_{X_S}I\simeq \CE_{S}\otimes_{X_S}I$. 
In particular, in this case 
the complex $\CC(\CE_{S}, \CE_{S'},I)$ is independent of $\CE_{S'}$.
\end{rema} 

A straightforward consequence of Proposition (\ref{infdef}) is the 
following statement analogous to \cite[Lemma 3.4]{RT}. 

\begin{coro}\label{infdefB}
Let $\CE_S$ be a flat family of stable ADHM sheaves on $X$ parameterized 
by an affine scheme $S$ of finite type over $\IC$, and let $S'$ 
be the trivial nilpotent thickening of $S$ determined by a coherent 
$\CO_S$-module $I$. Then there is a one-to-one correspondence 
between equivalence classes of extensions $\CE_{S'}$ of $\CE_S$ to 
$S'$ and hypercohomology classes in $\IH^1(X_S, \CC(\CE_S,I))$. 
\end{coro} 

{\it Proof.} Given the explicit construction of the extension class ${\mathfrak e}(\CE_S, \CE_{S'},I)\in \IH^1(X_S, \CC(\CE_S,I))$ in Proposition (\ref{infdef}.$i$) 
the proof of Corollary (\ref{infdefB}) is straightforward and details will be omitted. 

\hfill $\Box$

Next let us determine the obstructions in the deformation theory of 
stable ADHM sheaves. Recall \cite{artin-moduli}, \cite{LT} that a 
standard deformation situation consists of a 
sequence of closed embeddings of schemes over $S$
\be\label{eq:defsit}
S \subset S' \subset S'',
\ee
where $S$ is a scheme of finite type over $\IC$, and 
$S',S''$ are nilpotent thickenings of $S$. 
We will assume in the following that $S$ is also separated over 
$\IC$.

Let ${I_{S\subset S'}},{I_{S\subset S''}}$, $I_{S'\subset S''}$ 
be  the ideal sheaves corresponding to the closed embeddings 
$S\subset S'$, $S\subset S''$ and $S'\subset S''$ respectively. 
We will assume that  
that ${I_{S\subset S''}}\cdot {I_{S'\subset S''}}=0$ in $\CO_{S''}$, 
hence in particular $I_{S'\subset S''}^2=0$.
Therefore $I_{S'\subset S''}$ has an $\CO_S$-module structure, 
and we have an exact sequence of $\CO_S$-modules 
\be\label{eq:idealext}
0 \to {I_{S'\subset S''}} \to {I_{S\subset S''}} \to {I_{S\subset S'}} \to 0. 
\ee

Suppose we have a flat family $\CE_{S'}$ of stable ADHM sheaves 
 on $X$ parameterized by $S'$, which restricts 
to a given family $\CE_S$ over $S$. 
Recall that $E_{S}$, $E_{S'}$ are locally free $\CO_{X_S}$, 
respectively $\CO_{X_{S'}}$-modules. Moreover, since the projection 
morphism $S'\to S$ is finite, the base change theorem implies that 
$E_{S'}$ is also a locally free $\CO_{X_S}$-module.
Then the exact sequence of $\CO_{S}$-modules \eqref{eq:idealext} 
yields an exact sequence of complexes of $\CO_{X_{S}}$-modules 
\be\label{eq:obsA}
0\to \CC(\CE_{S}, \CE_{S'}, {I_{S'\subset S''}}) \to 
\CC(\CE_{S}, \CE_{S'}, I_{S\subset S''}) 
\to \CC(\CE_{S}, \CE_{S'}, I_{S\subset S'}) \to 0.
\ee
Note that we have a canonical isomorphism \[\CC(\CE_S, \CE_{S'}, 
I_{S'\subset S''}) \simeq \CC(\CE_S,\CE_S, I_{S'\subset S''}) \]
of $\CO_S$-modules since $I_{S'\subset S''}^2=0$. 
Therefore \eqref{eq:obsA} 
yields a long exact sequence of hypercohomology groups 
which reads in part 
\be\label{eq:extexseqC}
\begin{aligned} 
\cdots & \to \IH^1(X_S,\CC(\CE_S,\CE_{S'}, {I_{S\subset S''}}))
 \to \IH^1(X_S, \CC(\CE_S, \CE_{S'}, 
{I_{S\subset S'}}))
\\ & {\buildrel \partial\over \longto} 
\IH^2(X_S,\CC(\CE_S,\CE_S,
{I_{S'\subset S''}}))\to \cdots \\
\end{aligned}
\ee
According to Proposition (\ref{infdef}.$i$), the data $(\CE_S, \CE_{S'},
I_{S\subset S'})$ determines a hypercohomology class 
${\mathfrak e}(\CE_S,\CE_{S'},I_{S\subset S'})\in 
\IH^1(X_S, \CC(\CE_S, \CE_{S'}, I_{S\subset S'})$. 

\begin{defi}\label{obsclass}
Given a deformation situation \eqref{eq:defsit}, a 
flat family $\CE_S$ and an extension $\CE_{S'}$ to $S'$, 
we define the obstruction class 
\[
{\mathfrak {ob}}(\CE_{S'}, S',S'') \in \IH^2(X_S,
\CC(\CE_S, \CE_S, I_{S'\subset S''})) 
\]
to be \[{\mathfrak {ob}}(\CE_{S'}, S',S'') = 
\partial({\mathfrak e}(\CE_S,\CE_{S'},I_{S\subset S'})).\]
\end{defi}
 
\begin{prop}\label{infobs} 
$\CE_{S'}$ can be extended to a flat family of stable ADHM sheaves on 
$X$ parameterized by $S''$ if and only if 
${\mathfrak {ob}}(\CE_{S'},S',S'')=0$.
\end{prop}

{\it Proof.} This statement follows from Proposition (\ref{infdef}) 
by an argument analogous to \cite[Prop. 3.13]{RT} to which we refer the 
reader for more details. 

\hfill $\Box$

The main technical lemma needed in the proof of virtual smoothness is 
the following vanishing result. 

\begin{lemm}\label{hypervanA}
Suppose $X$ is a smooth projective curve over $\IC$, and 
let $\CE$ be a stable ADHM sheaf on $X$. Then 
$\IH^i(X,\CC(\CE))=0$, for all $i\geq 3$ and for all $i\leq 0$. 
\end{lemm} 

{\it Proof.} 
Let us first prove vanishing for $i\geq 3$. $\CC(\CE)$ is obtained 
by setting $S={\mathrm {Spec}}(\IC)$ and $I=\CO_S$ in Definition 
(\ref{defcomplex}). 
Since the degrees of the three terms in \eqref{eq:hypercohB} are 
$0,1,2$ respectively, it follows that 
all terms $E_1^{p,q}$, $p+q\geq 4$ in the standard hypercohomology 
spectral sequence are trivially zero.  Therefore 
$\IH^i(X,\CC(\CE))=0$ for $i\geq 4$. 
Moreover, the only nonzero term on the diagonal $p+q=3$ is 
\[
E_1^{2,1} = H^1(X, \lochom_{X}(E\otimes_X M_{12},E)).
\]
Let 
$d_1^{1,1} : E_1^{1,1} \to E_1^{2,1}$ be the differential
in the first term of the spectral sequence. 
We claim that this map is surjective if $\CE$ is stable. 
It suffices to prove that the dual map 
$({d_1^{1,1}})^\vee : ({E_1^{2,1}})^\vee \to ({E_1^{1,1}})^\vee$
is injective.
Using Serre duality, the dual differential is a linear map 
\[
\begin{aligned}
\ho_X(E,E\otimes M_{12}\otimes_X K_X) 
& \xymatrix{{}\ar[r]^-{({d_1^{1,1}})^\vee} & } 
\begin{array}{c}
\ho_{X}(E,E\otimes_X M_1\otimes_X K_X) \\ \oplus \\
\ho_{X}(E, E\otimes_X M_2 \otimes_X K_X) \\ \oplus \\ 
\ho_{X}(E_\infty,  
E\otimes_X M_{12}\otimes_X K_X)\\ \oplus \\ 
\ho_{X}(E,E_\infty\otimes_X K_X) \\ \end{array}\\
\end{aligned}
\]
which maps a global homomorphism 
$\alpha\in \ho_X(E,E\otimes_X M_{12}\otimes_X K_X)$ 
to 
\[
\begin{array}{c}
{}^t(-(\alpha\otimes 1_{M_2^{-1}}) \circ \Phi_2 + (\Phi_2\otimes 1_{M_1} \otimes 
1_{K_X}) \circ \alpha,\\
(\alpha\otimes 1_{M_1^{-1}}) \circ \Phi_1 - (\Phi_1\otimes 1_{M_2} \otimes 
1_{K_X}) \circ \alpha,\\
\alpha\circ \psi,\
(\phi\otimes 1_{K_X}) \circ \alpha). \\
\end{array}
\]
If $\alpha\in \mathrm{Ker}(({d_1^{1,1}})^\vee)$, it follows that 
$\mathrm{Ker}(\alpha)$ is $\Phi$-invariant and 
$\mathrm{Im}(\psi) \subseteq \mathrm{Ker}(\alpha)$. 
Moreover, since $\alpha$ is a morphism of locally free sheaves, 
$\mathrm{Ker}(\alpha)$ must be saturated subsheaf of $E$. 
Then the stability of $\CE$ implies that $\mathrm{Ker}(\alpha) 
= E$, hence $\alpha=0$. This proves the claim. 

The second assertion is trivial for $i<0$. The case $i=0$ follows 
again from the stability of $\CE$ by an identical argument. 
\hfill $\Box$

\subsection{ADHM Sheaves with Trivial Framing}\label{trivfram}
We conclude this section with some specific results for 
ADHM sheaves with $E_\infty =\CO_X$. Under the conditions of 
Definition (\ref{defcomplex}) note that we have a morphism 
of 
$\CO_{X_S}$-modules 
\[
{\mathcal Hom}_{X_S}((E_\infty)_S, (E_\infty)_S\otimes_{X_S} I) 
{\buildrel \kappa \over\longto}
\begin{array}{c}
{\mathcal Hom}_{X_S}(E_S^{(1)}\otimes_{X_S}(M_{12})_S, 
(E_\infty)_S\otimes_{X_S} I) \\
\oplus\\
{\mathcal Hom}_{X_S}((E_\infty)_S, E_S^{(2)}\otimes_{X_S} I)\\
\end{array}
\]
given by
\[
\alpha_\infty \to {}^t(-\alpha_\infty\circ \phi^{(1)}_S, 
(\psi_S^{(2)}\otimes 1_I)\circ \alpha_\infty)
\]
for any local section $\alpha_\infty$ of 
$
{\mathcal Hom}_{X_S}((E_\infty)_S, (E_\infty)_S\otimes_{X_S} I)$.
Moreover, it is straightforward to check that 
$d_2\circ  \kappa=0$, therefore $\kappa$ yields a morphism 
of complexes 
\be\label{eq:enhancedA}
{\mathcal Hom}_{X_S}((E_\infty)_S, (E_\infty)_S\otimes X_S I)[-1] 
\xymatrix{{}\ar[rr]^-{\kappa[-1]}& & {}} 
\CC(\CE_S^{(1)},\CE_S^{(2)},I).
\ee
Let 
\be\label{eq:enhancedB}
{\widetilde \CC}(\CE_S^{(1)}, \CE_S^{(2)}, I)
={\sf Cone}(\kappa[-1]). 
\ee 
Using notation conventions analogous to Definition 
(\ref{defcomplex}), we will 
denote $\wCC(\CE_S,\CE_S,I)$ by $\wCC(\CE_S,I)$
and $\wCC(\CE_S^{(1)},\CE_S^{(2)},\CO_S)$ by 
$\wCC(\CE_S^{(1)},\CE_S^{(2)})$.
If $I=\CO_S$, $\wCC(\CE_S,I)$ will be denoted by $\wCC(\CE_S)$. 

Note that there is a canonical 0-th hypercohomology class 
$1_{\CE_S}\in \IH^0(X_S,\wCC(\CE_S))$ determined by the pair 
\[
\left(1_{E_S}, 1_{(E_\infty)_S}\right) \in 
H^0(X_S, \lochom_{X_S}(E_S,E_S))\oplus 
H^0(X_S, \lochom_{X_S}((E_\infty)_S,(E_\infty)_S))
\]
It is straightforward to check that 
\[
{\widetilde d}_1\left(1_{E_S}, 1_{(E_\infty)_S}\right)=0
\]
therefore $1_{\CE_S}$ is indeed a well defined 0-th hypercohomology class.

By construction we also have an exact sequence of complexes of 
$\CO_{X_S}$-modules 
\be\label{eq:exseqdef}
0\to \CC(\CE_S^{(1)}, \CE_S^{(2)}, I) \to 
{\widetilde \CC}(\CE_S^{(1)}, \CE_S^{(2)}, I) \to 
\lochom_{X_S}((E_\infty)_S, (E_\infty)_S\otimes_{X_S} I)\to 0.
\ee

\begin{lemm}\label{splitlemma}
Under the conditions of Definition (\ref{defcomplex}),
suppose $\CE_S^{(1)}=\CE_S^{(2)}=\CE_S$ and $E_\infty =\CO_X$. 
Then the exact sequence \eqref{eq:exseqdef} has a canonical splitting, 
and we have isomorphisms of hypercohomology groups 
\be\label{eq:splithypercoh}
\IH^k(X_S, {\widetilde \CC}(\CE_S, \CE_S, I)) 
\simeq \IH^k(X_S, \CC(\CE_S, \CE_S, I)) \oplus 
H^k(X_S, p_S^*I)
\ee
for all $k\in \IZ$.
\end{lemm}

{\it Proof.} 
Since $(E_\infty)_S=\CO_{X_S}$, 
$\lochom_{X_S}((E_\infty)_S, (E_\infty)_S\otimes_{X_S} I)\simeq p_S^*I$, 
and we have a morphism 
\[
\begin{aligned}
p_S^*I & \longto \begin{array}{c}
{\mathcal Hom}_{X_S}(E_S, E_S\otimes_{X_S} I) 
\oplus
p_S^*I\\
\end{array}\\
\alpha_\infty & \longmapsto \left(1_{E_S}\otimes \alpha_\infty, \alpha_\infty\right)
\end{aligned}
\]
It is straightforward to check that this is a splitting of
\eqref{eq:exseqdef}.

\hfill $\Box$

Under the conditions of Proposition (\ref{infdef}), 
we have an exact sequence of complexes of $\CO_{X_S}$-modules 
\be\label{eq:exseqcompl}
0\to {\widetilde \CC}(\CE_S, \CE_{S'}, I) \to 
{\widetilde \CC}(\CE_S,\CE_{S'},\CO_{X_S}) \to {\widetilde \CC}(\CE_S)
\to 0.
\ee
Let 
\be\label{eq:enhestclas}
{\widetilde {\mathfrak e}}(\CE_S,\CE_{S'},I)
={\widetilde \partial}^1(1_{\CE_S}),
\ee
where 
\[
{\widetilde \partial}^1 : \IH^0(X_S,{\widetilde \CC}(\CE_S)) \to \IH^1(X_S,
{\widetilde \CC}(\CE_S,\CE_{S'},I))
\]
is the connecting isomorphism determined by \eqref{eq:exseqcompl}. 
We will also denote by 
\[
{\widetilde i}^1 : \IH^1(X_S, \CC(\CE_S,\CE_{S'},I)) \to 
\IH^1(X_S,\wCC(\CE_S,\CE_{S'},I))
\]
the natural morphism of hypercohomology groups determined by 
the canonical injective morphism of complexes 
$\CC(\CE_S,\CE_{S'},I)\hookrightarrow 
\wCC(\CE_S,\CE_{S'},I)$ and by 
\[
{\widetilde i}^2: \IH^2(X_S, \CC(\CE_S,\CE_S,I))\hookrightarrow 
\IH^2(X_S, \wCC(\CE_S,\CE_S,I)) 
\]
the injection determined by the splitting \eqref{eq:splithypercoh}. 

Now consider a deformation situation of the form \eqref{eq:defsit}, 
and let $\CE_{S'}$ be a flat family of stable ADHM sheaves 
parameterized by $S'$ which restricts to a given family 
$\CE_S$ on $X_S$. Then 
we have an exact sequence of complexes of $\CO_{X_S}$-modules 
\be\label{eq:enhexseq}
0\to \wCC(\CE_S, I_{S'\subset S''}) 
\to \wCC(\CE_S, \CE_{S'},I_{S\subset S''})
\to \wCC(\CE_S, \CE_{S'},I_{S\subset S'})
\to 0.
\ee

\begin{lemm}\label{trivfrobs}
Still assuming $E_\infty =\CO_X$, we have 
\be\label{eq:enhobsclass}
{\widetilde i}^2({\mathfrak {ob}}(\CE_{S'}, S',S'')) = ({\widetilde i}^1 
\circ {\widetilde \partial}^2)({\widetilde {\mathfrak e}}(\CE_{S}, \CE_{S'},
I_{S\subset S'})),
\ee 
where 
\[
{\widetilde \partial}^2 : \IH^1(X_S, \wCC(\CE_S,\CE_{S'}, I_{S\subset S'}))
\to \IH^2(X_S, \wCC(\CE_S,\CE_S,I_{S'\subset S''}))
\]
is the connecting morphism determined by \eqref{eq:enhexseq}.
\end{lemm}

{\it Proof.} This follows from a standard hypercohomology computation sing the 
detailed construction of the class ${\mathfrak e}(\CE_S,\CE_{S'},I_{S\subset S'})$ 
given in the proof of Proposition (\ref{infdef}). 

\hfill $\Box$

\section{Virtual Smoothness for ADHM Sheaves on Curves}\label{virsmoothsect}

From this point on we take $X$ to be a smooth projective curve of genus 
$g$ over $\IC$.
Then the Hilbert polynomial of any torsion-free $\CO_X$-module is determined 
by a pair of integers $(r,e)$ with $r\geq 1$, namely the rank 
and degree respectively. In the following we will denote the moduli 
space $\mfm_{ADHM}(\CX,P)$ by $\mfm_{ADHM}(\CX,r,e)$.  

In this section our goal 
is to prove that the moduli spaces ${\mfm}_{ADHM}(\CX,r,e)$ 
are virtually smooth, i.e. they carry natural 
perfect obstruction theories \cite{LT,BF}. 
Similar results have been 
previously obtained for moduli spaces of (decorated) sheaves 
in \cite{RT,OT,Minv,stabpairs-I}. 
Our treatment is closer to  \cite{RT}, relying on the construction 
of \cite{LT}, or, more precisely, 
the generalization presented in \cite[Sect 2.2]{degGW}. 
Moreover we will carry out our construction  
in the equivariant setting as in \cite{GP} with respect to 
an algebraic torus action on $\mfm_{ADHM}(\CX,r,e)$ satisfying 
the hypothesis of Lemma 
(\ref{Tactioncoro}.)
The results of this section are valid for any 
such action, in particular for those presented in 
Examples (\ref{torusactA}), (\ref{torusactB}). 

\subsection{Outline}\label{outline}
This subsection summarizes the main elements in the construction of a tangent-obstruction 
theory of a moduli stack following \cite[Sect. 1]{LT}, \cite[Sect. 1.2, Sect. 2.1]{degGW}.  Let $\mfm$ be a separated Deligne-Mumford stack of finite type over $\IC$. 
One first defines the perfect tangent-obstruction theory of an affine \'etale chart $\iota:S \to 
\mfm$ and then specifies a compatibility condition on overlaps. 
Let $\CE_S$ denote the object of $\mfm(S)$ corresponding to the morphism $\iota$
and ${\mathfrak {Mod}}_S$ denote the category of quasi-coherent $\CO_S$-modules. 
For any object $I$ of ${\mathfrak {Mod}}_S$ let $S_I$ denote the trivial nilpotent thickening of $S$ by $I$ 
defined above Definition (\ref{adhmfamB}).
The tangent-obstruction theory of the affine chart $\iota:S\to \mfm$ consists of two elements. 

The first element is the deformation functor ${\mathcal T}{\mfm}(\CE_S): {\mathfrak {Mod}}_S\to Sets$ which 
assigns to a $\CO_S$-module $I$ the set of equivalence classes of extensions 
of $\CE_S$ to $S_I$. Here two extensions of $\CE_S$ to $S_I$ are called equivalent if 
they are isomorphic as objects of $\mfm(S_I)$ and the isomorphism restricts to the 
identity $1_{\CE_S}$ on $S\subset S_I$. 
An important requirement is  that 
for any triple $(S,\CE_S,I)$ there exist a functor 
${\mathcal T}^1{\mfm}(\CE_S):{\mathfrak {Mod}}_S \to 
{\mathfrak {Mod}}_S $ 
such that ${\mathcal T}{\mfm}(\CE_S)(I)$ is 
canonically isomorphic 
to the space of global sections $\Gamma_S({\mathcal T}^1{\mfm}(\CE_S)(I))$
for any $I$. 
Moreover, ${\mathcal T}^1{\mfm}(\CE_S)$ 
is required to satisfy a natural base change property. 
The existence of this functor for ADHM sheaves follows from Proposition (\ref{infdef}), 
as stated in 
 Lemma (\ref{tangspace}) below. 
 
The next element is an obstruction theory for $\iota:S \to \mfm$. 
Suppose $S\subset S'\subset S''$ is a deformation situation
 as in \eqref{eq:defsit} and $\iota':S'\to \mfm$ is a morphism corresponding to an 
 object $\CE_{S'}$. Let $\CE_S = \CE_{S'}|_{S}$. 
 Then an obstruction theory \cite[Def 1.2]{LT}, \cite[Def. 1.11]{degGW},
  consists of a $\CO_S$-module $\CO b_{\CE_S}$ (depending 
 on $\CE_S$), 
 and an obstruction class 
 $${\mathfrak {ob}}'(\CE_{S'},S',S'')\in \Gamma_S(\CO b_{\CE_{S}}\otimes_S 
 \CI_{S'\subset S''})$$
 such that the vanishing of ${\mathfrak {ob}}'(\CE_{S'},S',S'')$ is a 
 necessary and sufficient condition  for the existence of an extension of $\CE_{S'}$ to $S''$.
 Moreover, ${\mathfrak {ob}}'(\CE_{S'},S',S'')$ is required to satisfy a base change property.
 Again, the existence of an obstruction theory for ADHM sheaves follows from Proposition
 (\ref{infobs}), as stated in Lemma (\ref{obsspace}). In particular the class 
 ${\mathfrak {ob}}'(\CE_{S'},S',S'')$ is naturally identified with the obstruction class 
  ${\mathfrak {ob}}(\CE_{S'},S',S'')$ introduced in Definition (\ref{obsclass}). 
 
Next note that given an obstruction theory of $\iota:S\to \mfm$, the $\CO_S$-module 
$\CO b_{\CE_S}$ determines a functor ${\mathcal T}^2{\mfm}(\CE_S):
{\mathfrak {Mod}}_S \to 
{\mathfrak {Mod}}_S$, ${\mathcal T}^2{\mfm}(\CE_S)(I) = \CO b_{\CE_S}\otimes_S I$. 
The tangent-obstruction complex of the chart $\iota:S\to \mfm$ is the complex 
${\mathcal T}^\bullet\mfm(\CE_S)=[{\mathcal T}^1\mfm(\CE_S)\to 
{\mathcal T}^2\mfm(\CE_S)]$ where the arrow is the zero map \cite[Def. 1.3]{LT}. 
Furthermore, the tangent-obstruction 
complex is called perfect if 
there exists a two-term 
complex of locally free $\CO_{S}$-modules $\IE_S^\bullet$ 
with the property that ${\mathcal T}^i\mfm(\CE_{S})(I)$, $i=1,2$, is the 
$i$-th cohomology sheaf of $\IE_S^\bullet\otimes_{S_\alpha} I$ 
for any $\CO_{S}$-module $I$. 
This property is checked in Lemma (\ref{perfect}) for the tangent-obstruction theory
of ADHM sheaves. As explained in more detail below the differential $\IE_S^1\to \IE_S^2$ 
is related to the Kuranishi map of the deformation theory \cite[Sect. 2]{LT}.

The above data defines a perfect tangent-obstruction theory for an affine 
\'etale chart. The gluing condition is formulated as follows 
\cite[Def. 1.2]{LT}, \cite[Def. 2.1]{degGW}. Suppose 
$\{\iota_\alpha :S_\alpha\to \mfm\}$ is a 
an \'etale affine cover of $\mfm$ consisting of finitely many charts. Suppose 
furthermore that each chart is equipped with a perfect tangent obstruction theory
with obstruction sheaf $\CO b_\alpha$ and obstruction assignment 
${\mathfrak {ob}}_\alpha$. Then the collection ${\CO b}_\alpha$ is required to 
descend to a sheaf on $\mfm$ and the obstruction assignments are required 
to agree on overlaps. More details will be provided in subsection (\ref{virsmooth}).

Finally note that a perfect obstruction theory determines a virtual cycle in the Chow 
group of the moduli space. First, the construction in \cite[Sect 3]{LT} based on relative 
Kuranishi maps produces a virtual normal cone $C_\alpha$ in the total space 
of the vector bundle $\mathrm{V}(\IE^2_\alpha)$ associated to 
the locally free sheaf $\IE_\alpha^2$. Very briefly, a relative Kuranishi map 
for the perfect tangent-obstruction complex $\IE_\alpha^\bullet$ is essentially a morphism 
$\kappa_\alpha : \hat {\mathrm{V}}(\IE_\alpha^1)\to \hat {\mathrm{V}}(\IE_\alpha^2)$ between 
the formal completions of $V(\IE^1_\alpha)$, $V(\IE^2_\alpha)$ along the zero sections, 
such that the inverse image $\kappa_\alpha^{-1}(0)$ is isomorphic to the formal completion $Z_\alpha$
of the diagonal in $S_\alpha \times S_\alpha$. It is also required to satisfy certain 
compatibility conditions with the obstruction classes which will not be reviewed here. 
More details can be found in \cite[Sect 2]{LT}. 

Given such a map let $\Gamma^{\kappa_\alpha}$ denote its graph and let 
${\mathcal N}^{\kappa_\alpha}$ denote the 
normal cone to be the scheme theoretic intersection 
\[
\Gamma^{\kappa_\alpha} \times_{\hat {\mathrm{V}}(\IE_\alpha^1)\times \hat {\mathrm{V}}(\IE_\alpha^2)} \hat {\mathrm{V}}(\IE_\alpha^1)
\]
where $ \hat {\mathrm{V}}(\IE_\alpha^1)$ is canonically embedded in 
$\hat {\mathrm{V}}(\IE_\alpha^1)\times \hat {\mathrm{V}}(\IE_\alpha^2)$ 
by the zero section. Note that ${\mathcal N}^{\kappa_\alpha}$ is a closed subcone 
of $Z_\alpha \times {\mathrm{V}}(\IE_\alpha^2)$, which is the total 
space of the normal bundle to ${\mathrm{V}}(\IE_\alpha^1)$ in 
${\hat {\mathrm{V}}(\IE_\alpha^1)\times \hat {\mathrm{V}}(\IE_\alpha^2)}$.
The cone $C_\alpha$ is the restriction of ${\mathcal N}^{\kappa_\alpha}$ to $S_\alpha\subset \kappa_\alpha^{-1}(0)$, that is 
$C_\alpha = {\mathcal N}^{\kappa_\alpha} 
\times_Z S_\alpha$. 

In principle applying the refined Gysin map 
$0_\alpha^! : A_*(\mathrm{V}(\IE_\alpha^2)) \to 
A_{*-\mathrm{rk}(\IE_\alpha^2)}(S_\alpha)$ to the virtual normal cone yields 
a cycle in the Chow group $A_*(S_\alpha)$. 
The main question is then how the resulting local cycles 
can be glued in order to obtain a cycle in $A_*(\mfm)$. An answer to this question is provided in \cite[Sect 2.2, Lemm. (2.6)-(2.8)]{degGW}, the final conclusion 
being that there exists indeed a well defined virtual cycle $[\mfm]^{vir}$. Omitting the 
details of the proof, which is long and very technical, note that the key element is the 
existence of a global obstruction sheaf $\CO b$ over $\mfm$ such that the obstruction 
classes are compatible on overlaps. It is also worth noting that the gluing algorithm 
is significantly simpler if the moduli stack $\mfm$ is a scheme, $\{S_\alpha\}$ is a 
Zariski open cover, and there exists a global 
two term complex of locally free sheaves $\IE^\bullet$ on $\mfm$ which restricts to 
$\IE_\alpha^\bullet$ on each $S_\alpha$. Then the above cone construction can be 
carried out globally over $\mfm$, producing the virtual cycle \cite[Sect 3]{LT}. 
Finally, note that although 
the local Kuranishi maps are not canonically determined by the perfect tangent-obstruction theory, the virtual cycle is independent of the choices involved in this construction. 

To summarize, in order to construct a virtual cycle it suffices to establish the 
existence of a perfect tangent-obstruction theory in each affine chart $S_\alpha \to 
\mfm$, then check that the there is a global obstruction sheaf, and the obstruction classes 
satisfy a gluing condition. This will be done in the next two subsections 
for moduli of ADHM sheaves on curves. 

\subsection{The Tangent-Obstruction Complex}\label{adhmtangobs}

The following lemmas establish the existence of a tangent-obstruction 
theory for any flat family of stable ADHM sheaves on $X$ parameterized 
by an affine scheme $S$ of finite type over $\IC$. The proofs consist 
of a straightforward verification of the conditions 
formulated in \cite[Sect 1]{LT} using Propositions (\ref{infdef}), (\ref{infobs}), 
Lemma (\ref{hypervanA}). Details will be omitted since similar arguments 
have been presented in the proof of \cite[Prop. 3.26]{RT} and \cite[Thm. 3.28]{RT}.

\begin{lemm}\label{tangspace} 
Let $\CE_S$ be a flat family of stable ADHM sheaves on $X$ 
parameterized by an affine 
scheme $S$ of finite type over $\IC$. Let $I$ be a coherent 
$\CO_S$-module. Then the following hold. 

$(i)$ In the notation of Remark (\ref{trivthick}), we have 
\be\label{eq:tangspaceA}
\IH^1(X_S, \CC(\CE_S,I)) 
\simeq \Gamma_S({\bf R}^1p_{S*}\CC(\CE_S,I)).
\ee

$(ii)$ Given a base change morphism  
$f:T\to S$ with $T$ affine of finite 
type over $\IC$, a coherent
$\CO_T$-module $J$ and a 
morphism of $\CO_{T}$-modules 
$\xi:f^*I\to J$, we have a canonical morphism 
\be\label{eq:tangspaceB}
c(f,\xi):f^*{\bf R}^1p_{S*} \CC(\CE_S,I)\to 
{\bf R}^1p_{T*}\CC(\CE_T,J)
\ee
of $\CO_{T}$-modules, where $\CE_T=f_X^*\CE_S$. 

$(iii)$ Given a sequence of  morphisms $U{\buildrel g\over \longto}T 
{\buildrel f\over \longto}S$, with $U,T$ affine of finite type 
over $\IC$, a coherent $\CO_U$-module $K$, 
a coherent $\CO_T$-module $J$, a coherent 
$\CO_S$-module $I$, and morphisms 
$\eta:g^*J\to K$, 
$\xi:f^*I\to J$, there is a commutative diagram of canonical 
morphisms of the form 
\be\label{eq:requireddiag}
\xymatrix{ 
g^*f^*{\bf R}p_{S*}\CC(\CE_S,I)
\ar[r]^-{\simeq} \ar[d]_-{g^*c(f,\xi)} & 
(f\circ g)^* {\bf R}^1p_{S*}\CC(\CE_S,I)
\ar[d]^-{c(f\circ g,\zeta)} \\
g^*{\bf R}^1p_{T*}\CC(\CE_T,J)
\ar[r]^-{c(g,\eta)} &
{\bf R}^1p_{U*}\CC(\CE_U,K),\\
}
\ee
where 
$\CE_T=f_X^*\CE_S$, $\CE_U=g_X^*\CE_T$ and 
$\zeta: (f\circ g)^*I \to K$ is defined by the composition 
\[
(f\circ g)^*I \simeq 
g^*f^*I {\buildrel g^*\xi\over \longto} g^*J {\buildrel \eta\over \longto} K.
\]

$(iv)$ 
Suppose there is a 
${\bf T}$-action on $S$ such that the family $\CE_S$ is ${\bf T}$-equivariant. 
Then ${\bf R}^1p_{S*}(\CE_S,\CE_S\otimes_{X_S} I)$
has a natural ${\bf T}$-equivariant structure for any ${\bf T}$-equivariant 
$\CO_S$-module $I$. Moreover, if all base change morphisms 
as well as the sheaves $I,J,K$ are {\bf T}-equivariant 
the canonical morphism $c(f,\xi)$ is {\bf T}-equivariant, 
and the diagram \eqref{eq:requireddiag} is ${\bf T}$-equivariant. 
 
\end{lemm}

\begin{lemm}\label{obsspace}
Let $\CE_S$ be a flat family of stable ADHM sheaves on $X$ 
parameterized by an affine scheme $S$ of finite type over 
$\IC$. Let 
$S\subset S'\subset S''$ be a deformation situation as 
in \eqref{eq:defsit}. 
Then the following hold. 

$(i)$ In the notation of Remark (\ref{trivthick}), we have 
\be\label{eq:obsspaceO}
\IH^2(X_S,\CC(\CE_S,I_{S' \subset S''})) \simeq 
\Gamma_S({\bf R}^2 p_{S*}\CC(\CE_S)\otimes_{\CO_S}I_{S'\subset S''}).
\ee

$(ii)$ Given a base change morphism $f:T\to S$ 
with $T$ affine of finite type over $\IC$
we have a canonical isomorphism 
\be\label{eq:canbasechange}
c(f):f^*{\bf R}^2p_{S*}\CC(\CE_S) {\buildrel \sim \over \longto} 
{\bf R}^2p_{T*}\CC(\CE_T) 
\ee
compatible with base change, where $\CE_T=f_X^*\CE_S$. 

$(iii)$ Suppose $S\subset S'\subset S''$ 
is a deformation situation of the form \eqref{eq:defsit}
with $S,S',S''$ affine of finite type over $\IC$.
Suppose moreover we have a base change diagram 
\[
\xymatrix{ 
T'' \ar[r]^{f''} \ar[d] & S''\ar[d]^{} \\
T \ar[r]^-{f} & S \\
}
\]
where $T$ is an affine scheme of finite type over $\IC$
and $T''\to T$ has a section 
$\sigma_T:T\to T''$ such that the  diagram 
\[
\xymatrix{
T'' \ar[r]^{f''} & S'' \\
T\ar[u]^{\sigma_T} \ar[r]^f & S \ar[u]^{\sigma_S} \\
}
\]
is commutative. Let $T'= S'\times_{S''} T''$, 
$f':T'\to S'$ be the natural projection, and $\CE_T = f_X^*\CE_S$.
Suppose furthermore we have an extension $\CE_{S'}$ of $\CE_S$ to 
$S'$, and let $\CE_{T'} = (f'_X)^* \CE_{S'}$. 
Then there is a canonical morphism of $\CO_{T}$-modules 
\be\label{eq:canmorphA}
c(f,S',S''): f^*{\bf R}^2p_{S*}(\CC(\CE_S)\otimes_S 
I_{S'\subset S''})
\to {\bf R}^2p_{T*}(\CC(\CE_T)\otimes_T I_{T'\subset T''})
\ee
such that 
\be\label{eq:obsspaceQ}
{\mathfrak {ob}}(\CE_{T'},T',T'') = (\Gamma_T(c(f,S,S''))\circ f^*)
{\mathfrak {ob}}(\CE_{S'},S',S'').
\ee

$(iv)$ Suppose there is a ${\bf T}$-action on $S$ such that the family 
$\CE_S$ is ${\bf T}$-equivariant. 
Then ${\bf R}^2p_{S*}(\CC(\CE_S,\CE_S\otimes_{X_S} I))$
has a natural ${\bf T}$-equivariant structure for any ${\bf T}$-equivariant 
$\CO_S$-module $I$. Moreover, statements $(i)$, $(ii)$ 
hold for ${\bf T}$-equivariant base change, and the canonical 
morphism \eqref{eq:canmorphA} is ${\bf T}$-equivariant if the deformation 
situations $S\subset S'\subset S''$, $T\subset T'\subset T''$ are
{\bf T}-equivariant. 
 
\end{lemm}

\hfill $\Box$

Next we prove that the tangent-obstruction theory of any affine 
chart $S\to \mfm_{ADHM}(\CX,r,e)$ of finite type over $\IC$
is perfect, verifying the conditions of 
\cite[Def. 1.3]{LT}.  We will need the following preliminary lemma. 
\begin{lemm}\label{EGAlemmaA}
Let $\rho:A\to B$ be a local morphism of Noetherian local rings, 
$k$ the residual field of $A$ and $M,N$ two finitely generated 
$B$-modules, $N$ flat over $A$. Let $u:M\to N$ a morphism of $B$-modules.
Then the following conditions are equivalent
\begin{itemize} 
\item[$(i)$] $u$ is injective and $\mathrm{Coker}(u)$ is flat over $A$ 
\item[$(ii)$] $u\otimes 1:M\otimes_A k \to 
N\otimes_Ak$ is injective. 
\end{itemize}
\end{lemm}

{\it Proof.} \cite[10.2.4]{EGAIIIa}.

\begin{lemm}\label{perfect}
Let $\CE_S$ be a flat family of stable ADHM sheaves parameterized by 
an affine scheme $S$ of finite type over $\IC$. 
Then there exists a two-term complex  
$\IE^\bullet = \left(\IE_{S}^1 {\buildrel {}\over \longto} \IE_{S}^2
\right)$ of coherent locally free $\CO_S$-modules with degrees $(1,2)$ such that 
\[
\begin{aligned}
\CH^1(\IE_S^\bullet\otimes_S I) 
& \simeq {\bf R}^1p_{S*}\CC(\CE_S,I) \\
\CH^2(\IE_S^\bullet\otimes_S I) 
& \simeq {\bf R}^2p_{S*}\CC(\CE_S)\otimes_{S} I \\
\end{aligned}
\]
for any coherent $\CO_S$-module $I$.

Moreover, if $S\to \mfm_{ADHM}(\CX,r,e)$ is a ${\bf T}$-equivariant chart, 
and $I$ is a {\bf T}-equivariant coherent $\CO_S$-module, 
the complex $\IE_{S}^\bullet$ can be chosen ${\bf T}$-equivariant as 
well.   
\end{lemm}

{\it Proof.} 
Note that  
$\CC(\CE_S,I)=\CC(\CE_S)\otimes_{X_S}I$ for any 
coherent $\CO_S$-module $I$.
According to \cite[Thm. 6.10.5]{EGAIIIb} or \cite[Prop. III.12.2]{har}
there exists a 
complex $\IF_{S}^\bullet$ of finitely generated locally free 
$\CO_S$-modules bounded above such that 
\[
\CH^k(\IF_{S}^\bullet\otimes_S I) = 
{\bf R}^kp_{S*}(\CC(\CE_S)\otimes_{X_S}I)
\]
for any $\CO_S$-module $I$ and for all $k\in \IZ$. 
Lemma (\ref{hypervanA}) and the base change theorem
\cite[Thm 7.7.5]{EGAIIIb}
\cite[Thm III.12.11]{har} 
imply that 
\[
{\bf R}^kp_{S*}(\CC(\CE_S,\CE_S))=0
\]
for all $k\geq 3$ and 
\[
{\bf R}^2p_{S*} (\CC(\CE_S)\otimes_{X_S}I) \simeq 
{\bf R}^2p_{S*} \CC(\CE_S)\otimes_S I 
\]
for any coherent $\CO_S$-module $I$. 
 
Since ${\IF}_{S}^\bullet$ is bounded above and has trivial 
cohomology in degrees $k\geq 3$, an easy induction argument 
based on \cite[Prop. III.9.1A.f]{har} shows that 
${\mathrm{Ker}}(\IF_{S}^2\to \IF_{S}^3)$ 
is a flat finitely generated $\CO_S$-module. 
Therefore  ${\mathrm{Ker}}(\IF_{S}^2\to \IF_{S}^3)$ is 
a locally free $\CO_S$-module, and we can truncate $\IF_{S}^\bullet$ 
such that $\IF_{S}^k=0$ for $k\geq 3$.

Moreover, 
theorem \cite[Thm. 6.10.5]{EGAIIIb} also implies that for any 
closed point $s\in S$, the cohomology of the complex of 
vector spaces $\IF_{S}^\bullet \otimes k(s)$ is isomorphic to the 
hypercohomology of the complex $\CC(\CE_S)|_{X_s}$. 
Since the later vanishes in degrees $k\leq 0$ according to Lemma 
(\ref{hypervanA}), it follows that the cohomology of 
$\IF_{S}^\bullet \otimes k(s)$ is trivial in degrees $k\leq 0$ for any 
closed point $s\in S$. Then Lemma (\ref{EGAlemmaA}) implies that we can 
truncate $\IF_S^\bullet$ to 
a locally free complex $\IE^\bullet_S$ satisfying the required properties. 

Now suppose $S\to \mfm_{ADHM}(\CX,r,e)$ is ${\bf T}$-equivariant, 
hence in particular 
there is a ${\bf T}$-action on $S$ such that $\CE_S$ has a ${\bf T}$-equivariant 
structure. Then all above arguments hold in the equivariant setting, 
as observed for example in 
\cite[Sect 2.2]{equivres}.    

\hfill $\Box$

In conclusion, Lemmas (\ref{tangspace}), (\ref{obsspace}), 
(\ref{perfect}) imply 
\begin{prop}\label{perftangobs} 
Let $\iota:S\to \mfm_{ADHM}(\CX,r,e)$ 
be a {\bf T}-equivariant affine  chart of the moduli space 
of stable ADHM sheaves on $X$. Let $\CE_S$ be the flat family of stable ADHM sheaves 
parameterized by $S$. Let $\CC(\CE_S,I) = \CC(\CE_S,\CE_S,I)$, 
$\CC(\CE_S)=\CC(\CE_S,\CO_S)$ be the deformation complex in 
Definition (\ref{defcomplex}) and  $\CC(\CE_S)=\CC(\CE_S,\CO_S)$. 
Then $S$ is equipped with a 
{\bf T}-equivariant perfect obstruction theory 
$[{\mathcal T}^1( \mfm_{ADHM}(\CX,r,e))(\CE_S)\to 
{\mathcal T}^2( \mfm_{ADHM}(\CX,r,e))(\CE_S)]$ 
where the arrow is the zero map, and
\[
\bal 
{\mathcal T}^1( \mfm_{ADHM}(\CX,r,e))(\CE_S)(I) & = {\bf R}^1p_{S*}\CC(\CE_S,I)\\
{\mathcal T}^2( \mfm_{ADHM}(\CX,r,e))(\CE_S)(I) &  = 
\CO b_{\CE_S}\otimes_S I
\\
\CO b_{\CE_S}  & = {\bf R}^2p_{S*}\CC(\CE_S)\\
\eal 
\]
for any $\CO_S$-module $I$.
Moreover, there is an isomorphism $ \Gamma_S(\CO b_{\CE_S})\simeq 
\IH^2(X_S, \CC(\CE_S))$ mapping the 
the obstruction 
class ${\mathfrak{ob}}'(\CE_{S'},S',S'')$  to the obstruction 
class ${\mathfrak{ob}}(\CE_{S'},S',S'')$ determined in Proposition (\ref{infobs}).
\end{prop}

\subsection{Virtual Smoothness}\label{virsmooth}
In this section we prove Theorem (\ref{virsmooththm}).

{\it Proof of Theorem (\ref{virsmooththm}).}
Given a torus action as in Lemma (\ref{Tactioncoro}) there exists a 
{\bf T}-equivariant affine Zariski open  cover $\{\iota_\alpha : S_\alpha \to \mfm_{ADHM}(\CX,r,e)\}$ of 
$\mfm_{ADHM}(\CX,r,e)$  
with $S_\alpha$ affine schemes of finite type over $\IC$.  
According to Proposition (\ref{perftangobs}), each chart $(S_\alpha, \iota_\alpha)$
carries a {\bf T}-equivariant perfect 
tangent-obstruction theory.
Therefore as explained in section (\ref{outline}), we just have to check the compatibility conditions 
formulated in \cite[Def. 2.1]{degGW}. Namely we have to check that there is a 
global obstruction sheaf $\CO b$ on $\mfm$ which restricts to the local obstruction 
sheaves on each $S_\alpha$, and that the obstruction classes agree on overlaps. 

According to Theorem (\ref{DMstack}), $\mfm_{ADHM}(\CX,r,e)$ is a fine quasi-projective 
moduli scheme equipped with a universal stable ADHM sheaf ${\mathfrak E}$ on 
$\mfm_{ADHM}(\CX,r,e)\times X$, flat over $\mfm_{ADHM}(\CX,r,e)$. 
Let $\CC(\mfe)$ is the deformation complex of the 
universal object defined in (\ref{defcomplex}) and 
${\mathfrak p} :  \mfm_{ADHM}(X,r,e)\times X \to \mfm_{ADHM}(X,r,e)$ 
be  the canonical projection. By construction, it is obvious that the restriction of $\CC(\mfe)$
to each affine open set $S_\alpha$ is isomorphic to the deformation complex $\CC(\mfe|_{S_\alpha})$. Let $\CO b = {\bf R}^2{\mathfrak p}_*\CC(\mfe)$. Then for each $S_\alpha$, 
the local obstruction sheaf $\CO b_{\CE_{S_\alpha}}$ given in Proposition (\ref{perftangobs}) is isomorphic to $\CO b|_{S_\alpha}$ since $S_\alpha$ is a 
Zariski open subset of the moduli space. 
Therefore $\CO b$ is indeed the global obstruction sheaf. 
One has also to  check that the obstruction assignments agree 
on overlaps. This follows by straightforward manipulations from 
Lemma (\ref{obsspace}) because the obstruction classes satisfy the base change 
property. 

Finally, note that there exists a global locally-free two term complex $\IE^\bullet$
on $\mfm$ which restricts to the complexes 
 $\IE^\bullet_{S_\alpha}$ on each Zariski open subset $S_\alpha$. 
The argument is identical to the proof of Lemma (\ref{perfect}), 
using \cite[Prop. 2.1.10]{huylehn} to establish the existence of a finite locally free 
complex $\IF$ on  $\mfm_{ADHM}(\CX,r,e)\times X$ analogous to $\IF_{S}$. 
Therefore, as explained at the end of section (\ref{outline}), in this case the virtual cycle is 
determined by a global cone over the moduli space. 

\hfill $\Box$

\section{Admissible Pairs}\label{admpairs}

Let $X$ be a smooth projective curve over a field $k$
over $\IC$, and $M_1,M_2$ 
fixed line bundles on $X$. 
Let $Y$ be the total space of the projective 
bundle $\mathrm{Proj}(\CO_X \oplus M_1\oplus M_2)$. 
In this section we introduce new objects -- called admissible pairs on $Y$
-- and prove that they are equivalent to stable pairs on $Y$ as defined in 
\cite{stabpairs-I}. In the next section we will show that 
admissible pairs are naturally related to stable ADHM sheaves 
on $X$ by the relative Beilinson spectral sequence. 

Let us first fix some notation. Let $\pi:Y\to X$ 
denote the canonical projection and $k_i=\mathrm{deg}(M_i)$, 
$i=1,2$.  Recall that \cite[Prop. II.7.1]{har} 
\be\label{eq:projdirimA}
\pi_*(\CO_Y(1)) = \CO_X \oplus M_1\oplus M_2.
\ee
Let $z_0\in H^0(Y,\CO_Y(1))$ denote the section corresponding to 
$1\in H^0(X,\CO_X)$ under the canonical isomorphism 
\[
H^0(Y,\CO_Y(1))\simeq H^0(X,\CO_X) \oplus H^0(X,M_1) 
\oplus H^0(X, M_2).\]
We have analogous canonical sections
\[z_i\in H^0(Y, \pi^*M_i^{-1}(1)),\] with $i=1,2$.

Let $D_\infty$ be the zero locus of $z_0$ on $Y$, which will be referred to as 
the divisor at infinity. Let $H\in A_2(Y)$ denote the cycle class 
of $D_\infty$ and let $F\in A_2(Y)$ denote the fiber class. 
Let $\sigma:X \to Y$ be the section determined by $z_1=z_2=0$, and 
let $\beta=\sigma_*[X]\in A_1(Y)$.
For future reference let us record the following easy lemma which follows from 
elementary intersection computations. 
\begin{lemm}\label{usefulstuff}
$(i)$ $\mathrm{Ker}(A_1(Y) {\buildrel  \cap H\over \longto} A_0(Y))
=\langle \beta\rangle$
 
$(ii)$ $\beta=H^2-(k_1+k_2)HF$ in the intersection ring of $Y$. 

$(iii)$ $c_1(Y) = 3H +(2-2g-k_1-k_2)F$.

$(iv)$ $\mathrm{Pic}(Y)\simeq \mathrm{Pic}(X) \times \IZ$ where 
the second direct summand is generated by the divisor class
$[D_\infty]$. 
\end{lemm}

\begin{defi}\label{admissible}
$(i)$ Let $d\in \IZ_{\geq 1}$, $n\in \IZ$. 
An admissible pair of type $(d,n)$ on $Y$ is a pair 
$(Q,\rho)$ consisting of a coherent $\CO_Y$-module 
$Q$ and a section $\rho\in H^0(Y,Q)$ satisfying the following conditions 
\begin{itemize} 
\item[$(a)$] $\rho$ is not identically zero.
\item[$(b)$] $Q$ is flat over $X$.
\item[$(c)$] $\ch_0(Q)=0$, $\ch_1(Q)=0$, $\ch_2(Q)=d\beta$  and 
$\chi(Q)=n$. 
\item[$(d)$] The cokernel $\mathrm{Coker}(\CO_Y {\buildrel \rho
\over \longto }Q)$ of the canonical morphism determined by $\rho$ 
is an $\CO_Y$-module of dimension 0.
\end{itemize}

$(ii)$ Two admissible pairs of type $(d,n)$ $(Q,\rho)$, $(Q',\rho')$ are 
isomorphic if there exists an isomorphism of 
$\CO_Y$-modules $u:Q\to Q'$ such that the following diagram is 
commutative
\[
\xymatrix{
\CO_Y \ar[r]^-{\rho} \ar[d]_-{1_{\CO_Y}} & Q \ar[d]^u \\
\CO_Y  \ar[r]^-{\rho'}& Q' \\
}
\]
\end{defi} 

\begin{rema}\label{admrem}
In the present context, a stable pair of type $(d,n)$ 
on $Y$ in the sense of 
\cite{stabpairs-I} is a pair $(Q,\rho)$ where $Q$ is a
coherent $\CO_{Y}$-module of pure dimension one 
with $\mathrm{ch}_2(Q)=d\beta$ and $\chi(Q)=n$
and $\rho\in H^0(Y,Q)$ is a nonzero section such that 
$\mathrm{Coker}(\CO_Y{\buildrel \rho\over \longto}
Q)$ is a $\CO_Y$-module of pure dimension zero. 
The following lemma shows that we can replace condition 
(\ref{admissible}.$i.b$) by 
\begin{itemize}
\item[$(b')$] $Q$ is of pure dimension one and $\mathrm{supp}(Q)$ is 
disjoint from $D_\infty$.
\end{itemize}
Therefore there is a one-to-one correspondence between admissible 
and  stable pairs of type $(d,n)$ on $Y$ whose supports are disjoint from $D_\infty$. 
\end{rema}

\begin{lemm}\label{flatness} 
Let $S$ be an affine scheme of finite type over $\IC$ and let $(Q_S, \rho_S)$ be a pair 
consisting of a coherent 
$\CO_{Y_S}$-module, flat over $S$, 
and a section $\rho_S \in H^0(Y_S,Q_S)$
such that the restriction $(Q_s,Y_s)=(Q_S|_{Y_s}, \rho_S|_{Y_s})$, 
$s\in Y$, satisfies 
conditions $(i.a)$, $(i.c)$ and $(i.d)$ of 
Definition (\ref{admissible}). Then $Q_s$ is flat over $X_s$ if and only if 
it is of pure dimension one and $\mathrm{supp}(Q_s)$ is disjoint from 
$(D_\infty)_s\subset Y_s$ for any $s\in Y$. 
\end{lemm}

{\it Proof.} 
The direct implication follows immediately from the local criterion of flatness 
\cite[Thm. 6.8]{eisenbud}
observing that the conditions of Definition (\ref{admissible}) and 
Lemma (\ref{usefulstuff}) rule out vertical components of $\mathrm{supp}(Q_s)$. 

Conversely, let $s\in S$ and 
suppose $Q_s$ is of pure dimension one and its 
support has no vertical components. Let $y\in Y_s$ be a closed 
point such that $(Q_s)_y\neq 0$, $x=\pi_s(y)$ and $\zeta \in \CO_{X_s,x}$ 
a uniformizing parameter, where $\pi_s:Y_s\to X_s$ is the 
natural projection.  Let $Z$ be the scheme theoretic support of $Q_s$ and
$\CI_Z$ be the ideal sheaf of $Z$ in $Y_s$. 
Since $Z$ is one dimensional, $\CI_{Z,y}$ has codimension 2 in 
$\CO_{Y_s,y}$. Moreover, since $Z$ has no vertical components, 
it follows that the ideal $\CJ\subset \CO_{Y_s,y}$ generated by 
$\zeta$ and $\CI_{Z,y}$  
is codimension 3 in $\CO_{Y_s,y}$. However note that 
$\CJ$ annihilates $\mathrm{Ker}((Q_s)_y{\buildrel 
\zeta\over \longto } (Q_s)_y)$. Therefore if $\mathrm{Ker}((Q_s)_y{\buildrel 
\zeta\over \longto } (Q_s)_y)$ is nontrivial, it must be a
dimension zero submodule of $Q_s$. This contradicts the
assumption that $Q_s$ is pure dimension $1$.

In conclusion, $\mathrm{Ker}((Q_s)_y{\buildrel 
\zeta\over \longto } (Q_s)_y)$ 
must be trivial for any $y\in Y_s$ such that $(Q_s)_y\neq 0$. 
Hence $Q_s$ is flat over $X_s$ by the local criterion of 
flatness. 

\hfill $\Box$

{\it Proof of Theorem (\ref{admstack})}
According to  
\cite{LePoitierA,LePoitierB,stabpairs-I} 
there exists a projective moduli space 
$\mfm_{St}(Y,d,n)$ parameterizing isomorphism 
classes of stable pairs of type $(d,n)$ on $Y$.
Since stable pairs have trivial stabilizers \cite[Sect 2.2]{stabpairs-I}, 
this is a fine moduli space. 
Moreover, according to \cite[Thm. 2.15]{stabpairs-I} the moduli space 
$\mfm_{St}(Y,d,n)$ has a natural perfect tangent-obstruction theory, 
therefore a virtual cycle $[\mfm_{St}(Y,d,n)]$. 

Lemma (\ref{flatness}) 
identifies $\mfm_{Adm}(Y,d,n)$ with an open subscheme of 
$\mfm_{St}(Y,d,n)$, therefore $\mfm_{Adm}(Y,d,n)$ is a 
quasi-projective scheme over $\IC$ equipped with an induced 
perfect tangent-obstruction theory.
Although the construction presented in \cite{stabpairs-I} is based on 
the formalism of \cite{BF}, it is fairly straightforward to 
provide an alternative construction of a perfect tangent-obstruction 
theory of $\mfm_{St}(Y,d,n)$ relying on \cite[Sect. 2.1]{degGW}.
Then \cite[Prop. 3]{BKP}
implies that both constructions yield identical virtual cycles. 

For future reference, let us recall the construction of tangent and 
respectively obstruction spaces, 
according to \cite[Sect 3]{RT}, \cite[Sect 2]{stabpairs-I}.
Given a flat family $(Q_S,\rho_S)$ of stable pairs, 
let $C(Q_S,\rho_S)$ denote the two term complex 
of $\CO_{Y_S}$-modules 
\be\label{eq:twotermO} 
\CO_{Y_S} {\buildrel \rho_S\over \longto} Q_S,
\ee
where the degrees are $0,1$ respectively. 
Let also $q_S:Y_S\to S$ denote the canonical projection.
Suppose that $S$ is affine of finite type over $\IC$, and 
let $S'$ be a trivial nilpotent extension of $S$ determined 
by a coherent $\CO_S$-module $I$. Then the tangent space to a
flat family $(Q_S,\rho_S)$ is the traceless
ext group 
\[
{\mathrm{Ext}}_{Y_S}^1(C(Q_S,\rho_S), C(Q_S,\rho_S)\otimes_{Y_S}I)_0
\simeq \Gamma_{S}
({\mathcal Ext}^1_{q_S}(C(Q_S,\rho_S), C(Q_S,\rho_S)\otimes_{Y_S}I)_0)
\]
Given a deformation situation of the form \eqref{eq:defsit}, 
a flat family $(Q_S, \rho_S)$ and an extension $(Q_{S'}, \rho_{S'})$, 
the obstruction space is the traceless
ext group 
\[
\begin{aligned}
{\mathrm{Ext}}_{Y_S}^2(C(Q_S,\rho_S), & C(Q_S,\rho_S)
\otimes_{Y_S}
I_{S'\subset S''})_0
\simeq \\ & \Gamma_{S}
({\mathcal Ext}^2_{q_S}(C(Q_S,\rho_S), C(Q_S,\rho_S))_0
\otimes_{S}
I_{S'\subset S''}).\\
\end{aligned}
\]
The obstruction class 
\[
{\mathfrak {ob}}(C(Q_{S'},\rho_{S'}),S,S'')
\in {\mathrm{Ext}}_{Y_S}^2(C(Q_S,\rho_S),C(Q_S,\rho_S)\otimes_{Y_S}
I_{S'\subset S''})_0
\]
is given by 
\[
{\mathfrak {ob}}(C(Q_{S'},\rho_{S'}),S,S'') ={\delta}^2 
({\mathfrak e}(C(Q_S,\rho_S),C(Q_{S'},\rho_{S'})), I_{S\subset S'})
\]
where 
\[
\delta^2:\mathrm{Ext}^1_{Y_S}(C(Q_S,\rho_S), C(Q{S'},\rho_{S'})\otimes_{Y_S}
I_{S\subset S'}) \to \mathrm{Ext}^2_{Y_S}(C(Q_S,\rho_S), 
C(Q_S,\rho_S)\otimes_{Y_S}I_{S'\subset S''}) 
\]
is a natural coboundary morphism, 
and 
\[
{\mathfrak e}(C(Q_S,\rho_S),C(Q_{S'},\rho_{S'}), I_{S\subset S'}) 
\in {\mathrm{Ext}}^1_{Y_S}(C(Q_S,\rho_S),C(Q_{S'},\rho_{S'}) \otimes_{Y_S}
I_{S\subset S'})
\]
is the extension class corresponding to $C(Q_{S'},\rho_{S'})$. 
The later is given in turn by 
\[
{\mathfrak e}(C(Q_S,\rho_S),C(Q_{S'},\rho_{S'}), I_{S\subset S'}) 
={\delta}^1(1_{C(Q_S,\rho_S)})
\]
where 
\[
\delta^1: \mathrm{Ext}^0_{Y_S}(C(Q_S,\rho_S), C(Q_S,\rho_S)) \to 
{\mathrm{Ext}}^1_{Y_S}(C(Q_S,\rho_S),C(Q_{S'},\rho_{S'}) \otimes_{Y_S}
I_{S\subset S'})
\]
is again a natural coboundary morphism. 

Note that 
$\delta^2({\mathfrak e}(C(Q_S,\rho_S),C(Q_{S'},\rho_{S'}), I_{S\subset S'})$
belongs to the traceless ext group 
${\mathrm{Ext}}_{Y_S}^2(C(Q_S,\rho_S),C(Q_S,\rho_S)\otimes_{Y_S}
I_{S'\subset S''})_0$, which is canonically identified with a subgroup 
of ${\mathrm{Ext}}_{Y_S}^2(C(Q_S,\rho_S),C(Q_S,\rho_S)\otimes_{Y_S}
I_{S'\subset S''})$ (see \cite[Thm. 3.23, Thm 3.28]{RT}.) 

Finally note that there is a ${\bf T}=\IC^\times \times \IC^\times$ 
action on $\mfm_{St}(Y,d,n)$ induced by the natural
fiberwise scaling action on $M_1\oplus M_2$.
Obviously this action preserves the open subspace $\mfm_{Adm}(Y,d,n)$. 
Moreover, the above
construction easily generalizes to the equivariant setting, yielding 
a  {\bf T}-equivariant perfect
tangent-obstruction theory on $\mfm_{Adm}(Y,d,n)$.

Note also that 
the fixed locus $\mfm_{Adm}(Y,d,n)^{\bf T}$ is naturally identified 
with a closed subset of the fixed locus $\mfm_{St}(Y,d,n)^{\bf T}$, 
which is proper over $\IC$. 
More precisely, it is identified with the closed subset of the fixed locus 
$\mfm_{St}(Y,d,n)^{\bf T}$ parameterizing {\bf T}-fixed 
stable pairs $(Q,\rho)$ on 
$Y$ such that the support of $Q$ contains no vertical components. 
Therefore $\mfm_{Adm}(Y,d,n)^{\bf T}$ is 
proper over $\IC$ as well, and has a virtual cycle as well as an 
equivariant virtual normal bundle,
determined by the perfect tangent-obstruction theory of $\mfm_{Adm}(Y,d,n)$. 

\hfill $\Box$

\section{ADHM Sheaves via 
Relative Beilinson Spectral Sequence}\label{adhmadmsect}

Our next goal is to prove that there is an isomorphism 
between the moduli space of stable ADHM sheaves on $X$ 
with data $\CX=(M_1,M_2,\CO_X)$ and the moduli space of admissible 
pairs on $Y$ using a relative Beilinson monad construction. 
Moreover will prove that this isomorphism is compatible 
with the natural ${\bf T}=\IC^\times \times \IC^\times$ actions
on these spaces and identifies the equivariant perfect obstruction 
theories. 

Recall that under the current notation conventions
we  set $X_S = S\times X$, $Y_S=S\times Y$ for any scheme $S$ over $\IC$, and 
let $p_X:X_S\to X$, $p_Y:Y_S \to Y$, 
$\pi_S: Y_S \to X_S$ denote the canonical projections. 
We will also set $F_S = p_X^*F$, $G_S= p_Y^*G$ for 
any $\CO_X$-module $F$, respectively $\CO_Y$-module $G$. 
Note that $Y_S$ is a projective scheme over $X_S$ and we have 
canonical isomorphisms 
$\CO_{Y_S}(1) \simeq 
p_Y^*\CO_{Y}(1)$, $\Omega^k_{Y_S/X_S}\simeq p_Y^*\Omega^k_{Y/X}$, 
$k=0,1,2$. Moreover, by base change we also have a canonical isomorphism 
\[
\pi_{S*}\CO_{Y_S}(1) \simeq \CO_{X_S} \oplus (M_1)_S \oplus (M_2)_S 
\]
Let 
\[
z_{S,0}\in H^0(Y_S,\CO_{Y_S}(1)),\qquad z_{S,i}\in H^0(Y_S,\pi_S^* 
(M_i)_S^{-1}(1)),
\] 
$i=1,2$, be canonical sections defined by analogy with 
$z_0,z_1,z_2$. 
Note that $z_{S,0}$, $z_{S,i}$, are canonically identified with 
$p_Y^*z_0$, $p_Y^*z_i$, respectively, for $i=1,2$.

\subsection{Relative Beilinson Monad for Admissible Pairs}

\begin{lemm}\label{relbeilinsonA}
Let $S$ be a scheme of finite type over $\IC$ 
and let $Q_S$ be an $\CO_{Y_S}$-module, flat over $S$
such that $Q_S|_{Y_s}$ satisfies conditions 
$(i.b)$, $(i.c)$ of Definition (\ref{admissible}).
Then there is a spectral sequence of $\CO_{Y_S}$-modules 
with 
\be\label{eq:specseqA}
\begin{aligned} 
E_1^{i,j} = 
\CO_{Y_S}(j) \otimes_{Y_S}\pi_S^* 
R^i\pi_{S*}(\Omega_{Y_S/X_S}^{-j}(-j)\otimes_{Y_S} 
Q_S), 
\end{aligned}
\ee
$i\geq 0$, $j=-2,-1,0$,
converging to $Q_S$ if $i=j=0$ and $0$ otherwise.  
\end{lemm}

{\it Proof.} This is the relative Beilinson spectral sequence for 
$Y_S/X_S$  \cite{orlov}.

\hfill $\Box$

Our next goal is to compute the terms in \eqref{eq:specseqA}.
We will need the following lemma.

\begin{lemm}\label{EGAlemmaB}
Let $\rho:A\to B$, $\sigma:B\to C$ be local morphisms of Noetherian local rings, 
$k$ the residual field of $A$ and $M$ 
a finitely generated $C$-module. 
Suppose that $B$ is flat over $A$. Then the following conditions are
equivalent. 
\begin{itemize}
\item[$(i)$] $M$ is flat over $B$ 
\item[$(ii)$] $M$ is flat over $A$ and $M\otimes_Ak$ is 
flat over $B\otimes_A k$. 
\end{itemize}
\end{lemm}

{\it Proof.} \cite[10.2.5]{EGAIIIa}

\hfill $\Box$

Given a flat family of admissible pairs $(Q_S,\rho_S)$ parameterized by 
$S$, let $C(Q_S,\rho_S)$ denote the two term complex of 
$\CO_{Y_S}$-modules 
\be\label{eq:twotermA}
\xymatrix{
\CO_{Y_S}\ar[r]^-{\rho_S} & 
Q_S\\}
\ee 
with degrees $(0,1)$. 

\begin{lemm}\label{relbeilinsonB}
Let $(Q_S,\rho_S)$ be a flat family of admissible pairs parameterized 
by $S$ and $C(Q_S,\rho_S)$ the complex defined in \eqref{eq:twotermA}. 
Then $C(Q_S,\rho_S)$ admits a canonical three term locally free resolution 
supported in degrees $-1,0,1$ of the form 
\be\label{eq:threetermA} 
\begin{aligned} 
0 & \to \CO_{Y_S}(-2)\otimes_{Y_S} 
\pi_S^*(E_S\otimes_{X_S} (M_{12})_S) \\
& {\buildrel \sigma_S\over \longto}  
\CO_{Y_S}(-1)\otimes_{Y_S} 
\pi_S^*(E_S\otimes_{X_S} (M_1)_S\oplus E_S\otimes_{X_S} 
(M_2)_S)\oplus \CO_{Y_S}  {\buildrel \tau_S\over \longto} 
\pi_S^*E_S \to 0, \\
\end{aligned}
\ee
where $E_S=R^0\pi_{S*}Q_S$ and 
\[
\begin{aligned}
\sigma_S & = {}^t(z_{S,0}\otimes \pi_S^*\Phi_{S,2}-z_{S,2},\, 
-z_{S,0}\otimes \pi_S^*\Phi_{S,1}+z_{S,1},\,0) \\
\tau_S & =(z_{S,0}\otimes \pi_S^*\Phi_{S,1}-z_{S,1},\, 
z_{S,0}\otimes \pi_S^*\Phi_{S,2}-z_{S,2},\, 
\pi_S^*\psi_S)\\
\end{aligned}
\]
for certain morphisms of $\CO_{X_S}$-modules $
\Phi_{S,i}:E_S\otimes_{X_S} (M_i)_S \to E_S$, $i=1,2$, 
$\psi_S: \CO_{X_S} \to E_S $ satisfying 
\be\label{eq:commrelA}
\Phi_{S,1}\circ (\Phi_{S,2}\otimes 1_{(M_1)_S}) - \Phi_{S,2}\circ 
(\Phi_{S,1}\otimes 1_{(M_2)_S}) =0,
\ee
and 
$\psi_S|_{X_s}: \CO_{X_s} \to E_S|_{X_s}$ is injective  
for each point $s\in S$. 
\end{lemm}

{\it Proof.} 
Let us first compute the terms in \eqref{eq:specseqA}. 
Since $Q_S$ is flat over $S$ and $Q_S|_{X_s}$ is flat over 
$X$ for any point $s\in S$, Lemma (\ref{EGAlemmaB}) 
implies that $Q_S$ is flat over $X_S$. Since it is also 
finite over $X_S$ according to (\ref{admrem}.$i$) the base change 
theorem \cite[Thm 7.7.5]{EGAIIIb}, 
\cite[Thm. III.3.4]{Banica} implies that 
\be\label{eq:dirimcompB}
R^i\pi_{S*}(\Omega_{Y_S/X_S}^{-j}(-j)\otimes_{Y_S} Q_S) =0 
\ee
for all $i\geq 1$ and all $j=-2,-1,0$. 
Moreover $R^0\pi_{S*}(\Omega_{Y_S/X_S}^{-j}(-j)\otimes_{Y_S} Q_S)$ 
are locally free $\CO_{X_S}$-modules for all $j=-2,-1,-0$ 
\cite[Cor. 7.9.10]{EGAIIIb}.

Then the column $i=0$ of the spectral sequence \eqref{eq:specseqA} 
reads 
\be\label{eq:specseqB}
\xymatrix{ 
0 \ar[d]^-{d_1^{0,-2}}  \\
\CO_{Y_S}(-2)\otimes_{Y_S} \pi_S^*R^0\pi_{S*}
(\Omega_{Y_S/X_S}^{2}(2)\otimes_{Y_S} Q_S)\ar[d]^-{d_1^{0,-1}} \\
\CO_{Y_S}(-1)\otimes_{Y_S} \pi_S^*R^0\pi_{S*}
(\Omega_{Y_S/X_S}^{1}(1)\otimes_{Y_S} Q_S)\ar[d]^-{d_1^{0,0}} \\
\pi_S^*R^0\pi_{S*}Q_S\ar[d]^-{d_1^{0,1}}\\ 
0 & 0, \\}
\ee
where we have used the notation $d_1^{i,j}:E_1^{i,j}\to E_1^{i,j+1}$, 
$i=0,1$, $j=-2,-1,0,1$ for differentials. All other columns are trivial. 

Note that the 
relative Euler sequence for $\pi_S:Y_S\to X_S$
\be\label{eq:relEulerB}
0 \to \Omega^1_{Y_S/X_S}(1) \to \pi_S^*(\CO_{X_S} \oplus 
(M_1 \oplus M_2)_S) 
{\buildrel \epsilon_S \over \longto} \CO_{Y_S}(1)\to 0
\ee
yields an exact sequence of $\CO_{X_S}$-modules of the form 
\[
\begin{aligned} 
0 & \to R^0\pi_{S*}(\Omega^1_{Y_S/X_S}(1)\otimes_{Y_S} Q_S)\\
& \to 
R^0\pi_{S*}(\pi_S^*(\CO_{X_S} \oplus (M_1 \oplus M_2)_S)
\otimes_{Y_S}Q_Sx) {\buildrel f_S\over \longto} R^0\pi_{S*}(Q_S(1))\to 0.\\
\end{aligned}
\] 
The morphism $f_S$ is of the form 
\[ 
f_S(t_0,t_1,t_2)= \varPsi_{S,0}(t_1)+ \varPsi_{S,1}(t_1)+\varPsi_{S,2}(t_2)
\]
where 
\be\label{eq:diffmap}
\begin{aligned}
\varPsi_{S,0}: R^0\pi_{S*}Q_S&\to R^0\pi_{S*}(Q_S(1))\\ 
\varPsi_{S,1}:R^0\pi_{S*}Q_S
\otimes_{X_S} (M_1)_S&\to 
R^0\pi_{S*}(Q_S(1))\\
\varPsi_{S,2}: 
R^0\pi_{S*}Q_S\otimes_{X_S} (M_2)_S &\to R^0\pi_{S*}(Q_S(1))\\
\end{aligned}
\ee
are induced by multiplication by $z_{S,0},z_{S,1},z_{S,2}$ respectively.
Note that multiplication by $z_{S,0}$ induces an isomorphism 
$Q_S \simeq  Q_S(1)$ 
since the support of $Q_S$ is disjoint from 
$S\times D_\infty\subset Y_S$. 
Therefore $\varPsi_{S,0}:R^0\pi_{S*}Q_S \to 
R^0\pi_{S*}(Q_S(1))$ is an isomorphism.

Then it follows that the morphism 
\[
\begin{aligned}
(M_1\oplus M_2)_S\otimes_{X_S} R^0\pi_{S*}Q_S& \to 
(\CO_{X_S} \oplus (M_1\oplus M_2)_S)\otimes_{X_S} R^0\pi_{S*}Q_S
\\
(t_1,t_2) & \to (-\varPsi_{S,0}^{-1}\circ \varPsi_{S,1}(t_1)-
\varPsi_{S,0}^{-1}
\circ \varPsi_{S,2}(t_2), t_1,t_2)
\end{aligned}
\]
maps 
$(M_1\oplus M_2)_S\otimes_{X_S} R^0\pi_{S*}Q_S$
isomorphically to $\mathrm{Ker}(f_S)$.
Therefore we obtain an identification 
\be\label{eq:identifA}
 R^0\pi_{S*}(\Omega^1_{Y_S/X_S}(1)\otimes_{Y_S} Q_S) \simeq 
(M_1\oplus M_2)_S\otimes_{X_S} R^0\pi_{S*}Q_S.
\ee
Note also that $\Omega^2_{Y_S/X_S}(2) \simeq 
\pi_S^*(M_{12})_S(-1)$. Therefore  \eqref{eq:specseqB} 
reduces to a complex of $\CO_{Y_S}$-modules of the form 
\be\label{eq:specseqC} 
\begin{aligned} 
0 & \to \CO_{Y_S}(-2) \otimes_{Y_S} \pi_S^*(R^0\pi_{S*}(Q_S(-1)) 
\otimes_{X_S} (M_{12})_S) \\
& {\buildrel d_1^{0,-1}\over 
\longto} \CO_{Y_S}(-1)\otimes_{Y_S} \pi_S^*(R^0\pi_{S*}Q_S 
\otimes_{X_S}(M_1\oplus M_2)_S)\\
& {\buildrel d_1^{0,0}\over 
\longto}  \pi_S^*R^0\pi_{S*}
Q_S \to 0. \\
\end{aligned}
\ee
The differentials $d_1^{0,-1}$, $d_1^{0,0}$ are determined by 
the morphisms $\kappa_{1S}, \kappa_{2S}$ in the Koszul complex 
\be\label{eq:koszulB} 
0 \to \CO_{Y_S}\boxtimes \Omega^2_{Y_S/X_S}(2) 
{\buildrel {\kappa_{1S}} \over \longto} \CO_{Y_S}(-1) 
\boxtimes \Omega^1_{Y_S/X_S}(1) {\buildrel {\kappa_{2S}}\over 
\longto} \CO_{Y_S} \boxtimes \CO_{Y_S}.
\ee 
The later are determined in turn by contraction 
with the section $\Sigma_S \in H^0(Y_S\times_{X_S} Y_S, 
\CO_{Y_S}(1)\boxtimes 
T_{Y_S/X_S}(-1))$ corresponding to the identity under the canonical 
identification 
\[
H^0(Y_S\times_{X_S} Y_S, \CO_{Y_S}(1)\boxtimes 
T_{Y_S/X_S}(-1)) \simeq {\mathrm{End}}_{X_S}(\CO_{X_S} 
\oplus (M_1^{-1} \oplus M_2^{-1})_S).
\]
Explicit expressions for $\kappa_{1S},\kappa_{2S}$
can be written using the isomorphisms
\[
\begin{aligned}
&\Omega^2_{Y_S/X_S}(2) \simeq \pi_S^*(M_{12})_S(-1), \\
&\Omega^1_{Y_S/X_S}(1) \simeq \mathrm{Ker}\left(\pi_S^*(\CO_{X_S} 
\oplus (M_1\oplus M_2)_S) {\buildrel \epsilon_S\over \longto} 
\CO_{Y_S}(1)\right).\\
\end{aligned}
\]
Then we have
\be\label{eq:koszulD}
\begin{aligned} 
\kappa_{1S}(u) =
(& z_{S,1}\boxtimes z_{S,2}-z_{S,2}\boxtimes z_{S,1}, \\
& -z_{S,1}\boxtimes z_{S,0} + z_{S,0}\boxtimes z_{S,1},
z_{S,2}\boxtimes z_{S,0} - z_{S,0}\boxtimes z_{S,2}),\\
\end{aligned}
\ee
where $z_{S,i}$, $i=1,2$, are identified with global sections of 
$\CO_{Y_S}(1)\boxtimes \pi_S^*(M_i^{-1})_S$,
$i=1,2$ via the obvious isomorphisms
\[
(\CO_{Y_S}(1)\otimes_{Y_S} \pi_S^*(M_i^{-1})_S 
\boxtimes \CO_{Y_S} \simeq 
\CO_{Y_S}(1)\boxtimes (\pi_S^*(M_i^{-1})_S).
\]
A priori, the morphism \eqref{eq:koszulD} takes values in 
$\CO_{Y_S}(-1)\boxtimes \pi_S^*(\CO_{X_S} \oplus (M_1)_S 
\oplus (M_2)_S)$. 
However one can easily check that 
\[
(1\boxtimes \epsilon_S)(\kappa_{1S}(u)) =0 
\]
for any local section $u$, hence $\kappa_{1S}$ takes values 
in $\CO_{Y_S}(1)\boxtimes \Omega^1_{Y_S/X_S}(1)$ as expected. 
The second morphism is locally given by 
\be\label{eq:koszulE} 
\kappa_{2S}(v_0,v_1,v_2) = p_{1S}^*(z_{S,0})v_0+p_{1S}^*(z_{S,1})v_1
+p_{1S}^*(z_{S,2})v_2. 
\ee
Using the identification \eqref{eq:identifA} and 
expressions \eqref{eq:koszulD}, \eqref{eq:koszulE}
the complex \eqref{eq:specseqC} becomes isomorphic to 
\be\label{eq:specseqCA}
\begin{aligned} 
0 & \to \CO_{Y_S}(-2) \otimes_{Y_S} \pi_S^*(
R^0\pi_{S*}(Q_S(-1))\otimes_{X_S} 
(M_{12})_S)\\
& {\buildrel {\widetilde  \sigma}_S\over 
\longto} \CO_{Y_S}(-1)\otimes_{Y_S} \pi_S^*(R^0\pi_{S*}Q_S  
\otimes_{X_S}(M_1\oplus M_2)_S)\\
& {\buildrel {\widetilde  \tau}_S\over 
\longto}  \pi_S^*R^0\pi_{S*}Q_S\to 0, \\
\end{aligned}
\ee
where 
\[
\begin{aligned}
{\widetilde  \sigma}_S& ={}^t(-z_{S,0}\otimes \pi_S^*\varPhi_{S,2}+
z_{S,2}\otimes \pi_S^*\varPhi_{S,0},\,
z_{S,0}\otimes \pi_S^* \varPhi_{S,1}-z_{S,1}\otimes \pi_S^*
\varPhi_{S,0})\\ 
{\widetilde \tau}_S& =(-z_{S,0}\otimes \pi_S^*
(\varPhi_{S,0}^{-1}\circ \varPsi_{S,1})+z_{S,1},\, 
-z_{S,0}\otimes \pi_S^*(\varPhi_{S,0}^{-1}\circ \varPhi_{S,2})+z_{S,2})
\end{aligned}
\]
and 
\be\label{eq:specseqCB}
\begin{aligned}
\varPhi_{S,0} : R^0\pi_{S*}(Q_S(-1)) & \to R^0\pi_{S*}Q_S\\
\varPhi_{S,1} : R^0\pi_{S*}(Q_S(-1))\otimes_{X_S} (M_1)_S & 
\to  R^0\pi_{S*}Q_S\\
\varPhi_{S,2} : R^0\pi_{S*}(Q_S(-1))\otimes_{X_S} (M_2)_S & 
\to  R^0\pi_{S*}Q_S\\
\end{aligned}
\ee
are induced by multiplication by $z_{S,0}, z_{S,1}, z_{S,2}$ respectively. 
Note that $\varPhi_{S,0}$ is an isomorphism since the support
of $Q_S$ is disjoint from $D_\infty$. Moreover we have obvious commutation 
relations 
\[
\varPsi_{S,i}\circ \varPhi_{S,0} = \varPsi_{S,0} \circ \varPhi_{S,i}
\]
for $i=1,2$. 

The complex \eqref{eq:specseqCA} is further isomorphic to 
\be\label{eq:specseqD}
\begin{aligned} 
0 & \to \CO_{Y_S}(-2) \otimes_{Y_S} \pi_S^*(
E_S\otimes_{X_S} 
(M_{12})_S)\\
& {\buildrel {\overline \sigma}_S\over 
\longto} \CO_{Y_S}(-1)\otimes_{Y_S} \pi_S^*(E_S  
\otimes_{X_S}(M_1\oplus M_2)_S)\\
& {\buildrel {\overline \tau}_S\over 
\longto}  \pi_S^*E_S\to 0, \\
\end{aligned}
\ee
where $E_S=R^0\pi_*Q_S$,  
\[
\begin{aligned}
{\overline \sigma}_S& ={}^t(-z_{S,0}\otimes \pi_S^*\Phi_{S,2}+z_{S,2},\,
z_{S,0}\otimes \pi_S^*\Phi_{S,1}-z_{S,1})\\ 
{\overline \tau}_S& =(-z_{S,0}\otimes\pi_S^*\Phi_{S,1}+z_{S,1}, \,
-z_{S,0}\otimes\pi_S^*\Phi_{S,2}+z_{S,2}),
\end{aligned}
\]
and 
\be\label{eq:thePhis}
\Phi_{S,i} = \varPsi_{S,0}^{-1} \circ \varPsi_{S,i}= 
\varPhi_{S,i}\circ \varPhi_{S,0}^{-1}
\ee 
for $i=1,2$. 
Note that $\Phi_{S,1},\Phi_{S,2}$ satisfy relation \eqref{eq:commrelA} 
since the 
morphisms $\varPhi_{S,1}, \varPhi_{S,2}$ in \eqref{eq:diffmap}
are induced by multiplication by $z_1,z_2$ respectively, therefore 
they commute.

Returning to the spectral sequence \eqref{eq:specseqB} 
note that all higher differentials vanish for degree reasons. 
Therefore the cohomology sheaves of 
the complex \eqref{eq:specseqD} must be trivial in degrees 
$-2,-1$ and isomorphic to 
$Q_S$ in degree $0$. Moreover, we claim that the natural evaluation map 
\[
\mathrm{ev}_S: \pi_S^*E_S \to Q_S 
\]
induces an isomorphism between the $0$-th cohomology sheaf of 
the complex \eqref{eq:specseqD} and $Q_S$. In order to prove this, 
note first that the evaluation map is surjective since its 
restriction to the fiber $Y_s$ is surjective for any $s\in S$. 
Furthermore, a simple computation shows that
\[
{\mathrm{Im}}({\overline \tau}_S)\subseteq \mathrm{Ker}(\mathrm{ev}_S).
\]
Therefore the evaluation map induces a surjective morphism 
of sheaves $\CH^0 \to Q_S$ where $\CH^0$ is the $0$-th cohomology sheaf 
of \eqref{eq:specseqD}. Since $\CH^0\simeq Q_S$, it follows that this 
morphism must be an isomorphism. This proves the claim. 

Let us denote by ${\bf B}(Q_S)$ the $[-1]$ shift 
of the complex \eqref{eq:specseqD} and by 
\[
(0,0,\mathrm{ev}_S):{\bf B}(Q_S) \to Q_S[-1]
\]
the quasi-isomorphism induced by the evaluation map. 

Next, note that by analogy with Lemma (\ref{relbeilinsonA}) 
we have a spectral sequence with first term 
\[
\CO_{Y_S}(j) \otimes_{Y_S} R^i\pi_{S*} \Omega^{-j}_{Y_S/X_S}(-j)
\]
converging to $\CO_{Y_S}$ for $i=j=0$ and $0$ otherwise. 
A simple calculation based on the relative Euler sequence and 
the base change theorem shows that this spectral sequence 
collapses to $\CO_{Y_S}$ for $i=j=0$ and $0$ otherwise. 

Now recall that for any coherent $\CO_{Y_S}$-module $F_S$, 
the derived direct image ${\bf R}p_{1S*} (\CK_{\Delta_S} 
\otimes_{Y_S\times_{X_S}Y_S} p_{2S}^*F_S)$ is computed by applying 
$\pi_{1S*}$ to an injective resolution of the complex 
$\CK_{\Delta_S} 
\otimes_{Y_S\times_{X_S}Y_S} p_{2S}^*F_S$. 
Since injective resolutions are functorial, it follows that the 
morphism $\rho_S: \CO_{Y_S} \to Q_S$ yields a morphism of relative 
Beilinson monads 
\[
\xymatrix{
\CO_{Y_S}[-1] \ar[r]^-{\pi_S^*\psi_S} & {\bf B}(Q_S)\\} 
\]
where $\psi_S: \CO_{X_S} \to E_S$ is the natural morphism 
of $\CO_{X_S}$-modules obtained by pushing-forward $\rho_S$.
Moreover, by construction we have 
diagram of complexes of $\CO_{Y_S}$-modules 
\[
\xymatrix{ 
\CO_{Y_S}[-1]\ar[d]^{\pi_S^*\psi_S} \ar[rr]^-{1_{\CO_{Y_S}}[-1]} & & \CO_{Y_S}[-1]
\ar[d]^-{\rho_S[-1]} \\
{\bf B}(Q_S) \ar[rr]^-{(0,0,\textrm{ev}_S)} & & Q_S[-1] \\}
\]
This yields a 
quasi-isomorphism of cones 
\[ 
{\sf Cone}(\xymatrix{\CO_{Y_S}[-1] \ar[r]^-{\pi_S^*\psi_S} & {\bf B}(Q_S)\\})
\to {\sf Cone}(\CO_{Y_S}{\buildrel -\rho_S\over \longto} Q_S)[-1]
\]
given by $(0,1_{\CO_{Y_S}}, \mathrm{ev}_S)$ in degrees $(-1,0,1)$ 
respectively.  
We will denote this quasi-isomorphism by $(0,1_{\CO_{Y_S}}, \textrm{ev}_S)$.
Note that 
\[
C(\CO_{Y_S}{\buildrel -\rho_S\over \longto} Q_S)[-1]= 
C(Q_S,\rho_S).
\] 
Furthermore, note that for each point the restriction 
$\psi_S|_{X_s}: \CO_{X_s} \to E_S|_{X_s}$ is nontrivial 
by construction since 
the morphism $\rho_S|_{Y_s}: CO_{Y_s}\to Q_S|_{Y_s}$ is nontrivial 
by definition. 
Therefore $\psi_S|_{X_s}$ is injective for any $s\in S$. 

In conclusion, the cone 
${\sf Cone}(\xymatrix{\CO_{Y_S}[-1] \ar[r]^-{\pi_S^*\psi_S} & {\bf B}(Q_S)\\})$
is the required resolution \eqref{eq:threetermA}. 

\hfill $\Box$

\begin{defi}\label{relmonad}
$(i)$ Given a flat family of admissible pairs $(Q_S,\rho_S)$ parameterized by 
a scheme $S$ of finite type over $\IC$, we define the relative monad 
${\bf B}(Q_S,\rho_S)$ associated to $(Q_S,\rho_S)$ to be the three term 
complex of locally free $\CO_{Y_S}$-modules 
\eqref{eq:threetermA}. 
In order to keep the notation short, we will denote the terms of this complex 
by 
\be\label{eq:threetermB} 
\begin{aligned}
0\to \CB_{-1}(Q_S,\rho_S) {\buildrel \sigma_S\over \longto} 
\CB_0(Q_S,\rho_S) {\buildrel \tau_S\over \longto} \CB_1(Q_S,\rho_S) 
\to 0.
\end{aligned}
\ee

$(ii)$
Conversely, given a flat family of stable ADHM sheaves $\CE_S$ on $X$ 
with $E_\infty=\CO_X$, 
Lemma (\ref{zerophi}) implies that $\phi_S=0$. Then the data 
$(E_S,\Phi_{S,1,2}, \psi_S)$  determines a 
complex ${\bf B}(\CE_S)$ of the form \eqref{eq:threetermA}. 
We will denote this complex by 
\be\label{eq:threetermBA}
\begin{aligned}
0\to \CB_{-1}(\CE_S) {\buildrel \sigma_S\over \longto} 
\CB_0(\CE_S) {\buildrel \tau_S\over \longto} \CB_1(\CE_S) 
\to 0.
\end{aligned}
\ee
Note that we have used the same symbols for the differentials 
in \eqref{eq:threetermB}, \eqref{eq:threetermBA} since the 
distinction will be clear from the context.
 
For further reference, let 
\be\label{eq:threetermBB}
\begin{aligned}
0\to \CB_{-1}'(\CE_S) {\buildrel \sigma'_S\over \longto} 
\CB_0'(\CE_S) {\buildrel \tau'_S\over \longto} \CB_1'(\CE_S) 
\to 0
\end{aligned}
\ee
be the three term subcomplex of ${\bf B}(\CE_S)$ obtained by removing 
the direct summand $\CO_{Y_S}$ in the degree $0$ term of \eqref{eq:threetermA}
and truncating $\sigma_S,\tau_S$ accordingly. 
We will denote this complex by ${\bf B}'(\CE_S)$.
\end{defi}

Standard properties of direct images and the base change theorem imply that 
relative monad complexes are compatible with base change i.e. the following 
holds. 
\begin{lemm}\label{relbasechange}
Let $(Q_S,\rho_S)$ be a flat family of admissible pairs parameterized by a 
scheme $S$ of finite type over $\IC$, and let $C(Q_S,\rho_S)$ be the 
complex \eqref{eq:twotermA}. Let $f:T\to S$ be a base 
change morphism, with $T$ of finite type over $\IC$, 
and let $(Q_T,\rho_T)=(f_Y)^*(Q_S,\rho_S)$.
Then we have the following commutative diagrams 
\be\label{eq:basechisom}
\begin{aligned}
& \quad \xymatrix{ 
f_Y^* {\bf B}(Q_S) \ar[rr]^{\simeq} \ar[d]_-{(0,0,{\textrm{ev}}_S)} 
& & {\bf B}(Q_T) \ar[d]^-{(0,0,{\textrm{ev}}_T)} \\
f_Y^* Q_S \ar[rr]^-{1_{Q_T}} & & Q_T \\}\\
& \xymatrix{ 
f_Y^* {\bf B}(Q_S,\rho_S) \ar[rr]^{\simeq} 
\ar[d]_-{(0,1_{\CO_{Y_S}},{\textrm{ev}}_S)} 
& & {\bf B}(Q_T,\rho_T) \ar[d]^-{(0,1_{\CO_{Y_T}},{\textrm{ev}}_T)} \\
f_Y^* C(Q_S,\rho_S) \ar[rr]^-{(1_{\CO_{Y_T}}, 1_{Q_T})} & & Q_T, \\}\\
\end{aligned}
\ee
where the top rows are canonical isomorphisms. 
\end{lemm}

\subsection{Admissible Pairs and ADHM Sheaves} 

Next we prove that there is a one-to-one correspondence 
between admissible pairs on $Y$ and stable ADHM sheaves on 
$X$ with $E_\infty=\CO_X$. We start with a preliminary lemma 
summarizing the relevant properties of the complexes 
${\bf B}(\CE_S)$, ${\bf B}'(\CE_S)$ defined in 
(\ref{relmonad}.$ii$).  

\begin{lemm}\label{relmonprop} 
Let $X$ be a smooth projective curve over an infinite field $k$ of 
characteristic 0, and let 
$\CE=(E,\Phi_{1,2}, \psi)$ be a stable ADHM sheaf on $X$ with 
$E_\infty = \CO_X$. Then the following hold
\begin{itemize} 
\item[$(i)$] $\sigma:\CB_{-1}(\CE)\to \CB_0(\CE)$ is an injective 
morphism of $\CO_Y$-modules 
and $\mathrm{Coker}(\sigma)$ is a coherent torsion free 
$\CO_Y$-module. 
\item[$(ii)$] The middle cohomology sheaf $\CH^0({\bf B}'(\CE))$ 
is trivial. 
\item[$(iii)$] The morphism of $\CO_Y$-modules 
$\tau:\CB_0(\CE)\to \CB_1(\CE)$ is surjective 
on the complement of a codimension three closed subset of $Y$. 
\item[$(iv)$] The restriction $\sigma|_{Y_x}: 
\CB_{-1}(\CE)|_{Y_x} \to {\CB}_0(\CE)|_{Y_x}$ 
to an arbitrary fiber $Y_x$, $x\in X$ is an injective 
morphism of $\CO_{Y_x}$-modules and $\mathrm{Coker}(\sigma|_{Y_x})$ 
is a coherent torsion free $\CO_{Y_x}$-module. 
\item[$(v)$] The middle cohomology sheaf 
$\CH^0({\bf B}'(\CE)|_{Y_x})$ is trivial for any point $x\in X$. 
\end{itemize}
\end{lemm}

{\it Proof.} 
Suppose $\sigma:\CB_{-1}(\CE)\to \CB_0(\CE)$ is not injective.  
Then $\mathrm{Ker}(\sigma)$ 
must be a nontrivial subsheaf of $\CB_{-1}(\CE)=\CO_Y(-2)\otimes_Y
\pi^*(E\otimes_X 
M_{12})$, which is a locally free $\CO_Y$-module.
Therefore $\mathrm{Ker}(\sigma)$  has to be torsion-free, 
hence locally free outside a codimension two closed subset 
of $Y$. Let $Y_\sigma$ be the maximal open subset of $Y$ such that 
${\mathrm{Ker}}(\sigma)|_{Y_\sigma}$ is a locally free
$\CO_{Y_\sigma}$-module. Since $Y_\sigma \subset Y$ is an 
open embedding, we have $\mathrm{Ker}(\sigma|_{Y_\sigma}) = 
\mathrm{Ker}(\sigma)|_{Y_\sigma}$, therefore 
$\mathrm{Ker}(\sigma|_{Y_\sigma})$ is a locally free 
$\CO_{Y_\sigma}$-module as well, and we have a 
short exact sequence of $\CO_{Y_\sigma}$-modules 
\[
0\to \mathrm{Ker}(\sigma|_{Y_\sigma}) \to \CO_{Y_\sigma}(-2)\otimes_{Y_\sigma}
\pi^*(E\otimes_X 
M_{12})|_{Y_\sigma} \to \mathrm{Coker} (\sigma|_{Y_\sigma})\to 0.
\]
in which the first two terms are locally free. Then
the associated ${\mathcal Tor}^{\CO_{Y_\sigma}}$ 
long exact sequence shows that 
\[
{\mathcal Tor}_1^{\CO_{Y_\sigma}}(\mathrm{Coker} (\sigma|_{Y_\sigma}), 
\CO_{y}) =0 
\] 
where $\CO_{y}$ is the structure sheaf 
of any point $y\in Y$. 
This implies that for any closed point $x\in X$, 
the induced map of $k(y)$-vector spaces 
\[
\sigma(y) : E(x)\otimes_{k(x)}k(y) \to (E(x)\oplus E(x) \oplus k(x))
\otimes_{k(x)}k(y)
\]
has a nontrivial kernel for any point $y\in Y_\sigma \cap Y_x$. 
However, as observed in the proof \cite[Lemma 2.7]{hilblect}, 
for any $x\in X$, $\sigma(y)$ 
may have a nontrivial kernel at most at finitely many 
closed points $y\in Y_x$.  
Hence $\sigma$ must be injective in order to avoid 
a contradiction. 

In order to prove that  
$\mathrm{Coker}(\sigma)$ is torsion free, note 
that $\mathrm{Coker}(\sigma)_y$ is a free $\CO_{Y,y}$-module 
at all points $y\in Y$ where $\sigma(y)$ is injective, according to 
Lemma (\ref{EGAlemmaA}). Therefore 
$\mathrm{Coker}(\sigma)$ is locally free on the complement of a closed 
subset of codimension two. Let $\omega_Y$ be the dualizing sheaf 
of $Y$; since $Y$ is smooth and projective $\omega_Y$ is locally free. 
The long exact ${\mathcal Ext}_Y^{\bullet}(\, {\underline {\ \ }}\, ,\omega_Y)$
sequence associated to the short exact sequence 
\[
\begin{aligned} 
0 \to \CB_{-1}(\CE) {\buildrel \sigma\over \longto} \CB_{0}(\CE)\to 
\mathrm{Coker}(\sigma)\to 0
\end{aligned}
\]
yields 
\[
{\mathcal Ext}^{q}_Y(\mathrm{Coker}(\sigma), \omega_Y) =0 
\]
for all $q\geq 2$ since the first two terms are locally free. 
Moreover, the local ext sheaf ${\mathcal Ext}^{1}_Y(\mathrm{Coker}(\sigma), \omega_Y)$
must be supported in codimension 2 since $\mathrm{Coker}(\sigma)$ 
is locally free outside a codimension 2 closed subset. 
Then \cite[Prop. 1.1.10]{huylehn} implies that $\mathrm{Coker}(\sigma)$ 
must be torsion free.
This proves (\ref{relmonprop}.$i$).

In order to prove that the middle cohomology sheaf $\CH^0({\bf B}'(\CE))$ 
is trivial, it suffices to prove that its stalk $\CH^0({\bf B}'(\CE))_y$ 
is trivial for any closed point $y\in Y$. 
This follows from a simple local computation. 
Given the explicit expressions of the differentials 
\[
\sigma' = {}^t(-z_0\otimes \pi^*\Phi_2+z_2, z_0\otimes \pi^*\Phi_1-z_1),
\qquad 
\tau' = (-z_0\otimes \pi^*\Phi_1 + z_1, -z_0\otimes \pi^*\Phi_2+z_2) 
\]
where $\Phi_1,\Phi_2$ satisfy \eqref{eq:commrelA}, it follows that 
the complex ${\bf B}'(\CE)_y$ is exact in degree $0$ if at least one of 
the morphisms of $\CO_{Y,y}$-modules 
\[
(z_0\otimes \pi^*\Phi_1 - z_1)_y,\qquad  (z_0\otimes \pi^*\Phi_2-z_2)_y
\]
is an isomorphism. Therefore $\CH^0({\bf B}'(\CE))$ 
must be supported on 
the subset of $Y$ where both endomorphisms 
fail to be isomorphisms, which is a codimension two 
condition on $Y$. 
However note that $\CH^0({\bf B}'(\CE))$ is a subsheaf of 
$\mathrm{Coker}(\sigma')$.
Since the last component of $\sigma$ is trivial, we have 
\[
\mathrm{Coker}(\sigma)\simeq \mathrm{Coker}(\sigma') \oplus \CO_Y,
\]
therefore $\mathrm{Coker}(\sigma')$ must be torsion free since 
$\mathrm{Coker}(\sigma)$ is torsion-free according to 
(\ref{relmonprop}.$i$). In conclusion the cohomology sheaf 
$\CH^0({\bf B}'(\CE))$ must be trivial. 

For the next claim, note that 
$\tau_y$ is surjective on stalks  
at a point $y\in Y$ if and only if the induced linear map 
of $k(y)$-vector spaces 
\[
\tau(y) : (E(x)\oplus E(x)\oplus k(x))\otimes_{k(x)}k(y) 
\to E(x)\otimes_{k(x)}k(y) 
\]
is surjective, where $x=\pi(y)$. 
If $\tau(y)$ fails to be surjective at some $y\in Y$, 
it follows that the dual map $\tau(y)^\vee$ will fail to be injective.  
As observed in the proof of
\cite[Lemma 2.7]{hilblect}, this 
 implies the existence of a nontrivial proper linear subspace 
$W\subset E(x)\otimes_{k(x)}k(y)$ such that 
\[
(\Phi_{1,2}(x)\otimes 1_{k(y)})(W) \subseteq W, 
\qquad \mathrm{Im}(\psi(x)\otimes 1_{k(y)})\subseteq W
\]
Recall that in Lemma (\ref{unstablemma}) we have constructed 
a canonical destabilizing subsheaf $E_0\subset E$. Note that  
$E_0(x)$ is a linear subspace of $E(x)$ for a generic closed 
point $x\in X$. Moreover, if a subspace $W\subset E(x)\otimes_{k(x)}k(y)$ 
as above exists for such a generic closed point $x\in X$, 
by construction we have $E_0(x)\otimes_{k(x)}k(y)\subset W(x)$. 
Therefore $\mathrm{rk}(E_0)<\mathrm{rk}(E)$, and $E_0$ 
must be a proper subsheaf of $E$ in this case. In conclusion, 
$\tau(y)$ must be surjective for any point $y\in Y_x$, for generic 
$x\in X$. This implies that there must exist a closed subset 
$S_\CE\subset X$ determined by $\CE$ such that 
$\tau(y)$, hence also $\tau_y$, is surjective for any 
$y\in \pi^{-1}(X\setminus S_\CE)$.
 
Now suppose $x\in S_\CE$. The condition that 
$\tau|_{Y_x}$ fail to be surjective is a codimension two 
condition along $Y_x$. Therefore $\tau_y$ may fail to be surjective 
at most along a codimension three locus in $Y$. Moreover, 
it is easy to check that $\tau$ is surjective along $D_\infty$, 
which is cut by $z_0=0$, hence the codimension three subset 
in question must be contained in $Y\setminus D_\infty$. 

The proofs of (\ref{relmonprop}.$iv$), (\ref{relmonprop}.$v$) are
analogous to the proofs of (\ref{relmonprop}.$i$), (\ref{relmonprop}.$ii$)
respectively, therefore will be omitted. 

\hfill $\Box$

\begin{prop}\label{correspA}
Let $X$ be a smooth projective curve over an infinite field $k$ 
of characteristic 0. 
Then there is a one-to-one correspondence between admissible pairs
$(Q,\rho)$ on $X$ with $\ch_2(Q)=d\beta$ and $\chi(Q) = n$ 
and stable ADHM sheaves $\CE$ on $X$ with $E_\infty=\CO_X$ and 
\be\label{eq:topinv}
r=d,\qquad  e= n +d(g-1).
\ee
Moreover, two admissible pairs $(Q,\rho)$, $(Q',\rho')$ are isomorphic 
if and only if the corresponding ADHM sheaves $\CE$, $\CE'$ are isomorphic.
\end{prop}

{\it Proof.} 
Given an admissible  pair $(Q,\rho)$, Lemma 
(\ref{relbeilinsonB}) provides a relative monad 
${\bf B}(Q,\rho)$. The data $(E,\Phi_{1,2},\psi)$ 
defines an ADHM sheaf 
$\CE$ with $E_\infty=\CO_X$ and $\phi=0$.

We claim that the resulting ADHM sheaf $\CE$ must be stable. 
Suppose there exists a nontrivial saturated proper subsheaf 
$E'\subset E$ such that 
\be\label{eq:unstableA}
\Phi_{i}(E'\otimes_X M_i) \subseteq E', \qquad 
\mathrm{Im}(\psi) \subseteq E'
\ee
for $i=1,2$.
Note that $E'$ must be torsion free, hence locally free, on $X$ 
and $\mathrm{rk}(E')<\mathrm{rk}(E)$. Let $\Phi_i'= 
\Phi_i|_{E'\otimes M_i}$, $i=1,2$ and let $\psi' :\CO_X \to E'$ denote 
the factorization of $\psi$ through $E'\subset E$; note that 
$\mathrm{Im}(\psi') = \mathrm{Im}(\psi)$ as subsheaves of $E'$.  
Note that the sheaf inclusions \eqref{eq:unstableA} 
yield the following exact sequences of $\CO_{X}$-modules 
\be\label{eq:unstableB}
\begin{aligned}
& {\mathcal Tor}_1^{\CO_X}(\mathrm{Coker}(\Phi'_i), \CO_x) 
\to \mathrm{Im}(\Phi'_i)\otimes_X \CO_x \to E'\otimes_X \CO_x \\
& {\mathcal Tor}_1^{\CO_X}(\mathrm{Coker}(\psi'), \CO_x) 
\to \mathrm{Im}(\psi') \otimes_X \CO_x \to E' \otimes_X \CO_x ,\\
\end{aligned}
\ee
where $\CO_x$ is the structure sheaf of an arbitrary closed point $x\in X$. 

According to \cite[Def. 1.1.4]{huylehn}, any coherent $\CO_X$-module 
$F$ admits a torsion sub-module $T_0(F)$ such that the quotient 
$F/T_0(F)$ is torsion free, hence locally free, on $X$. Then the long exact 
${\mathcal Tor}^{\CO_X}$ sequence associated to the short exact sequence 
\[
0 \to T_0(F) \to F \to F/T_0(F) \to 0 
\]
implies that 
\[
{\mathcal Tor}_1^{\CO_X}(F,\CO_x) \simeq 
{\mathcal Tor}_1^{\CO_X}(T_0(F), \CO_x).
\]
Therefore for any coherent $\CO_X$ module $F$, 
${\mathcal Tor}_1^{\CO_X}(F,\CO_x)$ 
is trivial unless $x\in X$ belongs to the support of $T_0(F)$, which 
consists of finitely many points on $X$. 

Applying this argument to $F=\mathrm{Coker}(\Phi'_i), \ 
\mathrm{Coker}(\psi')$ in \eqref{eq:unstableB} it follows that 
\be\label{eq:unstableC}
\mathrm{Im}(\Phi'_i(x)) \subseteq E'(x) \qquad 
\mathrm{Im}(\psi(x)) \subseteq E'(x) 
\ee
for all except finitely many closed points $x\in X$. 
Now for each point $y\in Y$, 
the morphism of $\CO_{Y,y}$-modules $\tau_y$ is surjective if and 
only if the linear map of $k(y)$-vector spaces 
\[
\tau(y): (E(x)\oplus E(x)\oplus k(x))\otimes_{k(x)}k(y) 
\to E(x)\otimes_{k(x)}k(y)
\] 
is surjective. 
This is equivalent to the dual 
linear map $\tau(y)^\vee$ being injective. 

However, as observed in the proof of lemma  \cite[Lemma 2.7]{hilblect}, 
the existence of a proper nontrivial subspace $E'(x)\subset E(x)$
satisfying conditions \eqref{eq:unstableC} implies that $\tau(y)^\vee$ 
will fail to be injective at finitely many closed points in $Y_x$. 
In conclusion, the existence of a destabilizing subsheaf $E'\subset E$ 
implies that $\tau$ will fail to be surjective at finitely many 
closed points $y\in Y_x$ for generic $x\in X$. 

Now, according to Lemma (\ref{relbeilinsonB}) 
the complex \eqref{eq:threetermA} 
is quasi-isomorphic to the two term complex $C(Q,\rho)$.
Therefore $\mathrm{Coker}(\tau) \simeq \mathrm{Coker}(\CO_Y
{\buildrel \rho\over \longto} Q)$. This implies that  
the morphism $\tau$ in \eqref{eq:threetermA}
must be surjective everywhere on $Y$ except the support 
of $\mathrm{Coker}(\CO_Y {\buildrel \rho 
\over \longto} Q)$, which is dimension zero by definition.
This contradicts the conclusion reached in the previous paragraph. 
Therefore a destabilizing proper subsheaf $E'\subset E$ as above 
cannot exist. 

The relations \eqref{eq:topinv} follow by a straightforward 
Grothendieck-Riemann-Roch computation using Lemma (\ref{usefulstuff}). 

Conversely, suppose we have a stable ADHM sheaf $(E,\Phi_{1,2},\phi,\psi)$
of type $(r,e)=(d,n+d(g-1))$ 
with $E_\infty=\CO_X$. According to (\ref{zerophi}) we must have $\phi=0$. 
Then we construct the complex ${\bf B}(\CE_S)$ of the form 
according to Definition (\ref{relmonad}.$ii$). 
Note that 
we have a commutative diagram of $\CO_Y$-modules 
with exact rows of the form 
\[
\xymatrix{ 
0\ar[r] & \textrm{Coker}(\sigma') \ar[r]\ar[d]^-{{\overline\tau}'} 
& \textrm{Coker}(\sigma)\ar[r]\ar[d]^-{{\overline\tau}} & 
\CO_Y \ar[d]\ar[r] & 0\\
0\ar[r] & \CB_1(\CE) \ar[r]^-{1}& \CB_1(\CE) \ar[r]& 0\ar[r] & 0\\}
\]
where ${\overline \tau}', {\overline \tau}$ are morphisms induced 
by $\tau', \tau$. 
Applying the snake lemma we obtain an exact sequence of 
$\CO_Y$-modules of the form 
\be\label{eq:truncatedD}
\begin{aligned} 
0\to \mathrm{Ker}({\overline \tau}') \to \mathrm{Ker}({\overline \tau}) 
\to \CO_Y \to \mathrm{Coker}({\overline \tau}') \to 
\mathrm{Coker}({\overline \tau}) \to 0.
\end{aligned}
\ee
However Lemma (\ref{relmonprop}.$ii$) implies 
that $\mathrm{Ker}({\overline \tau}')=0$ 
since the complex ${\bf B}'(\CE)$ is exact in degree 0. 
Therefore it follows that 
\[
\mathrm{Ker}({\overline \tau})\simeq \mathrm{Ker}(\tau)/\mathrm{Im}(\sigma)
\] 
is a subsheaf of $\CO_Y$. Moreover, it is easy to check that the 
morphism  
\[
\mathrm{Ker}(\tau)/\mathrm{Im}(\sigma)|_{D_\infty} \to {\CO_Y}|_{D_\infty}
\]
is an isomorphism. 
This implies that $\mathrm{Ker}({\overline \tau})$
is a torsion-free rank one sheaf on $Y$ with trivial determinant, 
therefore it 
must be the ideal 
sheaf of a closed subscheme $Z\subset Y$.
Moreover the support of $Z$ is disjoint from 
$D_\infty$. 

Set  $Q=\mathrm{Coker}({\overline \tau}')$ 
and let $\rho:\CO_Y\to Q$ be the morphism obtained from \eqref{eq:truncatedD}.
 Since $Q$ is the only nontrivial 
cohomology sheaf of the complex ${\bf B}'(\CE)$, 
a straightforward computation using the first relation in 
Lemma (\ref{usefulstuff}) yields 
\be\label{eq:cherncompA}
\ch_0(Q) =0,\qquad \ch_1(Q)=0,\qquad \ch_2(Q) = d\beta.
\ee
Hence $Q$ satisfies condition $(i.c)$ of Definition (\ref{admissible}). 
Moreover, we claim that the morphism $\rho$ is not identically zero. 
If it were trivial, we would obtain $\mathrm{Coker}({\overline \tau}') 
\simeq \mathrm{Coker}({\overline \tau})$ from \eqref{eq:truncatedD}. 
However $\mathrm{Coker}({\overline \tau})$ is zero dimensional 
according to Lemma (\ref{relmonprop}.$iii$), 
which would contradict \eqref{eq:cherncompA}. 

We also have to prove that $Q$ is flat over $X$. 
Note that according to Lemma (\ref{relmonprop}.$v$) 
the restriction 
\[
{\overline \tau}'|_{Y_x}: 
\mathrm{Coker}(\sigma')|_{Y_x}\to \CB_1(\CE)|_{Y_x} 
\]
is injective for all $x\in X$. Therefore $Q= \mathrm{Coker}({\overline
  \tau'})$
is flat over $X$ according to Lemma (\ref{EGAlemmaA}). 

The relations \eqref{eq:topinv} follow again by a standard 
Grothendieck-Riemann-Roch computation using Lemma (\ref{usefulstuff}). 

Finally, compatibility with isomorphisms follows from a  
routine verification. 

\hfill $\Box$

\subsection{Isomorphism of Moduli Spaces}\label{isomsection}

In this section we prove Theorem (\ref{correspB}). 
We first prove that the relative Beilinson spectral sequence 
yields a {\bf T}-equivariant 
isomorphism of algebraic spaces 
\[
{\mathfrak f}:\mfm_{Adm}(Y,d,n) 
\simeq \mfm_{ADHM}(\CX,d,n+d(g-1)).
\]  
According to Lemma (\ref{relbeilinsonB}), a 
flat family of admissible pairs $(Q_S,\rho_S)$ parameterized by a
scheme $S$ of finite type over $\IC$ determines an ADHM sheaf 
$\CE_S=(E_S,\Phi_{S,1,2}, 0 , \psi_S)$ on $X_S$. 
$E_S$ is a locally free $\CO_{X_S}$-module, 
hence it is flat over $S$. Lemma (\ref{relbasechange})
and Proposition (\ref{correspA}) 
imply that $\CE_S|_{X_s}$ is a stable ADHM sheaf 
for any point $s\in S$. Moreover,
Lemma (\ref{relbasechange}) and Proposition (\ref{correspA}) 
further imply that this correspondence is compatible 
with isomorphisms and base change, therefore we obtain a 
functor ${\mathfrak f}:\mfm_{Adm}(Y,d,n) \to \mfm_{ADHM}(\CX,d,n+d(g-1))$. 

Conversely, suppose  $\CE_S=(E_S,\Phi_{S,1,2}, \phi_S , \psi_S)$ is a flat 
family of stable ADHM sheaves on $X$ parameterized by $S$. 
According to Lemma (\ref{zerophi}), the restriction 
$\phi_S|_{X_s}$ must be trivial for each closed point 
$s\in S$. Therefore  $\phi_S=0$.

Consider the three-term complex ${\bf B}(\CE_S)$ defined in (\ref{relmonad}.
$ii$). The restriction of ${\bf B}(\CE_S)$ to each fiber $Y_s$, $s\in S$ 
is isomorphic to ${\bf B}(\CE_s)$, where $\CE_s=\CE_S|_{Y_s}$. 
Therefore Lemma (\ref{relmonprop}.$i$) and Lemma (\ref{EGAlemmaA}) 
imply that $\sigma_S : \CB_{-1}(\CE_S)\to \CB_0(\CE_S)$ is 
an injective morphism of $\CO_{Y_S}$-modules, and 
$\mathrm{Coker}(\sigma_S)$ is flat over $S$. 
The same holds for $\sigma_{S}': \CB'_{-1}(\CE_S)\to \CB'_0(\CE_S)$
by a similar argument. 
Let 
\[
{\overline \tau'}_S: \mathrm{Coker}(\sigma'_S) \to \CB_1(\CE_S)
\]
be the morphism of $\CO_{Y_S}$-modules induced by $\tau'_S$. 
Then applying again (\ref{relmonprop}.$v$) and Lemma (\ref{EGAlemmaA})
it follows that $\mathrm{Ker}({\overline \tau'}_S)=0$
and $Q_S=\mathrm{Coker}({\overline \tau'}_S)$ is 
flat over $S$. Furthermore, since 
\be\label{eq:restr-Q}
Q_S|_{Y_s} \simeq \mathrm{Coker}({\overline \tau'}_S|_{Y_s}),
\ee
Proposition (\ref{correspA}) and Lemma (\ref{EGAlemmaB}) imply that 
$Q_S$ is flat over $X_S$ and 
\[
\mathrm{ch}_0(Q_S|_{Y_s})=0, \qquad 
\mathrm{ch}_1(Q_S|_{Y_s})=0,\qquad 
\mathrm{ch}_2(Q_S|_{Y_s})=d\beta, \qquad \chi(Q_S|_{Y_s})=n
\]
for any $s\in S$. 

Now, applying the snake lemma as in  the proof of the inverse 
implication of Proposition (\ref{correspA}) 
(above equation \eqref{eq:truncatedD}) 
we obtain an exact sequence of 
$\CO_{Y_S}$-modules of the form
\be\label{eq:famtruncatedB}
\begin{aligned} 
0 \to \mathrm{Ker}({\overline \tau}_S) 
\to \CO_{Y_S} {\buildrel \rho_S\over \longto} Q_S \to 
\mathrm{Coker}({\overline \tau}_S) \to 0.
\end{aligned}
\ee
Using \eqref{eq:famtruncatedB} and the fact that tensor product 
is right exact, we have isomorphisms of $\CO_{Y_s}$-modules 
\[
\mathrm{Coker}(\rho_S|_{Y_s}) \simeq \mathrm{Coker}(\rho_S)|_{Y_s} 
\simeq \mathrm{Coker}({\overline \tau}_S)|_{Y_s}  
\simeq \mathrm{Coker}({\overline \tau}_S|_{Y_s})
\simeq \mathrm{Coker}({\tau}_S|_{Y_s})
\]
for any $s\in S$. Then Proposition (\ref{correspA}) 
implies that $ \mathrm{Coker}(\rho_S|_{Y_s})$ is a sheaf of 
pure dimension 0 for any $s\in Y$. Therefore we conclude that 
$(Q_S,\rho_S)$ is indeed a flat family of admissible pairs of 
type $(d,n)$ parameterized by $S$.

It is again routine to check that this construction is compatible with 
isomorphisms and base change, therefore we obtain an inverse functor 
${\mathfrak g} : \mfm_{ADHM}(\CX,d,n+d(g-1)) \to \mfm_{Adm}(Y,d,n)$. 
By construction, the functors ${\mathfrak f}, \ {\mathfrak g}$ 
are inverse to each other. Therefore we have proven that the two 
stacks are indeed isomorphic. It is also easy to check that the 
functors ${\mathfrak f}, \ {\mathfrak g}$ are {\bf T}-equivariant.

\subsection{Comparison of Tangent-Obstruction Theories}\label{compattangobs}
Next we prove that the perfect 
tangent-obstruction theories of the isomorphic moduli spaces 
$\mfm_{Adm}(Y,d,n)$, 
$\mfm_{ADHM}(\CX,d,n+d(g-1))$, 
are compatible with the respect to the isomorphism ${\mathfrak f}$
in the sense of definition
\cite[Def 4.1]{degGW}. 
More precisely, let 
$(\iota_\alpha :S_\alpha\to \mfm_{Adm}(Y,d,n))$, $\alpha \in \Lambda$, be a 
finite affine open cover of $\mfm_{Adm}(Y,d,n)$ , each $S_\alpha$ being 
preserved by the torus action. Note that $({\mathfrak f}\circ \iota_\alpha 
: S_\alpha \to \mfm_{ADHM}(\CX,d,n+d(g-1))$, $\alpha \in \Lambda$, 
is a similar 
cover of $\mfm_{ADHM}(\CX,d,n+d(g-1))$. 
Let $\CQ_\alpha=(Q_\alpha, \rho_\alpha)$, respectively 
$\CE_\alpha$ denote corresponding
flat families of admissible pairs, respectively stable ADHM sheaves 
for all $\alpha \in \Lambda$. 
Then we will prove the following 
\begin{itemize}
\item[$(I)$] For each $\alpha\in \Lambda$ and any {\bf T}-equivariant 
coherent $\CO_{S_\alpha}$-module $I$, there is a natural 
{\bf T}-equivariant isomorphism
of tangent-obstruction 
theories associated to the data $(\CQ_\alpha,
I)$, $(\CE_\alpha, I)$. These isomorphisms are furthermore compatible 
with base change morphisms of the form $f_\alpha : T_\alpha \to S_\alpha$.
\item[$(II)$] There exist {\bf T}-equivariant 
perfect tangent-obstruction data $(\IF_\alpha^\bullet)$, $(\IE_\alpha^\bullet)$
and {\bf T}-equivariant quasi-isomorphisms 
\[
{\sf q}_\alpha : \IF_\alpha^{\bullet} \to \IE_\alpha^{\bullet}
\]
of complexes of $\CO_{S_\alpha}$-modules such that 
\begin{itemize} 
\item[$(a)$] The induced isomorphisms
in cohomology
\[
\CH^k({\sf q}_\alpha)(I) : \CH^k(\IF_\alpha^\bullet \otimes_{S_\alpha} I) 
{\buildrel \sim \over \longto} \CH^k(\IE_\alpha^\bullet \otimes_{S_\alpha} I) 
\]
$k=1,2$ agree with the natural isomorphisms constructed at point $(I)$ above
for all $\alpha\in \Lambda$. 
\item[$(b)$] For $k=2$, the isomorphisms 
\[
\CH^2({\sf q}_\alpha): \CH^2(\IF_\alpha^\bullet) 
{\buildrel \sim \over \longto} \CH^k(\IE_\alpha^\bullet) 
\]
determine a global isomorphism of obstruction sheaves. 
\item[$(c)$] Given any deformation situation $S_\alpha\subset S'_\alpha
\subset S_\alpha''$ as in \eqref{eq:defsit}, the obstructions 
${\mathfrak {ob}}_\alpha(\CQ_\alpha, S_\alpha', S_\alpha'')$, 
${\mathfrak {ob}}_\alpha(\CE_\alpha, S_\alpha', S_\alpha'')$
agree i.e. 
\[
\CH^2({\sf q}_\alpha)(I_{S_\alpha'\subset S''_\alpha})
({\mathfrak {ob}}_\alpha(\CQ_\alpha, S_\alpha',
S_\alpha''))
={\mathfrak {ob}}_\alpha(\CE_\alpha, S_\alpha', S_\alpha'')
\]
for all $\alpha \in \Lambda$.
\end{itemize}
\end{itemize}
Since all the computations below extend trivially to the equivariant setting, for brevity
 we will not explicitely state this in the following lemmas. 
Moreover compatibility with base change will be always analogous to the 
condition formulated in 
Lemma (\ref{tangspace}.$ii$) and will not be made explicit in the following. 

\begin{lemm}\label{tracelemma}
Let $(Q_S,\rho_S)$ be a flat family of admissible pairs 
parameterized by an affine scheme $S$ of finite type over $\IC$, 
and let $I$ be a coherent $\CO_S$-module. 
Then we have canonical isomorphisms compatible with base change 
\be\label{eq:traceisomA}
{\mathcal Ext}^k_{q_S}(C(Q_S,\rho_S), C(Q_S,\rho_S)\otimes_{Y_S}q_S^*I)_0
\simeq {\mathcal Ext}^k_{q_S}(C(Q_S,\rho_S), 
C(Q_S,\rho_S)(-1)\otimes_{Y_S}q_S^*I)
\ee
for any $k\in \IZ$, where $q_S:Y_S\to S$ is the natural projection. 
\end{lemm}

{\it Proof.} 
According to Lemma (\ref{relbeilinsonB}), the relative Beilinson
monad ${\bf B}(Q_S,\rho_S)$ is a locally free resolution of 
$C(Q_S,\rho_S)$ and there is a canonical quasi-isomorphism 
\[
(0,1_{\CO_{Y_S}}, \mathrm{ev}_S) : {\bf B}(Q_S,\rho_S)\to 
C(Q_S,\rho_S).
\] 
Therefore $C(Q_S\rho_S)$ can be canonically replaced 
by ${\bf B}(Q_S,\rho_S)$  
in Lemma (\ref{tracelemma}). 
Since the complex ${\bf B}(Q_S,\rho_S)$ is locally free, 
the trace map yields a morphism of complexes of $\CO_{Y_S}$-modules 
\be\label{eq:tracemapA}
\xymatrix{
0 \ar[d] & 0\ar[d] \\
{\mathcal Hom}_{Y_S}({\bf B}(Q_S,\rho_S), 
{\bf B}(Q_S,\rho_S)(-1)\otimes_{Y_S}q_S^*I) \ar[r]^-{\textrm{tr}(-1)}\ar[d]
& q_S^*I(-1)\ar[d]\\
{\mathcal Hom}_{Y_S}({\bf B}(Q_S,\rho_S), 
{\bf B}(Q_S,\rho_S)\otimes_{Y_S}q_S^*I) \ar[r]^-{\textrm{tr}}\ar[d]
& q_S^*I\ar[d]\\
{\mathcal Hom}_{Y_S}({\bf B}(Q_S,\rho_S), 
{\bf B}(Q_S,\rho_S)\otimes_{Y_S}q_S^*I)\otimes_{Y_S}\CO_{(D_\infty)_S} 
\ar[r]^-{\textrm{tr}_\infty}\ar[d]
& q_S^*I\otimes_{Y_S}\CO_{(D_\infty)_S}\ar[d]\\
0 & 0. \\}
\ee
A straightforward calculation shows that there is a canonical 
quasi-isomorphism 
\[
{\bf B}(Q_S,\rho_S) \otimes_{Y_S}\CO_{(D_\infty)_S} \to \CO_{(D_\infty)_S}
\]
and the bottom horizontal row
of \eqref{eq:tracemapA} is a quasi-isomorphism of complexes as well.
In order to keep the notation short we will denote by $\CH(-1), \CH,
\CH_\infty$ the terms in the left column of \eqref{eq:tracemapA}, 
starting from the top. Then note 
 that
the diagram \eqref{eq:tracemapA} yields a commutative diagram of 
$\CO_{S}$-modules of the form 

\be\label{eq:tracemapC}
\xymatrix{
{\bf R}^{k-1}_{q_S*}(\CH_\infty) \ar[rr]^-{{\bf R}^{k-1}q_{S*}
(\textrm{tr}_\infty)} \ar[d] & & {\bf R}^{k-1}q_{S*}
(q_S^*I\otimes_{Y_S}\CO_{(D_\infty)_S})\ar[d]\\
{\bf R}^kq_{S*}(\CH(-1)) \ar[rr]^-{{\bf R}^{k}q_{S*}(\textrm{tr}(-1))} \ar[d]
&  &  {\bf R}^{k}q_{S*}
(q_S^*I(-1))\ar[d]\\ 
{\bf R}^kq_{S*}(\CH) \ar[rr]^-{{\bf R}^{k}q_{S*}(\textrm{tr})} \ar[d]
&  &  {\bf R}^{k}q_{S*}
(q_S^*I)\ar[d]\\
{\bf R}^{k}_{q_S*}(\CH_\infty) \ar[rr]^-{{\bf R}^{k}q_{S*}
(\textrm{tr}_\infty)} & & {\bf R}^{k}q_{S*}
(q_S^*I\otimes_{Y_S}\CO_{(D_\infty)_S})\\
}
\ee
in which the columns are exact and the horizontal morphisms of the form 
${\bf R}^{k}q_{S*}
(\mathrm{tr}_\infty)$
are isomorphisms. It is also straightforward to check that 
${\bf R}^{k}q_{S*}(q_S^*I(-1))=0$
for all $k\in\IZ$. 
Then applying the snake lemma to the commutative diagram 
\eqref{eq:tracemapC}, we obtain isomorphisms 
\[
{\mathrm {Ker}}({\bf R}^{k}q_{S*}(\textrm{tr})) \simeq 
{\bf R}^kq_{S*}{\mathcal Hom}_{Y_S}({\bf B}(Q_S,\rho_S), 
{\bf B}(Q_S,\rho_S)(-1)\otimes_{Y_S}q_S^*I)
\]
for all $k\in \IZ$. Since ${\bf B}(Q_S,\rho_S)$ is a locally free 
resolution of $C(Q_S,\rho_S)$, this proves Lemma (\ref{tracelemma}.$i$).

Compatibility with base change follows from Lemma (\ref{relbasechange}) 
standard properties of direct images.  

\hfill $\Box$

\begin{lemm}\label{tangobscompA}
$(i)$ Let $(Q_S,\rho_S)$, $(Q_S',\rho_S')$ be flat families of admissible 
pairs parameterized by a scheme $S$ of finite type over $\IC$. 
Let $\CE_S$, $\CE'_S$ be the corresponding flat families 
of stable ADHM sheaves. Then we have  canonical isomorphisms compatible with base 
change 
\be\label{eq:tangobsisomA}
\begin{aligned}
{\mathcal Ext}_{q_S}^k(C(Q_S,\rho_S), 
C(Q'_S,\rho'_S)\otimes_{Y_S}q_S^*I) & \simeq 
{\bf R}^kp_{S*} \wCC(\CE_S,\CE'_S,I)\\
\end{aligned}
\ee
for any $k\in \IZ$ and for any coherent 
$\CO_S$-module $I$, where $q_S: Y_S\to S$ 
is the projection morphism $q_S=p_S\circ \pi_S : Y_S \to S$
and $\wCC(\CE_S,\CE'_S,I)$ is the complex defined in \eqref{eq:enhancedB}.
\end{lemm}

{\it Proof.} 
We will use again the 
canonical quasi-isomorphism 
\[
(0,1_{\CO_{Y_S}}, \mathrm{ev}_S) : {\bf B}(Q_S,\rho_S)\to 
C(Q_S,\rho_S)
\] 
constructed in the proof of Lemma (\ref{relbeilinsonB}). This yields 
canonical isomorphisms 
\be\label{eq:relextisom}
\begin{aligned}
{\mathcal Ext}^k_{q_S}& 
(C(Q_S,\rho_S), C(Q'_S,\rho'_S)\otimes_{Y_S} q_S^*I) 
\simeq \\
& {\bf R}^kq_{S*}(C(Q'_S,\rho'_S)\otimes_{Y_S}{\bf B}(Q_S,\rho_S)^\vee 
\otimes_{Y_S} q_S^*I)\\
\end{aligned}
\ee
for all $k\in \IZ$. 
Next note that 
\[
C(Q'_S,\rho'_S)\otimes_{Y_S} q_S^*I = {\sf Cone}\bigg(
q_S^*I \xymatrix{{}\ar[rr]^-{-\rho_S'\otimes 1_{q_S^*I}} & & {}}
Q'_S\otimes_{Y_S} q_S^*I \bigg)[-1]
\]
Therefore 
\be\label{eq:tensorcone}
\begin{aligned}
& C(Q'_S,\rho_S') \otimes_{Y_S} 
{\bf B}(Q_S,\rho_S)^\vee\otimes_{Y_S}  q_S^*I =\\ 
& {\sf Cone}\left( \
\begin{CD}
{\overline {{\bf B}(Q_S,\rho_S)^\vee}}\otimes_{Y_S}q_S^*I \\
@VV{-1_{{\bf B}(Q_S,\rho_S)}
\otimes \rho_S'\otimes 1_{q_s^*I}}V \\
{{\bf B}(Q_S,\rho_S)^\vee} \otimes_{Y_S} Q'_S \otimes_{Y_S} 
q_S^*I\\
\end{CD}\qquad
\right)[-1],
\end{aligned}
\ee
where ${\overline {{\bf B}(Q_S,\rho_S)^\vee}}$ is the complex obtained 
by flipping the sign of all differentials in ${\bf B}(Q_S,\rho_S)^\vee$.
Here we are using the sign conventions of \cite[Ch. 1.3]{duality}
as stated at the end of the introduction.  

Using equation \eqref{eq:threetermA} and the projection formula
for the flat morphism $\pi_S: Y_S \to X_S$, it is 
straightforward to check that all terms in the complex \eqref{eq:tensorcone}
are acyclic with respect to pushforward by $\pi_S$. 
This implies that the 
Grothendieck spectral sequence for the composition $q_S = p_S \circ \pi_S$ 
collapses to an isomorphism 
\be\label{eq:Gspecseq}
\begin{aligned}
& {\bf R}^kq_{S*}(C(Q'_S,\rho'_S)\otimes_{Y_S}{\bf B}(Q_S,\rho_S)^\vee  
\otimes_{Y_S} q_S^*I) \simeq \\
& {\bf R}^kp_{S*} \pi_{S*}(C(Q'_S,\rho'_S)\otimes_{Y_S}
{\bf B}(Q_S,\rho_S)^\vee \otimes_{Y_S} q_S^*I)\\
\end{aligned}
\ee
for each $k\in \IZ$. Moreover, using \eqref{eq:tensorcone} 
we obtain 
\be\label{eq:tensorconeB}
\begin{aligned}
& \pi_{S*} ({\bf B}(Q_S,\rho_S)^\vee \otimes_{Y_S} 
C(Q'_S,\rho'_S)\otimes_{Y_S} q_S^*I)=\\ 
& {\sf Cone}
\left(\
\begin{CD}
\pi_{S*}{\overline {{\bf B}(Q_S,\rho_S)^\vee}}\otimes_{Y_S}q_S^*I\\
 @VV{\pi_{S*}(-1_{{\bf B}(Q_S,\rho_S)}
\otimes \rho_S'\otimes 1_{q_s^*I})}V\\
\pi_{S*}
{\bf B}(Q_S,\rho_S)^\vee \otimes_{Y_S} Q'_S
 \otimes_{Y_S} q_S^*I\\
\end{CD}\qquad \qquad\right)[-1].\\
\end{aligned}
\ee
Next, we claim  
there is a quasi-isomorphism of complexes of $\CO_{X_S}$-modules 
\be\label{eq:quasi-A}
\CO_{X_S}\oplus 
{\mathcal Hom}_{X_S}(E_S\otimes_{X_S}(M_{12})_S, p_S^*I)[-1] 
\to \pi_{S*}({\overline {{\bf B}(Q_S,\rho_S)^\vee}}\otimes_{Y_S} q_S^*I).
\ee
Since the complex ${\bf B}(Q_S,\rho_S)$ is locally free, this 
follows from a direct computation analogous to the proof of Lemma 
(\ref{relbeilinsonB}). Details will be omitted. 

The remaining direct image $\pi_{S*}({\bf B}(Q_S,\rho_S)^\vee\otimes_{Y_S} Q'_S\otimes_{Y_S}q_S^*I)$ in the right hand side of equation \eqref{eq:tensorconeB}
is also amenable to a direct explicit computation. Collecting all terms, 
it follows that there is indeed a canonical isomorphism of complexes 
of the form \eqref{eq:tangobsisomA}.

The base change property follows again from Lemma (\ref{relbasechange})
via the base change theorem and
general properties of direct images.

\hfill $\Box$ 

Entirely analogous computations also prove 

\begin{lemm}\label{tangobscompB}
Under the conditions of Lemma (\ref{tangobscompA})
we also have canonical isomorphisms compatible with base change 
\be\label{eq:tangobsisomK}
\begin{aligned}
{\mathcal Ext}_{q_S}^k(C(Q_S,\rho_S), 
C(Q'_S,\rho'_S)(-1)\otimes_{Y_S}q_S^*I) & \simeq 
{\bf R}^kp_{S*} \CC(\CE_S,\CE'_S,I)\\
\end{aligned}
\ee
for any $k\in \IZ$, 
where $\CC(\CE_S,\CE'_S,I)$ is the complex defined in \eqref{eq:hypercohB}.
\end{lemm}

\begin{lemm}\label{compcobdmaps}
Consider a deformation situation of the form \eqref{eq:defsit}. 
Let $(Q_S,\rho_S)$ be a flat family of admissible pairs 
parameterized by an affine scheme $S$ of finite type over $\IC$, 
and let $(Q_{S'}, \rho_{S'})$ be an 
extension to $S'$. Let $\CE_S, \CE_{S'}$ respectively be 
the corresponding flat families of stable ADHM sheaves
parameterized by $S,S'$. Then $\CE_{S'}$ is an extension 
of $\CE_S$ to $S'$, and we have commutative diagrams compatible with base change 
of the form 
\be\label{eq:cobddiagA}
\xymatrix{
{\mathrm{Ext}}^0_{Y_S}(C(Q_S,\rho_S), C(Q_{S},\rho_{S})) 
\ar[r]^{\simeq} 
\ar[d]_-{\delta^1} 
& \IH^0(X_S,  
\wCC(\CE_{S}, \CE_{S})) \ar[d]^-{{\widetilde \partial}^1}\\
{\mathrm{Ext}}^1_{Y_S}(C(Q_S,\rho_S), C(Q_{S'},\rho_{S'})\otimes_{Y_S}
I_{S\subset S'}) \ar[r]^-{\simeq} 
& \IH^0(X_S,  
\wCC(\CE_{S}, \CE_{S'}, I_{S\subset S'})),
\\}
\ee
\be\label{eq:cobddiagB} 
\xymatrix{ 
{\mathrm{Ext}}^1_{Y_S}(C(Q_S,\rho_S), C(Q_{S'},\rho_{S'})\otimes_{Y_S}
I_{S\subset S'}) \ar[r]^-{\simeq} 
\ar[d]_-{\delta^2} 
& \IH^1(X_S,
\wCC(\CE_{S}, \CE_{S'}, I_{S\subset S'})) \ar[d]^-{{\widetilde \partial}^2}\\
{\mathrm{Ext}}^2_{Y_S}(C(Q_S,\rho_S), C(Q_{S},\rho_{S})\otimes_{Y_S}
I_{S'\subset S''}) \ar[r]^-{\simeq} 
& \IH^2(X_S, 
\wCC(\CE_{S},I_{S'\subset S''})),
\\}
\ee
where the horizontal arrows are induced by the isomorphisms 
\eqref{eq:tangobsisomK}, and the vertical arrows are 
coboundary morphisms. 
\end{lemm}

{\it Proof.}
Using the isomorphisms \eqref{eq:relextisom},  
the coboundary morphisms $\delta^1$, $\delta^2$ are determined
by the following exact sequences of $\CO_{Y_S}$-modules 
\[
\begin{aligned}
0& \to C(Q_{S'}, \rho_{S'}) 
\otimes_{Y_S}{\bf B}(Q_S,\rho_S)^\vee
\otimes_{Y_S} I_{S\subset S'} \to 
C(Q_{S'}, \rho_{S'}) \otimes_{Y_S}
{\bf B}(Q_S,\rho_S)^\vee\\
& \to C(Q_S,\rho_S)\otimes_{Y_S} {\bf B}(Q_S,\rho_S) ^\vee\to 0.\\
0&\to C(Q_{S}, \rho_{S}) 
\otimes_{Y_S}{\bf B}(Q_S,\rho_S)^\vee\otimes_{Y_S} I_{S'\subset S''} \to 
C(Q_{S'}, \rho_{S'})
\otimes_{Y_S}{\bf B}(Q_S,\rho_S)^\vee\otimes_{Y_S} I_{S\subset S''}  \\ 
& \to C(Q_{S'}, \rho_{S'})
\otimes_{Y_S} {\bf B}(Q_S,\rho_S) ^\vee\otimes_{Y_S} I_{S\subset S'}\to 0.\\  
\end{aligned}
\]
All above complexes consist of acyclic terms with respect to 
push forward by $\pi_S:Y_S\to X_S$. Therefore applying $\pi_{S*}$ 
to the above exact sequences we obtain the exact sequences of 
$\CO_{X_S}$-modules 
\[ 
\begin{aligned} 
0 & \to \wCC(\CE_S,\CE_{S'},I_{S\subset S'}) \to 
\wCC(\CE_S, \CE_{S'}) \to \wCC(\CE_S) \to 0\\
0 & \to \wCC(\CE_S,I_{S'\subset S''}) \to 
\wCC(\CE_S, \CE_{S'}, I_{S\subset S''}) \to 
\wCC(\CE_S,\CE_{S'}, I_{S\subset S'}) \to 0\\
\end{aligned}
\]
Now the Lemma (\ref{compcobdmaps}) follows from the fact that there is 
an isomorphism of derived functors 
\[
{\bf R}\Gamma_{Y_S} \simeq {\bf R}\Gamma_{S}\circ {\bf R}q_{S*} 
\simeq {\bf R}\Gamma_S \circ {\bf R}p_{S*} \circ {\bf R}\pi_{S*} 
\]
Compatibility with base change follows 
again by analogy with Lemma (\ref{tangobscompA}). 

\hfill $\Box$.

\noindent
Finally, we also have 
\begin{lemm}\label{perfcomp}
Let $(Q_S,\rho_S)$ be a flat family of admissible pairs on $Y$ 
parameterized by an affine scheme $S$ of finite type over $\IC$
and let $\CE_S$ be the corresponding family of stable ADHM sheaves.
Then there exist two-term complexes of coherent locally free 
$\CO_S$-modules $\IF_S^\bullet=(F^1_S\to F^2_S)$, $\IE_S^\bullet
=(E_S^1\to E_S^2)$ such that 
\be\label{eq:perfcompA}
\begin{aligned}
\CH^1(\IF_S^\bullet\otimes_S I) & \simeq {\mathcal Ext}_{q_S}^1(C(Q_S,\rho_S), 
C(Q_S,\rho_S)(-1)\otimes_{Y_S}I)\\
\CH^2(\IF_S^\bullet\otimes_S I) & \simeq {\mathcal Ext}_{q_S}^2(C(Q_S,\rho_S), 
C(Q_S,\rho_S)(-1))\otimes_{S}I,\\
\end{aligned}
\ee
\be\label{eq:perfcompB} 
\begin{aligned}
\CH^1(\IE_S^\bullet\otimes_S I ) & \simeq {\bf R}^1p_{S*} \CC(\CE_S,I)\\
\CH^2(\IE_S^\bullet\otimes_S I ) & \simeq 
{\bf R}^2p_{S*}\CC(\CE_S)\otimes_S I\\
\end{aligned}
\ee
for any coherent $\CO_S$-module $I$, 
and a quasi-isomorphism ${\sf q}:\IF_S^\bullet \to 
\IE_S^\bullet$ such that the ${\sf q}\otimes 1_I$ induces the isomorphisms 
\eqref{eq:tangobsisomK} in cohomology, for $k=1,2$. 
\end{lemm}

{\it Proof.} 
Lemmas (\ref{relbeilinsonB}), (\ref{tangobscompB}) 
imply that we have an isomorphism 
\[
\begin{aligned}
& {\bf R}q_{S*}{\bf R}{\mathcal Hom}_{Y_S}
(C(Q_S,\rho_S), C(Q_S,\rho_S)(-1)) \simeq \\
& {\bf R}q_{S*}(B(Q_S,\rho_S)\otimes_{Y_S} B(Q_S,\rho_S)^\vee\otimes_{Y_S}I)
\simeq 
{\bf R}p_{S*} \CC(\CE_S,I) \\
\end{aligned}
\]
in the derived category $D^-(S)$, for any coherent $\CO_S$-module $I$. 
In particular, there exist complexes 
of $\CO_S$-modules ${\widetilde \IF}_S^\bullet$, ${\widetilde \IE}_S^\bullet$
bounded above such that  
\[
\begin{aligned} 
\CH^k({\widetilde \IF}_S^\bullet\otimes_S I) & 
\simeq {\mathcal Ext}_{q_S}^k(C(Q_S,\rho_S), 
C(Q_S,\rho_S)(-1)\otimes_{Y_S}I)\\
\CH^k(\IE_S^\bullet\otimes_{S}I) & \simeq {\bf R}^kp_{S*} \CC(\CE_S,I)\\
\end{aligned}
\]
for $k=1,2$, and a quasi-isomorphism 
${\widetilde {\sf q}}:{\widetilde \IF}_S^\bullet 
\to {\widetilde \IE}_S^\bullet$ which induces 
the isomorphisms \eqref{eq:tangobsisomK}. 
Then, employing the inductive 
construction of \cite[Lemma III.12.3]{har}, we obtain complexes 
${\IF'}_S^\bullet$, ${\IE'}_S^\bullet$ of coherent locally 
free $\CO_S$-modules, bounded above,
and a natural quasi-isomorphism
${\sf q}':{\IF'}_S^\bullet\to {\IE'}_S^\bullet$
with the same properties. Moreover, using the same 
argument as in the proof of Lemma (\ref{perfect}),  
the complexes ${\IF'}_S^\bullet$, 
${\IE'}_S^\bullet$ can be 
simultaneously truncated to two-term complexes of coherent 
locally-free modules  
${\IF}_S$, $\IE_S$. The quasi-isomorphism ${\sf q}'$ also yields 
by truncation a quasi-isomorphism
${\sf q}: \IF_S \to \IE_S$. 

\hfill $\Box$ 

Now we can conclude the proof of Theorem (\ref{correspB}). 
The tangent-obstruction theory of the moduli space of ADHM sheaves 
has been constructed in (\ref{tangspace}), (\ref{obsspace}), 
the relevant deformation theory results being proven in 
Proposition (\ref{infdef}), Corollary (\ref{infdefB})
and Proposition (\ref{infobs}). The
tangent-obstruction theory of the moduli space 
of admissible pairs has been presented in the proof of Theorem 
(\ref{admstack}) (see section (\ref{admpairs}).) 
Then statement (I) in section 
(\ref{compattangobs}) follows from Lemmas 
(\ref{tracelemma}), 
(\ref{tangobscompA}), (\ref{tangobscompB}), (\ref{compcobdmaps}). 
 Statements $(II.a)-(II.c)$ follow from Lemmas (\ref{tangobscompB})
and (\ref{perfcomp}). 

\hfill $\Box$

\bibliography{adhmref.bib}

\begin{thebibliography}{10}

\bibitem{dimred}
L.~{\'A}lvarez-C{\'o}nsul and O.~Garc{\'{\i}}a-Prada.
\newblock Dimensional reduction and quiver bundles.
\newblock {\em J. Reine Angew. Math.}, 556:1--46, 2003.

\bibitem{HKquivers}
L.~{\'A}lvarez-C{\'o}nsul and O.~Garc{\'{\i}}a-Prada.
\newblock Hitchin-{K}obayashi correspondence, quivers, and vortices.
\newblock {\em Comm. Math. Phys.}, 238(1-2):1--33, 2003.

\bibitem{holchains}
L.~{\'A}lvarez-C{\'o}nsul, O.~Garc{\'{\i}}a-Prada, and A.~H.~W. Schmitt.
\newblock On the geometry of moduli spaces of holomorphic chains over compact
  {R}iemann surfaces.
\newblock {\em IMRP Int. Math. Res. Pap.}, pages Art. ID 73597, 82, 2006.

\bibitem{artin-moduli}
M.~Artin.
\newblock Algebraization of formal moduli {I}.
\newblock In {\em Global Analysis}, pages 21--71. Tokyo, 1969.

\bibitem{artin-versal}
M.~Artin.
\newblock Versal deformations and algebraic stacks.
\newblock {\em Invent. Math.}, 27:165--189, 1974.

\bibitem{Banica}
C.~B{\u{a}}nic{\u{a}} and O.~St{\u{a}}n{\u{a}}{\c{s}}il{\u{a}}.
\newblock {\em Algebraic methods in the global theory of complex spaces}.
\newblock Editura Academiei, Bucharest, 1976.
\newblock Translated from the Romanian.

\bibitem{BF}
K.~Behrend and B.~Fantechi.
\newblock The intrinsic normal cone.
\newblock {\em Invent. Math.}, 128(1):45--88, 1997.

\bibitem{benzvi-2006}
D.~Ben-Zvi and T.~Nevins.
\newblock Perverse bundles and {C}alogero-{M}oser spaces.
\newblock {\em Compos. Math.}, 144(6):1403--1428, 2008.

\bibitem{logflips}
A.~Bertram.
\newblock Stable pairs and log flips.
\newblock In {\em Algebraic geometry---Santa Cruz 1995}, volume~62 of {\em
  Proc. Sympos. Pure Math.}, pages 185--201. Amer. Math. Soc., Providence, RI,
  1997.

\bibitem{BCK-I}
A.~Bertram, I.~Ciocan-Fontanine, and B.~Kim.
\newblock Two proofs of a conjecture of {H}ori and {V}afa.
\newblock {\em Duke Math. J.}, 126(1):101--136, 2005.

\bibitem{BCK-II}
A.~Bertram, I.~Ciocan-Fontanine, and B.~Kim.
\newblock Gromov-{W}itten invariants for abelian and nonabelian quotients.
\newblock {\em J. Algebraic Geom.}, 17(2):275--294, 2008.

\bibitem{BDW}
A.~Bertram, G.~Daskalopoulos, and R.~Wentworth.
\newblock Gromov invariants for holomorphic maps from {R}iemann surfaces to
  {G}rassmannians.
\newblock {\em J. Amer. Math. Soc.}, 9(2):529--571, 1996.

\bibitem{pairsI}
S.~B. Bradlow and G.~D. Daskalopoulos.
\newblock Moduli of stable pairs for holomorphic bundles over {R}iemann
  surfaces.
\newblock {\em Internat. J. Math.}, 2(5):477--513, 1991.

\bibitem{surfhiggs}
S.~B. Bradlow, O.~Garc{\'{\i}}a-Prada, and P.~B. Gothen.
\newblock Surface group representations and {${\rm U}(p,q)$}-{H}iggs bundles.
\newblock {\em J. Differential Geom.}, 64(1):111--170, 2003.

\bibitem{triples}
S.~B. Bradlow, O.~Garc{\'{\i}}a-Prada, and P.~B. Gothen.
\newblock Moduli spaces of holomorphic triples over compact {R}iemann surfaces.
\newblock {\em Math. Ann.}, 328(1-2):299--351, 2004.

\bibitem{BP}
J.~Bryan and R.~Pandharipande.
\newblock The local {G}romov-{W}itten theory of curves.
\newblock {\em J. Amer. Math. Soc.}, 21(1):101--136 (electronic), 2008.
\newblock With an appendix by Bryan, C. Faber, A. Okounkov and Pandharipande.

\bibitem{hyperquot}
K.~L. C.-H.~Liu and S.-T. Yau.
\newblock ${S}^1$-fixed-points in hyper-{Q}uot-schemes and an exact mirror
  formula for flag manifolds from the extended mirror principle diagram.
\newblock arxiv:math/0401367.

\bibitem{nonabhilb}
I.~Ciocan-Fontanine, D.-E. Diaconescu, B.~Kim, and D.~Maulik.
\newblock work in progress.

\bibitem{vf-dg}
I.~Ciocan-Fontanine and M.~Kapranov.
\newblock Virtual fundamental classes via dg-manifolds.
\newblock {\em Geom. Topol.}, 13(3):1779--1804, 2009.

\bibitem{Cirafici:2008sn}
M.~Cirafici, A.~Sinkovics, and R.~J. Szabo.
\newblock {Cohomological gauge theory, quiver matrix models and
  Donaldson-Thomas theory}.
\newblock {\em Nucl. Phys.}, B809:452--518, 2009.

\bibitem{duality}
B.~Conrad.
\newblock {\em Grothendieck duality and base change}, volume 1750 of {\em
  Lecture Notes in Mathematics}.
\newblock Springer-Verlag, Berlin, 2000.

\bibitem{pseudored}
B.~Conrad, O.~Gabber, and G.~Prasad.
\newblock {\em Pseudo-reductive groups}, volume~17 of {\em New Mathematical
  Monographs}.
\newblock Cambridge University Press, Cambridge, 2010.

\bibitem{Demazure}
M.~Demazure.
\newblock Structures alg{\'e}briques, cohomologie des groupes.
\newblock In {\em SGA 3: Sch{\'e}mas en groupes}, volume 151 of {\em Lecture
  Notes in Mathematics}, pages 1--37. Springer Verlag, Berlin, 1970.

\bibitem{chamberI}
D.-E. Diaconescu.
\newblock Chamber structure and wallcrossing in the {A}{D}{H}{M} theory of
  curves {I}.
\newblock arXiv:0904.4451.

\bibitem{freefermions}
R.~Dijkgraaf, L.~Hollands, P.~Sulkowski, and C.~Vafa.
\newblock {Supersymmetric Gauge Theories, Intersecting Branes and Free
  Fermions}.
\newblock {\em JHEP}, 02:106, 2008.

\bibitem{Donagi:2007hi}
R.~Donagi and E.~Sharpe.
\newblock {GLSM's for partial flag manifolds}.
\newblock {\em J. Geom. Phys.}, 58:1662--1692, 2008.

\bibitem{DM-quivers}
M.~Douglas and G.~Moore.
\newblock D-branes, quivers, and {A}{L}{E} instantons.
\newblock http://arxiv.org/abs/hep-th/9603167.

\bibitem{eisenbud}
D.~Eisenbud.
\newblock {\em Commutative algebra with a view toward algebraic geometry},
  volume 150 of {\em Graduate Texts in Mathematics}.
\newblock Springer-Verlag, New York, 1995.

\bibitem{FG-virtual}
B.~Fantechi and L.~G{\"o}ttsche.
\newblock Riemann-{R}och theorems and elliptic genus for virtually smooth
  schemes.
\newblock {\em Geom. Topol.}, 14(1):83--115, 2010.

\bibitem{tensors}
T.~L. Gomez and I.~Sols.
\newblock Stable tensors and moduli space of orthogonal sheaves.
\newblock math.AG/0103150.

\bibitem{homquiv}
P.~B. Gothen and A.~D. King.
\newblock Homological algebra of twisted quiver bundles.
\newblock {\em J. London Math. Soc. (2)}, 71(1):85--99, 2005.

\bibitem{GP}
T.~Graber and R.~Pandharipande.
\newblock Localization of virtual classes.
\newblock {\em Invent. Math.}, 135(2):487--518, 1999.

\bibitem{EGAI}
A.~Grothendieck.
\newblock \'{E}l\'ements de g\'eom\'etrie alg\'ebrique. {I}. {L}e langage des
  sch\'emas.
\newblock {\em Inst. Hautes \'Etudes Sci. Publ. Math.}, (4):228, 1960.

\bibitem{EGAIIIa}
A.~Grothendieck.
\newblock \'{E}l\'ements de g\'eom\'etrie alg\'ebrique. {III}. \'{E}tude
  cohomologique des faisceaux coh\'erents. {I}.
\newblock {\em Inst. Hautes \'Etudes Sci. Publ. Math.}, (11):167, 1961.

\bibitem{EGAIIIb}
A.~Grothendieck.
\newblock \'{E}l\'ements de g\'eom\'etrie alg\'ebrique. {III}. \'{E}tude
  cohomologique des faisceaux coh\'erents. {II}.
\newblock {\em Inst. Hautes \'Etudes Sci. Publ. Math.}, (17):91, 1963.

\bibitem{resdual}
R.~Hartshorne.
\newblock {\em Residues and duality}.
\newblock Lecture notes of a seminar on the work of A. Grothendieck, given at
  Harvard 1963/64. With an appendix by P. Deligne. Lecture Notes in
  Mathematics, No. 20. Springer-Verlag, Berlin, 1966.

\bibitem{har}
R.~Hartshorne.
\newblock {\em Algebraic geometry}.
\newblock Springer-Verlag, New York, 1977.
\newblock Graduate Texts in Mathematics, No. 52.

\bibitem{HT}
K.~Hori and D.~Tong.
\newblock {Aspects of non-Abelian gauge dynamics in two-dimensional N = (2,2)
  theories}.
\newblock {\em JHEP}, 05:079, 2007.

\bibitem{HV}
K.~Hori and C.~Vafa.
\newblock {Mirror symmetry}.
\newblock hep-th/0002222.

\bibitem{framed}
D.~Huybrechts and M.~Lehn.
\newblock Framed modules and their moduli.
\newblock {\em Internat. J. Math.}, 6(2):297--324, 1995.

\bibitem{stpairs}
D.~Huybrechts and M.~Lehn.
\newblock Stable pairs on curves and surfaces.
\newblock {\em J. Algebraic Geom.}, 4(1):67--104, 1995.

\bibitem{huylehn}
D.~Huybrechts and M.~Lehn.
\newblock {\em The geometry of moduli spaces of sheaves}.
\newblock Cambridge Mathematical Library. Cambridge University Press,
  Cambridge, second edition, 2010.

\bibitem{INOV}
A.~Iqbal, N.~Nekrasov, A.~Okounkov, and C.~Vafa.
\newblock {Quantum foam and topological strings}.
\newblock {\em JHEP}, 04:011, 2008.

\bibitem{instsheaves}
M.~Jardim.
\newblock Moduli spaces of framed instanton sheaves on projective spaces.
\newblock arXiv:0801.2550.

\bibitem{Kapustin:2000ek}
A.~Kapustin, A.~Kuznetsov, and D.~Orlov.
\newblock {Noncommutative instantons and twistor transform}.
\newblock {\em Commun. Math. Phys.}, 221:385--432, 2001.

\bibitem{Keel-Mori}
S.~Keel and S.~Mori.
\newblock Quotients by groupoids.
\newblock {\em Ann. of Math. (2)}, 145(1):193--213, 1997.

\bibitem{BKP}
B.~Kim, A.~Kresch, and T.~Pantev.
\newblock Functoriality in intersection theory and a conjecture of {C}ox,
  {K}atz, and {L}ee.
\newblock {\em J. Pure Appl. Algebra}, 179(1-2):127--136, 2003.

\bibitem{ALEquiver-I}
P.~B. Kronheimer and H.~Nakajima.
\newblock Yang-{M}ills instantons on {ALE} gravitational instantons.
\newblock {\em Math. Ann.}, 288(2):263--307, 1990.

\bibitem{LePoitierA}
J.~Le~Poitier.
\newblock Syst{\`e}mes coh{\'e}rents et structures de niveau.
\newblock {\em Ast{\'e}risque}, 214:143, 1993.

\bibitem{LePoitierB}
J.~Le~Poitier.
\newblock Faisceaux semi-stable et syst{\`e}mes coh{\'e}rents.
\newblock In {\em Vector Bundles in Algebraic Geometry (Durham, 1993)}, volume
  208 of {\em London Math. Soc. Lecture Notes Ser.}, pages 179--239. Cambridge
  University Press, Cambridge, 1995.

\bibitem{degGW}
J.~Li.
\newblock A degeneration formula of {GW}-invariants.
\newblock {\em J. Differential Geom.}, 60(2):199--293, 2002.

\bibitem{LT}
J.~Li and G.~Tian.
\newblock Virtual moduli cycles and {G}romov-{W}itten invariants of algebraic
  varieties.
\newblock {\em J. Amer. Math. Soc.}, 11(1):119--174, 1998.

\bibitem{fixed-quot}
B.~H. Lian, C.-H. Liu, K.~Liu, and S.-T. Yau.
\newblock The {$S\sp 1$} fixed points in {Q}uot-schemes and mirror principle
  computations.
\newblock In {\em Vector bundles and representation theory (Columbia, MO,
  2002)}, volume 322 of {\em Contemp. Math.}, pages 165--194. Amer. Math. Soc.,
  Providence, RI, 2003.

\bibitem{MNOP-I}
D.~Maulik, N.~Nekrasov, A.~Okounkov, and R.~Pandharipande.
\newblock Gromov-{W}itten theory and {D}onaldson-{T}homas theory. {I}.
\newblock {\em Compos. Math.}, 142(5):1263--1285, 2006.

\bibitem{MNOP-II}
D.~Maulik, N.~Nekrasov, A.~Okounkov, and R.~Pandharipande.
\newblock Gromov-{W}itten theory and {D}onaldson-{T}homas theory. {II}.
\newblock {\em Compos. Math.}, 142(5):1286--1304, 2006.

\bibitem{Minv}
T.~Mochizuki.
\newblock A theory of the invariants obtained from the moduli stacks of stable
  objects on a smooth polarized surface.
\newblock arXiv:math/0210211.

\bibitem{MP}
D.~R. Morrison and M.~Ronen~Plesser.
\newblock {Summing the instantons: Quantum cohomology and mirror symmetry in
  toric varieties}.
\newblock {\em Nucl. Phys.}, B440:279--354, 1995.

\bibitem{GIT}
D.~Mumford, J.~Fogarty, and F.~Kirwan.
\newblock {\em Geometric invariant theory}, volume~34 of {\em Ergebnisse der
  Mathematik und ihrer Grenzgebiete (2) [Results in Mathematics and Related
  Areas (2)]}.
\newblock Springer-Verlag, Berlin, third edition, 1994.

\bibitem{ALEquiver-II}
H.~Nakajima.
\newblock Instantons on {ALE} spaces, quiver varieties, and {K}ac-{M}oody
  algebras.
\newblock {\em Duke Math. J.}, 76(2):365--416, 1994.

\bibitem{hilblect}
H.~Nakajima.
\newblock {\em Lectures on {H}ilbert schemes of points on surfaces}, volume~18
  of {\em University Lecture Series}.
\newblock American Mathematical Society, Providence, RI, 1999.

\bibitem{instcountA}
H.~Nakajima and K.~Yoshioka.
\newblock Instanton counting on blowup. {I}. 4-dimensional pure gauge theory.
\newblock {\em Invent. Math.}, 162(2):313--355, 2005.

\bibitem{instcountB}
H.~Nakajima and K.~Yoshioka.
\newblock Instanton counting on blowup. {II}. {$K$}-theoretic partition
  function.
\newblock {\em Transform. Groups}, 10(3-4):489--519, 2005.

\bibitem{Nekrasov:1998ss}
N.~Nekrasov and A.~S. Schwarz.
\newblock {Instantons on noncommutative ${R}^4$ and (2,0) superconformal six
  dimensional theory}.
\newblock {\em Commun. Math. Phys.}, 198:689--703, 1998.

\bibitem{Nekrasov:2009ui}
N.~A. Nekrasov and S.~L. Shatashvili.
\newblock {Quantum integrability and supersymmetric vacua}.
\newblock {\em Prog. Theor. Phys. Suppl.}, 177:105--119, 2009.

\bibitem{Nekrasov:2009uh}
N.~A. Nekrasov and S.~L. Shatashvili.
\newblock {Supersymmetric vacua and Bethe ansatz}.
\newblock {\em Nucl. Phys. Proc. Suppl.}, 192-193:91--112, 2009.

\bibitem{OT}
C.~Okonek and A.~Teleman.
\newblock Comparing virtual fundamental classes: gauge theoretical
  {G}romov-{W}itten invariants for toric varieties.
\newblock {\em Asian J. Math.}, 7(2):167--198, 2003.

\bibitem{OP}
A.~Okounkov and R.~Pandharipande.
\newblock The local {D}onaldson-{T}homas theory of curves.
\newblock {\em Geom. Topol.}, (14):1503--1567, 2010.

\bibitem{orlov}
D.~O. Orlov.
\newblock Projective bundles, monoidal transformations, and derived categories
  of coherent sheaves.
\newblock {\em Izv. Ross. Akad. Nauk Ser. Mat.}, 56(4):852--862, 1992.

\bibitem{stabpairs-II}
R.~Pandharipande and R.~P. Thomas.
\newblock The 3-fold vertex via stable pairs.
\newblock {\em Geom. Topol.}, 13(4):1835--1876, 2009.

\bibitem{stabpairs-I}
R.~Pandharipande and R.~P. Thomas.
\newblock Curve counting via stable pairs in the derived category.
\newblock {\em Invent. Math.}, 178(2):407--447, 2009.

\bibitem{stabpairs-III}
R.~Pandharipande and R.~P. Thomas.
\newblock Stable pairs and {BPS} invariants.
\newblock {\em J. Amer. Math. Soc.}, 23(1):267--297, 2010.

\bibitem{gract-stacks}
M.~Romagny.
\newblock Group actions on stacks and applications.
\newblock {\em Michigan Math. J.}, 53(1):209--236, 2005.

\bibitem{decorated}
A.~Schmitt.
\newblock A universal construction for moduli spaces of decorated vector
  bundles over curves.
\newblock {\em Transform. Groups}, 9(2):167--209, 2004.

\bibitem{modquivers}
A.~Schmitt.
\newblock Moduli for decorated tuples of sheaves and representation spaces for
  quivers.
\newblock {\em Proc. Indian Acad. Sci. Math. Sci.}, 115(1):15--49, 2005.

\bibitem{frHitchin}
A.~H.~W. Schmitt.
\newblock Framed {H}itchin pairs.
\newblock {\em Rev. Roumaine Math. Pures Appl.}, 45(4):681--711 (2001), 2000.

\bibitem{GIT-decorated}
A.~H.~W. Schmitt.
\newblock {\em Geometric invariant theory and decorated principal bundles}.
\newblock Zurich Lectures in Advanced Mathematics. European Mathematical
  Society (EMS), Z\"urich, 2008.

\bibitem{equiv-compl}
H.~Sumihiro.
\newblock Equivariant completion.
\newblock {\em J. Math. Kyoto Univ.}, 14:1--28, 1974.

\bibitem{szendroi-2005}
B.~Szendr{\H{o}}i.
\newblock Sheaves on fibered threefolds and quiver sheaves.
\newblock {\em Comm. Math. Phys.}, 278(3):627--641, 2008.

\bibitem{verlinde}
M.~Thaddeus.
\newblock Stable pairs, linear systems and the {V}erlinde formula.
\newblock {\em Invent. Math.}, 117(2):317--353, 1994.

\bibitem{RT}
R.~P. Thomas.
\newblock A holomorphic {C}asson invariant for {C}alabi-{Y}au 3-folds, and
  bundles on {$K3$} fibrations.
\newblock {\em J. Differential Geom.}, 54(2):367--438, 2000.

\bibitem{equivres}
R.~W. Thomason.
\newblock Equivariant resolution, linearization, and {H}ilbert's fourteenth
  problem over arbitrary base schemes.
\newblock {\em Adv. in Math.}, 65(1):16--34, 1987.

\end{thebibliography}

\bibliographystyle{abbrv}
\end{document}